\documentclass[preprint-12pt,times,sort&compress]{elsarticle}
\usepackage[letterpaper,top=2.5cm,bottom=2.5cm,left=2.5cm,right=2.5cm]{geometry}

\pdfoutput=1
\journal{}
\usepackage{hyperref}
\usepackage{amsmath}
\usepackage{amsbsy}
\usepackage{mathtools}
\usepackage{makecell}
\usepackage{booktabs}
\usepackage{multirow}
\usepackage[table,xcdraw,dvipsnames]{xcolor}
\usepackage{graphicx}
\usepackage{subcaption}
\usepackage{float}
\usepackage{overpic} 
\usepackage{rotating}
\usepackage[export]{adjustbox}
\usepackage{tikz}
\usepackage{pgfplots}
\pgfplotsset{compat=newest}

\DeclareSymbolFont{rmlargesymbols}{OMX}{cmex}{m}{n}
\DeclareMathSymbol{\rmsum}{\mathop}{rmlargesymbols}{80}
\DeclareMathSymbol{\rmintop}{\mathop}{rmlargesymbols}{82}
\renewcommand{\sum}{\rmsum}
\renewcommand{\int}{\rmintop\nolimits}

\DeclareSymbolFont{greeksymbolsptm}{OML}{ptm}{m}{n}
\DeclareMathSymbol{\sigma}{\mathalpha}{greeksymbolsptm}{27}

\DeclareSymbolFont{greeksymbolsit}{OML}{mdbch}{i}{n}
\DeclareMathSymbol{\lambda}{\mathalpha}{greeksymbolsit}{21}

\DeclareSymbolFont{greeksymbols}{OML}{mdbch}{m}{n}
\DeclareMathSymbol{\Omega}{\mathalpha}{greeksymbols}{10}

\DeclareSymbolFont{greeksymbolsant}{OML}{antt}{i}{n}
\DeclareMathSymbol{\nu}{\mathalpha}{greeksymbolsant}{23}

\DeclareMathAlphabet{\altmathcal}{OMS}{cmsy}{m}{n}
\newcommand{\kr}{\altmathcal{K}}		
\newcommand{\Nr}{\altmathcal{N}}		
\renewcommand{\O}{\altmathcal{O}}		
\newcommand{\Or}[1]{\O\pra{#1}}			

\newcommand{\K}{\textbf{K}}				
\newcommand{\M}{\textbf{M}}				
\newcommand{\LL}{\textbf{L}}			
\renewcommand{\H}{\textbf{H}}			
\newcommand{\D}{\rm{D}}					
\newcommand{\DD}{\textbf{D}}			
\newcommand{\T}{\textbf{T}}				
\newcommand{\V}{\textbf{V}}				
\newcommand{\uu}{\textbf{u}}			
\newcommand{\uuh}{\uu^{\.h}}			
\newcommand{\uh}{u^{\.h}}				
\newcommand{\ue}{u^{\.e}}				
\newcommand{\w}{\textbf{w}}				
\newcommand{\vv}{\textbf{v}}			
\newcommand{\e}{\textbf{e}}				
\newcommand{\rr}{\textbf{r}}			
\newcommand{\R}{\text{I\!R}}			
\newcommand{\om}{\omega}				
\newcommand{\Om}{\Omega}				
\newcommand{\sg}{\sigma}				
\newcommand{\lm}{\lambda}				
\newcommand{\lmh}{\lm^h}				
\newcommand{\lme}{\lm^e}				
\newcommand{\EV}{ \text{EVerr}}			
\newcommand{\EFL}{\text{EFerr}_{L^2}}	
\newcommand{\EFE}{\text{EFerr}_{E}}		
\newcommand{\ti}{T}						
\newcommand{\tiN}{\ti_{\N_0}}			
\newcommand{\tih}{\widetilde{\ti}}		
\newcommand{\tia}{\ti_{\rm av}}			
\newcommand{\tiax}[1]{\ti_{\rm av,#1}}	

\newcommand{\Uu}{U}						
\newcommand{\UU}{\textbf{\Uu}}			
\newcommand{\n}{n} 						
\newcommand{\p}{p} 						
\newcommand{\B}[2]{B_{#1,#2}}			
\newcommand{\x}{x} 						
\newcommand{\y}{y}						
\newcommand{\z}{z}						
\newcommand{\xx}{\textbf{x}} 			
\newcommand{\X}{\Xi} 					
\renewcommand{\ne}{n_e} 				
\newcommand{\N}{N} 				 		
\newcommand{\Nit}{\Nr_{\text{iter}}}	
\newcommand{\Nsh}{\Nr_{\text{shift}}}	
\newcommand{\Nfa}{\Nr_{\text{fact}}}	
\newcommand{\Nfb}{\Nr_{\fb}}			
\newcommand{\Nmv}{\Nr_{\text{m--v}}}	
\newcommand{\nz}{N_{nz}}				
\renewcommand{\l}{\ell} 				
\newcommand{\BS}{\textrm{blocksize}}	
\newcommand{\fb}{\textrm{f/\.b}}		
\newcommand{\Fb}{\textrm{F/\.b}}		
\newcommand{\FL}{\textrm{FLOPs}}		

\renewcommand{\.}{\!\:}					
\newcommand{\pr}{^{\,\prime}}			
\newcommand{\pra}[1]{\left(#1\right)}	
\newcommand{\bra}[1]{\left[#1\right]}	
\newcommand{\norm}[1]{\left\lVert#1\right\rVert} 
\newcommand{\normM}[1]{\norm{#1}_{\.\M}}
\newcommand{\normL}[1]{\norm{#1}_{\.L^2}}
\newcommand{\normE}[1]{\norm{#1}_{\.E}}

\newcommand{\az}[1]{{\color{black}#1}}
\newcommand{\fig}[1]{\mbox{\az{Fig.} #1}}
\newcommand{\figs}[2]{\mbox{\az{Figs.} #1 and #2}}

\newcommand{\eq}[1]{\mbox{Eq. #1}}

\newcommand{\eqz}[2]{\mbox{Eqs. #1\.--\.#2}}
\newcommand{\tab}[1]{\mbox{\az{Table} #1}}

\newcommand{\Sec}[1]{\mbox{Section #1}}

\newcommand{\Rem}[1]{\mbox{Remark #1}}
\newcommand{\Lem}[1]{\mbox{Lemma #1}}

\newdefinition{rem}{Remark}
\newdefinition{lem}{Lemma}

\definecolor{drot}{rgb}{0.7,0,0.1}

\definecolor{C1}{rgb}{   0,    0.45,    0.74}
\definecolor{C2}{rgb}{0.85,    0.33,    0.10}
\definecolor{C3}{rgb}{0.93,    0.69,    0.13}
\definecolor{C4}{rgb}{0.49,    0.18,    0.56}
\definecolor{C5}{rgb}{0.47,    0.67,    0.19}
\definecolor{C6}{rgb}{0.30,    0.75,    0.93}
\definecolor{C7}{rgb}{0.64,    0.08,    0.18}
\definecolor{C8}{rgb}{1.00,       0,    1.00}
\definecolor{C9}{rgb}{   0,       0,       0}
\definecolor{C10}{rgb}{  0,       0,    1.00}

\pgfplotscreateplotcyclelist{mycolorlist}{%
C1, mark=none\\%
C2, mark=none\\%
C4, mark=none\\%
C5, mark=none\\%
C6, mark=none\\%
C8, mark=none\\%
C9, mark=none\\%
C10, mark=none\\%
}

\begin{document}
\baselineskip14pt
\sloppy

\begin{frontmatter}

\title{Refined isogeometric analysis for generalized Hermitian eigenproblems}

\author[a1]{Ali Hashemian\corref{cor1}}
\cortext[cor1]{Corresponding author}
\ead{ahashemian@bcamath.org}

\author[a2,a1,a3]{David Pardo}
\author[a4,a5,a6]{Victor M. Calo} 

\address[a1]{BCAM -- Basque Center for Applied Mathematics, Bilbao, Basque Country, Spain}
\address[a2]{University of the Basque Country UPV/EHU, Bilbao, Basque Country, Spain}
\address[a3]{Ikerbasque -- Basque Foundation for Sciences, Bilbao, Basque Country, Spain}
\address[a4]{Western Australian School of Mines, Curtin University, Perth, Australia}
\address[a5]{Commonwealth Scientific and Industrial Research Organisation (CSIRO), Perth, Australia}
\address[a6]{Curtin Institute for Computation, Curtin University, Perth, Australia}

\begin{abstract} 

  We use the refined isogeometric analysis (rIGA) to solve generalized Hermitian eigenproblems $({\K\uu=\lm\M\uu})$. The rIGA framework conserves the desirable properties of maximum-continuity isogeometric analysis (IGA) discretizations while reducing the computation cost of the solution through partitioning the computational domain by adding zero-continuity basis functions.  As a result, rIGA enriches the approximation space and decreases the interconnection between degrees of freedom.  We compare computational costs of rIGA versus those of IGA when employing a Lanczos eigensolver with a shift-and-invert spectral transformation.  When all eigenpairs within a given interval ${[\lm_s,\lm_e]}$ are of interest, we select several shifts ${\sg_k\in[\lm_s,\lm_e]}$ using a spectrum slicing technique.  For each shift $\sg_k$, the cost of factorization of the spectral transformation matrix ${\K-\sg_k \M}$ drives the total computational cost of the eigensolution.   Several multiplications of the operator matrices ${(\K-\sg_k\M)^{-1}\M}$ by vectors follow this factorization.  Let $\p$ be the polynomial degree of basis functions and assume that IGA has maximum continuity of ${p-1}$, while rIGA introduces $C^{\.0}$ separators to minimize the factorization cost. For this setup, our theoretical estimates predict computational savings to compute a fixed number of eigenpairs of up to ${\O(\p^{\.2})}$ in the asymptotic regime, that is, large problem sizes. Yet, our numerical tests show that for moderately-sized eigenproblems, the total computational cost reduction is $\O(\p)$. Nevertheless, rIGA improves the accuracy of every eigenpair of the first $\N_0$ eigenvalues and eigenfunctions. Here, we allow $\N_0$ to be as large as the total number of eigenmodes of the original maximum-continuity IGA discretization.

\end{abstract}


\end{frontmatter}


\section{Introduction} 
\label{sec.Introduction}

Hughes et al.~\cite{Hughes2005} introduced isogeometric analysis (IGA), a widely-used numerical technique for the solution approximation of partial differential equations (PDEs). IGA delivered useful solutions to many scientific and engineering problems  (see, e.g.,~\cite{ Bazilevs2006, GOMEZ2008, BUFFA2010, AURICCHIO2012, Espath2016, CASQUERO2017, SIMPSON2018, Puzyrev2019}).
IGA uses spline basis functions, which are standard in computer-aided design (CAD), as shape functions of finite element analysis (FEA). These functions can have high continuity (up to ${C^{\p\.-1}}$ where $\p$ is the polynomial order of spline bases) across the element interfaces.

Compared to the minimal interconnection of elements in traditional FEA, maximum-continuity IGA discretizations strengthen the interconnection between elements. This increased interconnectivity degrades the performance of sparse direct solvers (see, e.g.,~\cite{ Collier2012}).  Garcia et al.~\cite{ Garcia2017} introduced the refined isogeometric analysis (rIGA) to ameliorate this performance degradation and to exploit the recursive partitioning capability of multifrontal direct solvers~\cite{ Duff1983}.  The rIGA framework preserves some of the desirable properties of maximum-continuity IGA discretizations while partitioning the computational domain into macroelement blocks that are weakly interconnected by ${C^{\.0}}$ separators.  As a result, the matrix factorization (e.g., LU or Cholesky) step has a lower computational cost.  The performance of the rIGA framework on preconditioned conjugate gradient solvers as well as its applicability to mechanical and electromagnetic problems are also studied in~\cite{ Garcia2018, Garcia2019}.

The application of IGA in eigenanalysis is a well-studied topic in the literature (see, e.g.,~\cite{ Cottrell2006, HUGHES2008, Hughes2014, Mazza2019}).  Improving the efficiency of the system integration and accuracy of the spectral approximation of the IGA approach are also of interest (see, e.g.,~\cite{ Puzyrev2017, Deng2018, Hosseini2018, Deng20182, Calo2019, Deng2019, Deng20192, Barton2020}). 
Herein, we investigate the beneficial effect of using refined isogeometric analysis in eigenproblems. We compare the computational cost and accuracy of the resulting eigenpairs that both refined and maximum-continuity IGA produce.  We first review some numerical aspects of eigenanalysis to perform a detailed comparison of the methods and their results.

Eigenanalysis is a computationally expensive proposition, especially when seeking for a large number of eigenpairs (i.e., eigenvalues and eigenvectors) on a multidimensional domain.  Frequently-used Krylov eigensolvers such as Lanczos and Arnoldi project onto  Krylov subspaces.  The convergence rate of these iterative algorithms degrades when computing a large interval of eigenvalues, particularly when the eigenmodes are not well separated.  Eigenvalue clustering and repetition are common in multidimensional PDEs.  Let us consider the discrete system ${\K\uu=\lm\M\uu}$ as a generalized Hermitian eigenproblem~(GHEP), where the term \textit{generalized} distinguishes it from the standard Hermitian eigenproblem ${\textbf{A}\uu=\lm\uu}$.  A well-established recommendation in the literature (e.g.,~\cite{ Ericsson1980, NourOmid1987, Grimes1994, Demmel2000, XUE2011}) is to first perform a spectral transformation (ST), and then solve the shifted eigenproblem ${(\K-\sg\M)\.\uu=(\lm-\sg)\.\M\uu}$. This spectral shift results in a fast convergence when calculating eigenvalues near the shift $\sg$.  A more efficient way in eigenpairs approximation is to solve a  \mbox{``shift-and-invert"} spectral transformed eigenproblem ${(\K-\sg\M)^{-1}\M\uu=\theta\.\uu}$, with ${\theta=1/(\lm-\sg)}$.  This approach preserves the separation of eigenvalues near $\sg$ to reach an accurate eigensolution.  In many practical cases, we seek all eigenpairs $\lm_i$ and $\uu_i$ within a given (large) interval ${\lm_i\in[\lm_s,\lm_e]}$, where either $\lm_s$ or $\lm_e$ can be infinite.  Hence, we select several shifts ${\sg_k\in[\lm_s,\lm_e]}$ to preserve the convergence rate for eigenvalues far from $\sg$.  We employ a ``spectrum slicing technique" (see, e.g.,~\cite{ Campos2012}) to dynamically select $\sg_k$s and find all eigenpairs with the true multiplicities and without incurring in unnecessary computational efforts.  \Sec{\ref{sec.Eigensolution}} provides more algorithmic details of the eigenanalysis.

The factorization of the ST matrix ${\K-\sg_k\M}$ for each $\sg_k$ is a major component of computational effort that the eigenanalysis requires, especially when dealing with moderate to large algebraic systems.  Once we compute this factorization, the computation of the Krylov subspace requires several multiplications of the operator matrix ${(\K-\sg_k\M)^{-1}\M}$ by vectors.  These multiplications consist of two steps, namely the forward/backward eliminations of the respective factorized forms of the ST matrices, and the products of matrix $\M$ by vectors. Let $\N$ be the total number of degrees of freedom in the system.  Using maximum-continuity IGA, the computational cost of factorization is $\O(\N^{1.5}\p^3)$ and $\O(\N^{2}\p^3)$ for 2D and 3D systems, respectively.  The cost of performing forward/backward eliminations is $\O(\N\p^2)$ and $\O(\N^{1.33}\p^2)$ for 2D and 3D systems, respectively, while it is $\O(\N\p^2)$ and $\O(\N\p^3)$ for multiplying the $\M$ matrix by vectors in 2D and 3D cases, respectively (cf.~\cite{ Collier2012, Collier2013}).  

In \Sec{\ref{sec.CostEstimate}}, we show that rIGA computes a fixed number of eigenpairs faster than maximum-continuity IGA.  When using multifrontal direct solvers, rIGA reduces the factorization cost by up to ${\O(\p^{\.2})}$ for large domains.  Also, the cost of the forward/backward eliminations decreases by $\O(\p)$, since the factorized form of the ST matrix has fewer nonzero entries in rIGA versus its IGA counterpart (see~\cite{ Garcia2017}). Nevertheless, the matrix multiplication of $\M$ by vectors is slightly more expensive for rIGA as it has more nonzero entries than its {\em smooth} counterpart.  There are other costs such as vector--vector products which grow when using rIGA. However, their contribution to the total cost is irrelevant compared to those of the factorization and matrix--vector multiplications, as we mentioned above.  

In practice, the numerical tests show that the total computational cost of the eigensolution decreases by a factor of $\O(\p)$ when employing the rIGA discretization.  While our theoretical analysis shows that an improvement of ${\O(\p^{\.2})}$ is asymptotically possible. That is, for sufficiently large problems,  the matrix factorization governs the solution cost.  Additionally, rIGA approximates better the first $\N_0$ eigenvalues and eigenfunctions, where $N_0$ can be as large as to the total number of degrees freedom in the smooth IGA discretization.  The improved accuracy is a consequence of the continuity reduction of the basis functions that enriches the Galerkin space and modifies the approximation properties of the smooth IGA approach.

The organization of the remainder of this paper is as follows.  \Sec{\ref{sec.Preliminaries}} defines the problem. \Sec{\ref{sec.IGArIGA}} briefly revisits the notation and definitions of smooth (maximum-continuity) IGA and rIGA frameworks.  \Sec{\ref{sec.Eigensolution}} describes the eigensolution algorithm for finding the eigenpairs, while \Sec{\ref{sec.CostEstimate}} derives theoretical cost estimates of the eigensolution under IGA and rIGA discretizations.  We provide implementation details in \Sec{\ref{sec.Implementation}}, numerical cost evaluations in \Sec{\ref{sec.Results}}, and accuracy assessments in \Sec{\ref{sec.Accuracy}}.  Finally, \Sec{\ref{sec.Conclusions}} draws conclusions.


\section{Preliminaries} 
\label{sec.Preliminaries}

\subsection{Model problem} 
\label{sub:Definition}

We consider the eigensolution of the Laplace operator in the unit square:
\begin{equation}
	\begin{aligned}
		&\mbox{Find~} \lm\in\R^+ \mbox{~and~} u\.(\xx)\in H_0^1(\Om), \mbox{~satisfying}\\
		&\left\{\begin{array}{rr}
		\nabla^2u+\lm u=0\,,  &\\
		u\.(\xx)|_{\.\partial\Om}=0\,, & \xx\in\Om:(0,1)^{\.d}\,,\end{array}\right.
	\end{aligned}
	\label{eq.helmholtz}
\end{equation}
where ${d=2\rm{~or~}3}$ is the space dimension.  The exact eigenvalues $\lme$ and eigenfunctions $\ue$ in 2D and 3D are
\begin{equation}
	\begin{aligned}
		2\D&: & \lme_{ij}  &= \pi^{\,2}\pra{i^{\,2}+j^{\,2}}\,,  		& \ue_{ij}  &= 2\sin\pra{i\.\pi x}\sin\pra{j\.\pi y}\,,\\
		3\D&: & \lme_{ijk} &= \pi^{\,2}\pra{i^{\,2}+j^{\,2}+k^{\.2}}\,,	& \ue_{ijk} &= 2\sqrt{2}\sin\pra{i\.\pi x}\sin\pra{j\.\pi y}\sin\pra{k\.\pi \z}\,. 
	\end{aligned}
	\label{eq.exact}
\end{equation}

Considering the arbitrary test function ${v\.(\xx)\in H_0^1(\Om)}$, a weak form of \eq{\eqref{eq.helmholtz}} is 

\begin{equation}
	\int_\Om \nabla u \cdot \nabla v\,d\Om = \int_\Om \lm\.u\.v\,d\Om\,.
	\label{eq.weakform}
\end{equation}
 
\noindent 
Using symmetric bilinear forms

\begin{equation}
	\begin{aligned}
		a\.(u,v) &:= \int_\Om \nabla u \cdot \nabla v\,d\Om \,,\\
		   (u,v) &:= \int_\Om u\.v\,d\Om \,, 
	\end{aligned}
	\label{eq.biliniear}
\end{equation}

\noindent 
we write the weak formulation as 

\begin{equation}
	a\.(u,v) = \lm\.(u,v)\,.
	\label{eq.weakform2}
\end{equation}

Introducing a Galerkin discretization of the continuous eigenproblem leads to the following discrete form (superscript $h$ refers to the numerical computed eigenpairs)~\cite{ Strang1973},

\begin{equation}
	a\pra{\uh,v^{\.h}} = \lmh\pra{\uh,v^{\.h}}\,.
	\label{eq.weakform3}
\end{equation}


\subsection{Eigenvalue and eigenfunction errors} 
\label{sub:errors}

Considering the exact eigenpairs of \eq{\eqref{eq.exact}} and their numerical counterparts obtained by \eq{\eqref{eq.weakform3}}, in order to assess the accuracy of the spectral approximation, we study the eigenvalue error and the eigenfunction $L^2$ and energy error norms.  Using the Pythagorean eigenvalue error theorem (see~\cite{ Strang1973}), for the ${i\.}$th discrete eigenmode, we reach the following relation between the eigenvalue error and the eigenfunction errors in ${L^2}$ and energy norms:
\begin{equation}
	\lmh_i-\lme_i + \lme_i\normL{\.\uh_i-\ue_i}^{\.2} = \normE{\.\uh_i-\ue_i}^{\.2} \,,
	\label{eq.Pythagorean}
\end{equation}
where
\begin{equation}
	\begin{aligned}
		& \normL{\.\uh_i-\ue_i}^{\.2} =  \pra{\.\uh_i-\ue_i,\.\uh_i-\ue_i} \,,\\ 
		& \normE{\.\uh_i-\ue_i}^{\.2} = a\pra{\.\uh_i-\ue_i,\.\uh_i-\ue_i} \,.
	\end{aligned}
	\label{eq.norms}
\end{equation}
Based on the above theorem, we define for the $i\.$th eigenmode of our spectral approximation, its normalized eigenvalue error ($\EV$) and eigenfunction $L^2$ and energy norm errors ($\EFL$ and $\EFE$, respectively) as
\begin{equation}
	\begin{aligned}
		\EV  &:= \dfrac{\lmh_i-\lme_i}{\lme_i} 	\,,\\ 
		\EFL &:= \normL{\.\uh_i-\ue_i}^{\.2} 	\,,\\ 
		\EFE &:= \dfrac{\normE{\.\uh_i-\ue_i}^{\.2}}{\lme_i} = \EV+\EFL \,.
	\end{aligned}
	\label{eq.eigerrors}
\end{equation}


\section{Refined isogeometric analysis} 
\label{sec.IGArIGA}

We first review some basic concepts of maximum-continuity IGA discretizations.
For the sake of simplicity, we assume that the computational grid consists of 
a tensor-product mesh in ${[0,1]^d}$ with the same number of equally-spaced elements in each direction. For descriptions of how to use tensor-product discretization relying on mapped geometries refer to~\cite{Hughes2005,Cottrell2009,petiga,petigamf}.

\subsection{IGA discretization} 
\label{sub:IGA}

\begin{figure}[h]
	\centering
	\begin{subfigure}{.47\textwidth}\centering
		\begin{overpic}[width=1\linewidth,trim={6cm 0.25cm 10cm -0.75cm},clip]{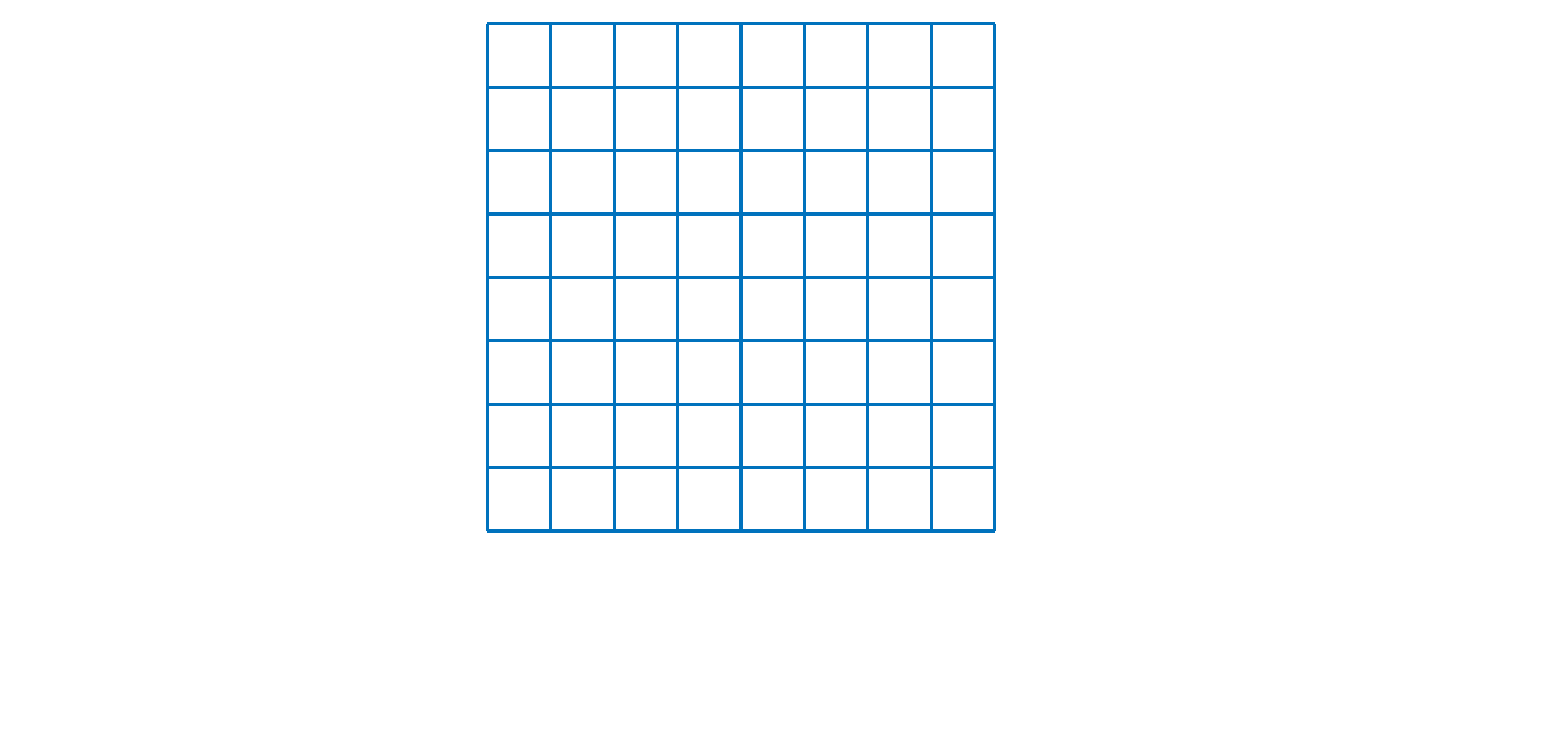}
			\put(16.75,4){
				\begin{overpic}[width=0.71\linewidth,trim={1cm 1.5cm 2cm 8cm},clip]{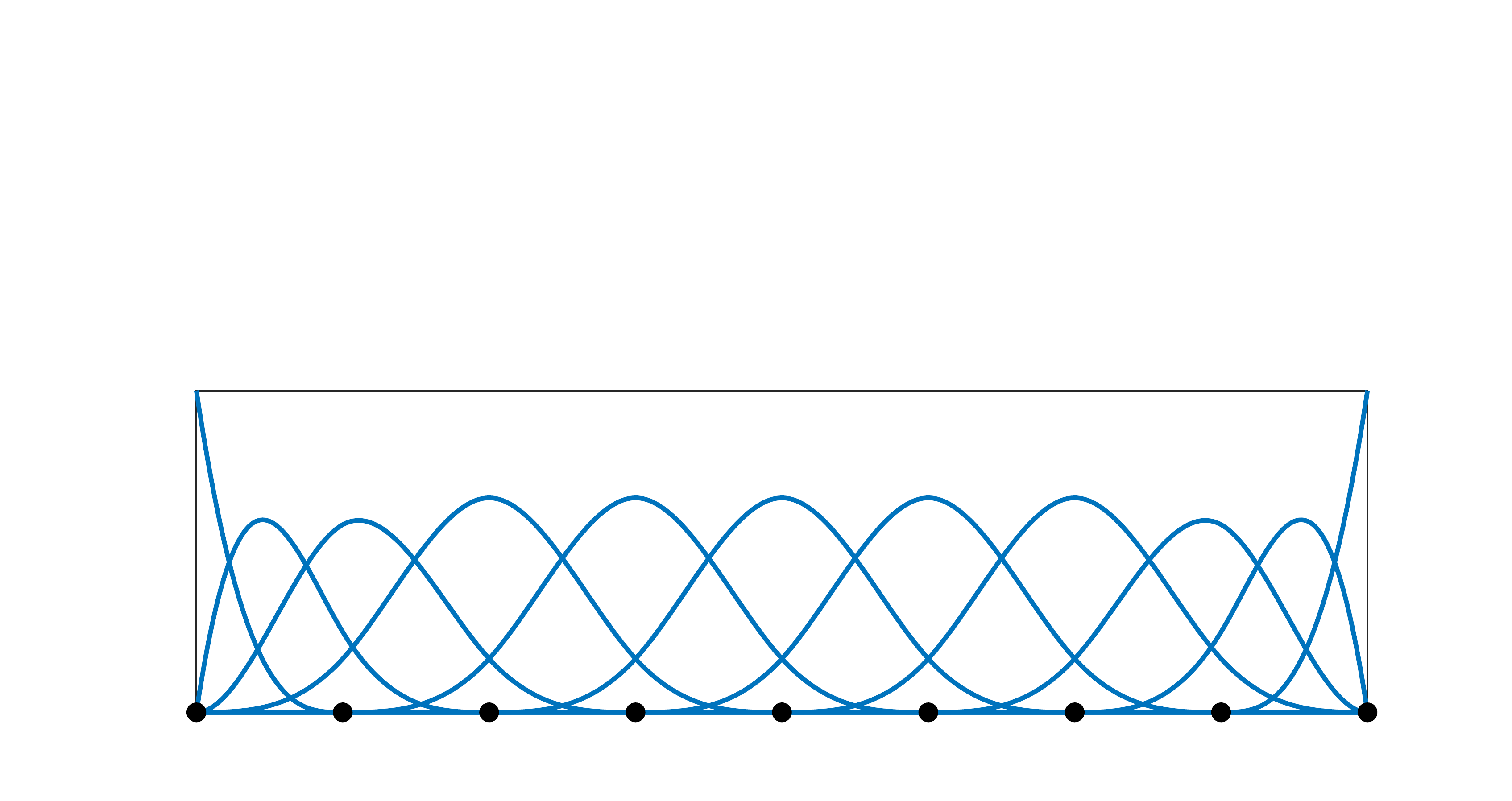}
					\put(9,-3){\footnotesize $\x_{\p}$}
					\put(19,-3){\footnotesize $\x_{\p+1}$}
					\put(94,-3){\footnotesize $\x_{\n}$}
					\put(52,-6){\footnotesize $\x$}
					\begin{turn}{90}
						\put(6,-7){\footnotesize $\B{i}{\p}(\x)$}
					\end{turn}
				\end{overpic}}
			\begin{turn}{-90}
				\put(-93.5,5){
					\begin{overpic}[width=0.71\linewidth,trim={1cm 1.5cm 2cm 8cm},clip]{figs/BasisN8L0-eps-converted-to}
						\put(9,-3){\footnotesize $\y_{\p}$}
						\put(19,-3){\footnotesize $\y_{\p+1}$}
						\put(94,-3){\footnotesize $\y_{\n}$}
						\put(52,-6){\footnotesize $\y$}
						\begin{turn}{90}
							\put(6,-7){\footnotesize $\B{j}{\p}(\y)$}
						\end{turn}
					\end{overpic}}
			\end{turn}
		\end{overpic}
	\end{subfigure}
	\begin{subfigure}{.50\textwidth}\centering
		\begin{overpic}[width=1\linewidth,trim={1.5cm -0.25cm 1.5cm -0.75cm},clip]{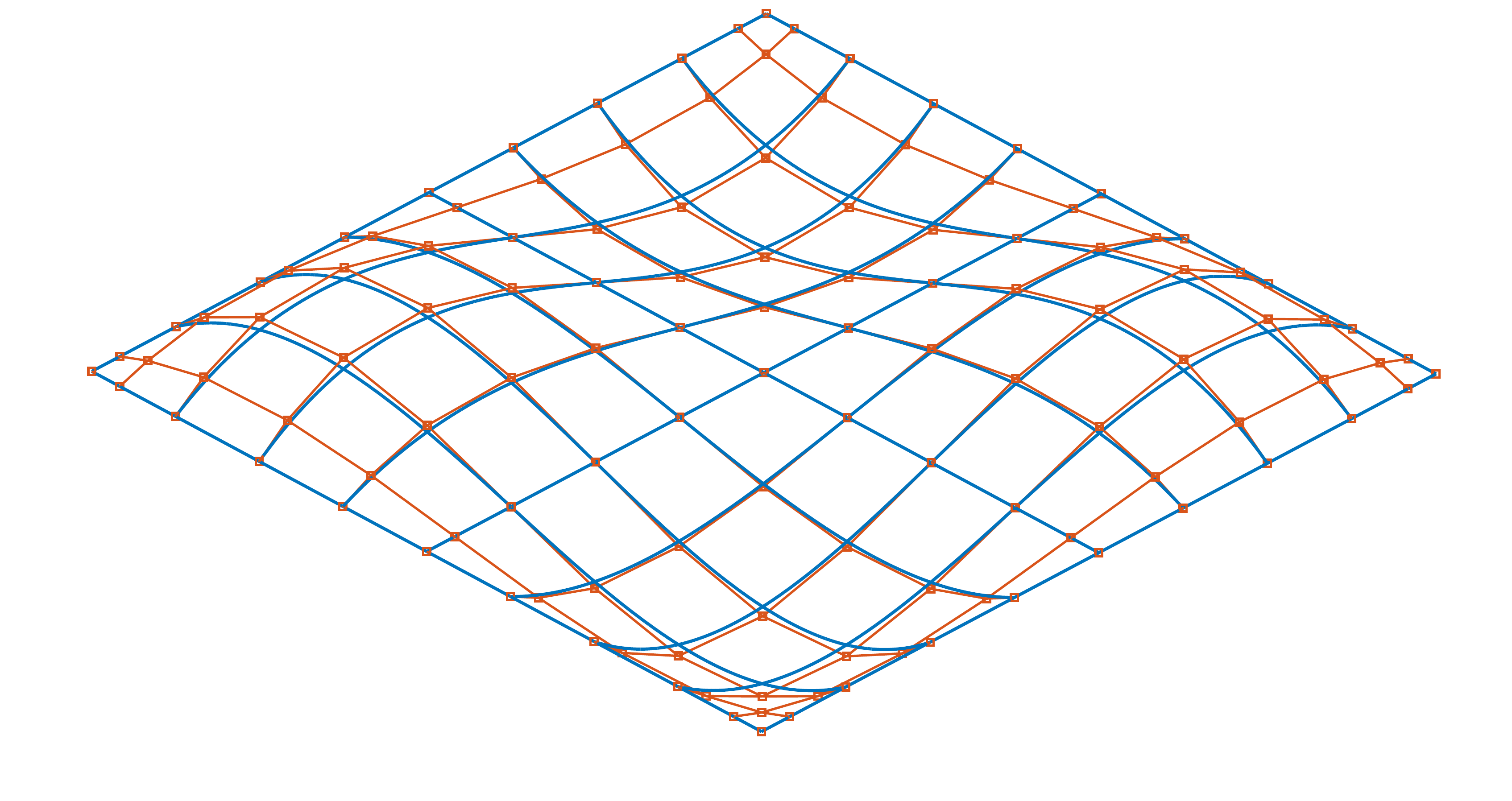}
			\put(50.5,-10){\vector(-2,1){10}}
			\put(50.5,-10){\vector(2,1){10}}
			\put(50.5,-10){\vector(0,1){10}}
			\put(61,-4){\footnotesize $\x$}
			\put(41.5,-3){\footnotesize $\y$}
			\put(51,0.5){\footnotesize $\uh(\x,\y)$}
			\put(13.5,36.5){\line(-1,1){10}}
			\put(40,45){\line(-1,1){10}}
			\put(26,57){\footnotesize $\uh(\x,\y)$}
			\put(2,48){\footnotesize $\Uu_{ij}$}
		\end{overpic}
	\end{subfigure}
	\caption{Left: a uniform bicubic ${8\times8}$ computational mesh in 2D with its respective basis functions.
			 Right: a spatial representation of the approximate solution field $\uh(\x,\y)$ with the net of control variables $\Uu_{ij}$ (red) and the projection of elements on the solution surface (blue). 
			 Control variables coincide with the solution surface on the domain boundaries, but not necessarily at interior points.}
	\label{fig.N8L0}	
\end{figure}
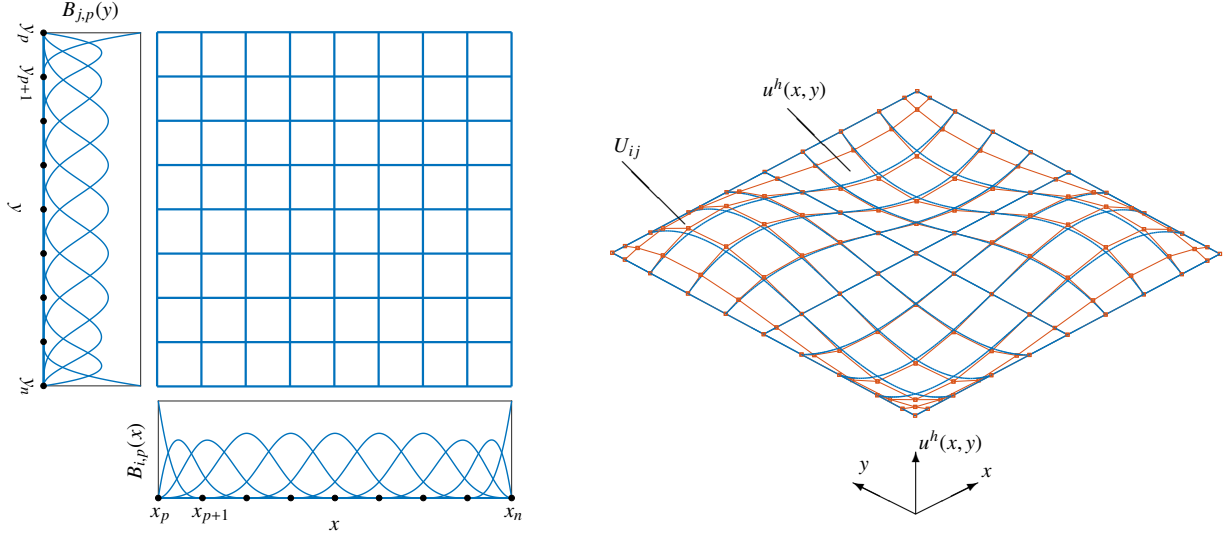

We discretize our computational domain with a uniform mesh of ${\ne^{\.d}}$ elements, being $\ne$ the number of elements in each direction (see \fig{\ref{fig.N8L0}} for the 2D case).  The approximate solution field ${\uh(\xx)}$ is described by the B-spline representation as
\begin{equation}
	\begin{aligned}
		2\D&: & \uh(\x,\y)	 &=\sum_{i\.=\.0}^{\n-1} \sum_{j\.=\.0}^{\n-1} \B{i}{\p}(\x)\B{j}{\p}(\y)\,\Uu_{ij}\,, \\
		3\D&: & \uh(\x,\y,\z)&=\sum_{i\.=\.0}^{\n-1} \sum_{j\.=\.0}^{\n-1} \sum_{k\.=\.0}^{\n-1} \B{i}{\p}(\x)\B{j}{\p}(\y)\B{k}{\p}(\z)\,\Uu_{ijk}\,, \\
	\end{aligned}
	\label{eq.u}
\end{equation}

\noindent where $\UU$ is the ${d\.}$th-order tensor of control variables (i.e., degrees of freedom), 
$\p$ is the polynomial degree of the basis functions (equal in all directions), and ${\n=\ne+\p}$ is the number of control variables in each direction (with maximum continuity of $C^{\p-1}$ in basis functions). The parameter space is then characterized by the knot sequence $\X$ with single multiplicities, given for the $\x$ direction by 

\begin{equation}
	\X=[\underbrace{0,0,...,0}_{\p+1},\x_{\p+1},\x_{\p+2}...,\x_{\n-1},\underbrace{1,1,...,1}_{\p+1}]\,,
	\label{eq.U}
\end{equation}

\noindent and the B-spline basis functions $\B{i}{\p}$ are expressed by the Cox--De Boor recursion formula~\cite{ TheNURBSBook} as 

\begin{equation}
	\begin{aligned}
		\B{i}{0}(\x)  &= \left\{\begin{array}{ll} 
			1 & \x_i\leq \x < \x_{i+1}\,, \\ 
			0 & \mbox{otherwise}\,, \end{array} \right. \\
		\B{i}{\p}(\x) &= \dfrac{\x-\x_i}{\x_{i+\p}-\x_i}\B{i}{{\p-1}}(\x)+
		\dfrac{\x_{i+\p+1}-\x}{\x_{i+\p+1}-\x_{i+1}}\B{{i+1}}{{\p-1}}(\x)\,.
	\end{aligned}
	\label{eq.N}
\end{equation}

Applying the weak form given by \eq{\eqref{eq.weakform3}} to the above discretization, the generalized Hermitian eigenproblem~(GHEP) is
\begin{equation}
	\K\uuh=\lmh\M\uuh\,,
	\label{eq.GHEP}
\end{equation}
where $\K$ and $\M$ are real symmetric (or Hermitian) stiffness and mass matrices, respectively.

\begin{rem}
  The eigenvectors $\uuh$ in \eq{\eqref{eq.GHEP}} are the vectorial form of the tensor of control variables $\UU$. 	Hence, we obtain the tensor form by ${\UU = \mbox{Tensor}\.(\uuh)}$~\cite{ Gao2015}, and compute the numerical eigenfunction $\uh$ for each eigenmode by \eq{\eqref{eq.u}}.
  \label{rem.matvec}
\end{rem}

\subsubsection{Construction of multidimensional mass and stiffness matrices}

We define 1D stiffness and mass matrices in $x$, $y$, and $z$ directions as
\begin{equation}
  \begin{aligned}
    K^{\x}_{ij}  &:= \int_x \B{i}{\p}\pr(\x) \B{j}{\p}\pr(\x) \,d\x \,, &
    K^{\y}_{ij}  &:= \int_y \B{i}{\p}\pr(\y) \B{j}{\p}\pr(\y) \,d\y \,, &
    K^{\z}_{ij}  &:= \int_z \B{i}{\p}\pr(\z) \B{j}{\p}\pr(\z) \,d\z \,, \\
    M^{\x}_{ij}  &:= \int_x \B{i}{\p}   (\x) \B{j}{\p}   (\x) \,d\x \,, &
    M^{\y}_{ij}  &:= \int_y \B{i}{\p}   (\y) \B{j}{\p}   (\y) \,d\y \,, &
    M^{\.\z}_{ij}&:= \int_z \B{i}{\p}   (\z) \B{j}{\p}   (\z) \,d\z \,.
  \end{aligned}
  \label{eq.1Dmatrix}
\end{equation}
We build the 2D and 3D system matrices with the following formulae
(see, e.g.,~\cite{ Gao2015}):
\begin{equation}
  \begin{aligned}
    2\D&: & \K &= \M^{\.\y}\otimes\K^{\x}+\K^{\y}\otimes\M^{\x} \,,
    & \M &= \M^{\.\y}\otimes\M^{\x} \,,\\ 
    3\D&: & \K &= \M^{\z}\otimes\M^{\.\y}\otimes\K^{\x}+\M^{\z}\otimes\K^{\y}\otimes\M^{\x}+\K^{\z}\otimes\M^{\.\y}\otimes\M^{\x} \,,
    & \M &= \M^{\z}\otimes\M^{\.\y}\otimes\M^{\x} \,,
  \end{aligned}
  \label{eq.2Dmatrix}
\end{equation}
where $\otimes$ indicates the Kronecker product.	\label{lem.tensorproduct}


\subsection{Refined IGA discretizations} 
\label{sub:rIGApartitioning}

\begin{figure}[h]
	\centering
	\begin{subfigure}{.325\textwidth}\centering
		\begin{overpic}[width=.952\linewidth,trim={5.9cm 1cm 11.5cm -0.18cm} ,clip]{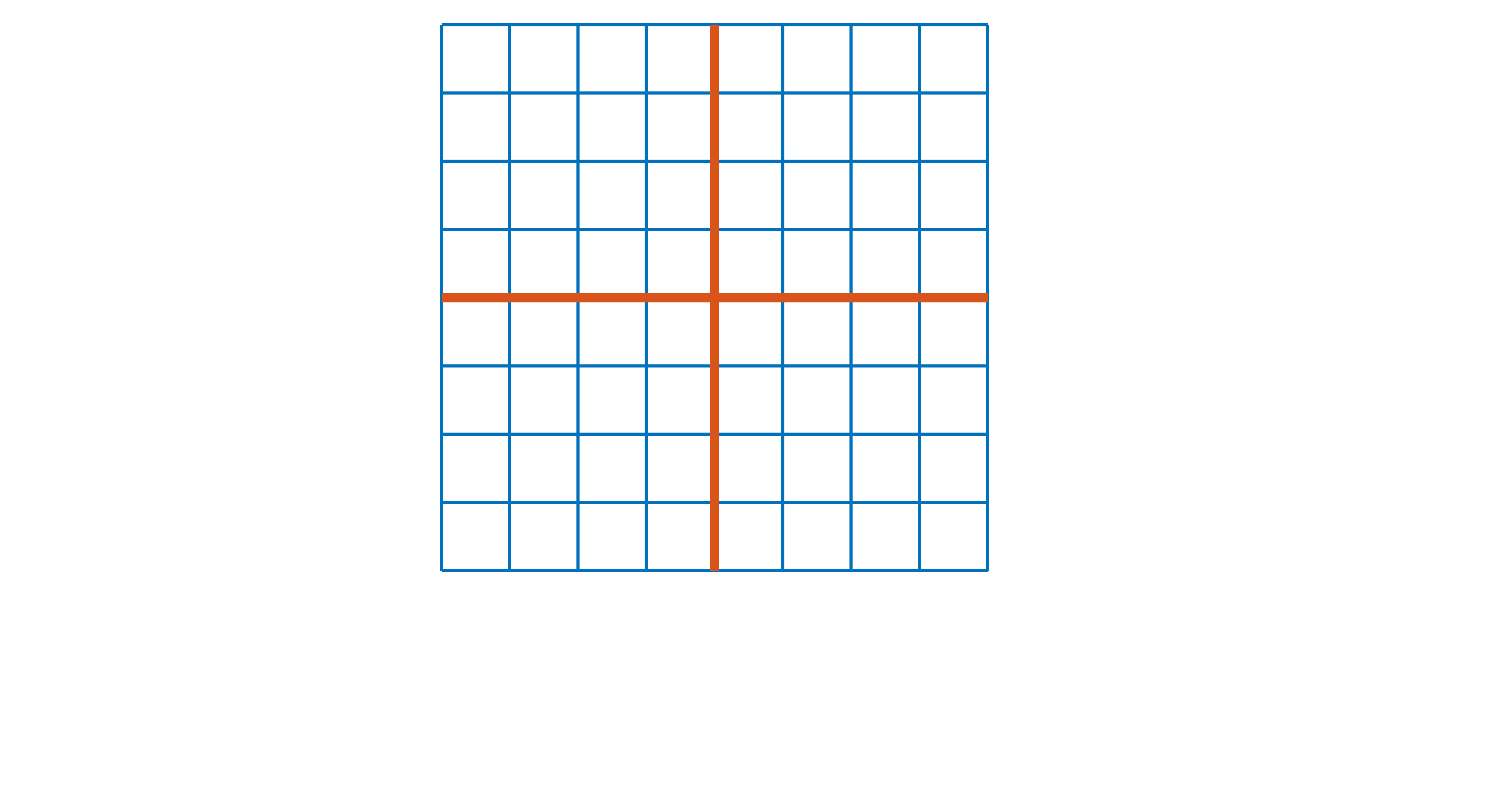}
			\put(20,0){\begin{overpic}[width=0.77\linewidth,trim={3cm 1.5cm 2cm 8cm},clip]{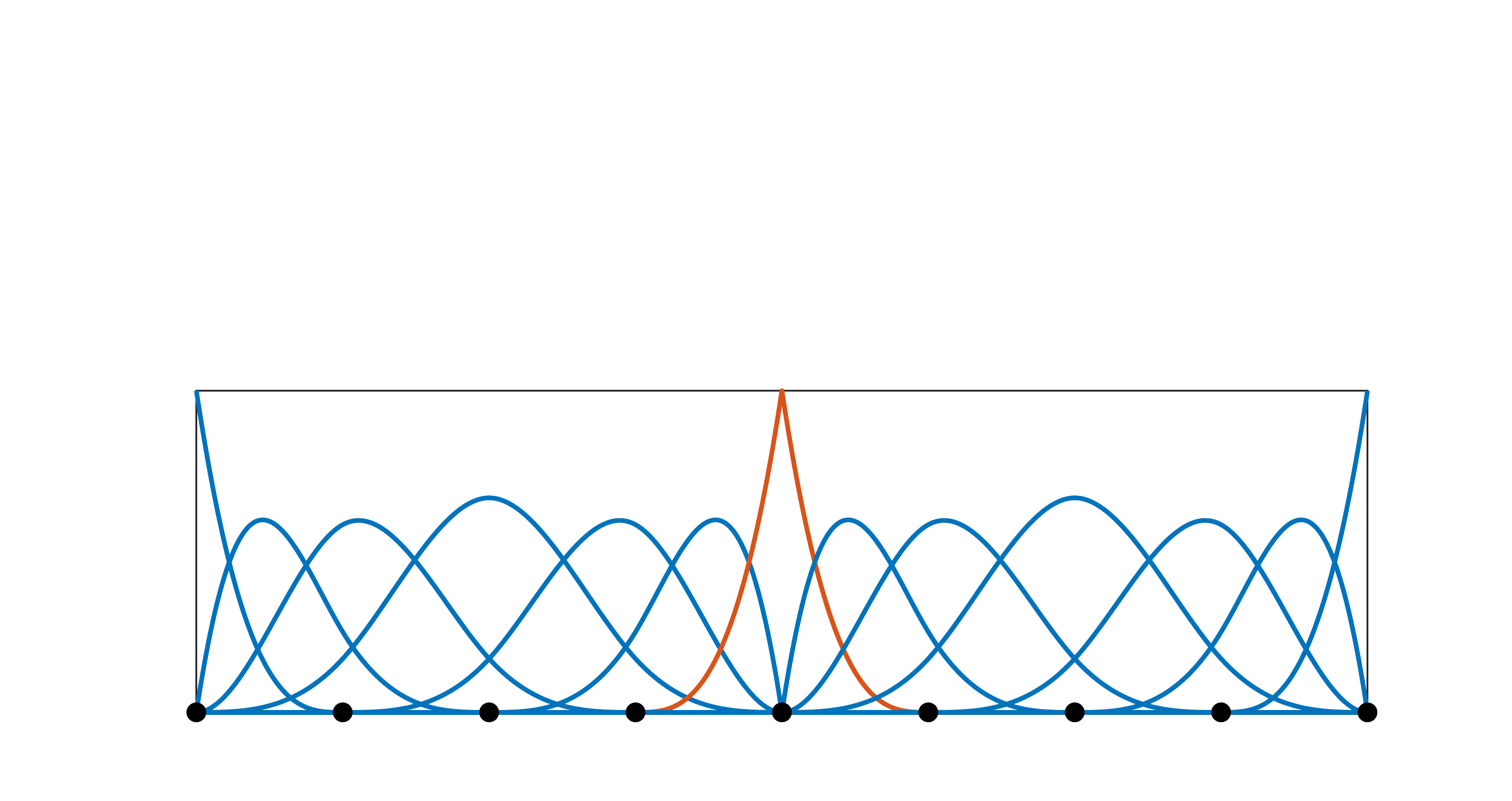}\end{overpic}}
			\begin{turn}{-90}
				\put(-100,0){\begin{overpic}[width=0.77\linewidth,trim={3cm 1.5cm 2cm 8cm},clip]{figs/BasisN8L1-eps-converted-to}\end{overpic}}
			\end{turn}
		\end{overpic}
	\caption*{\qquad\qquad (a) $\l=1~(\BS=4)$}
	\end{subfigure}
	\begin{subfigure}{.325\textwidth}\centering
		\begin{overpic}[width=.952\linewidth,trim={5.9cm 1cm 11.5cm -0.18cm}, clip]{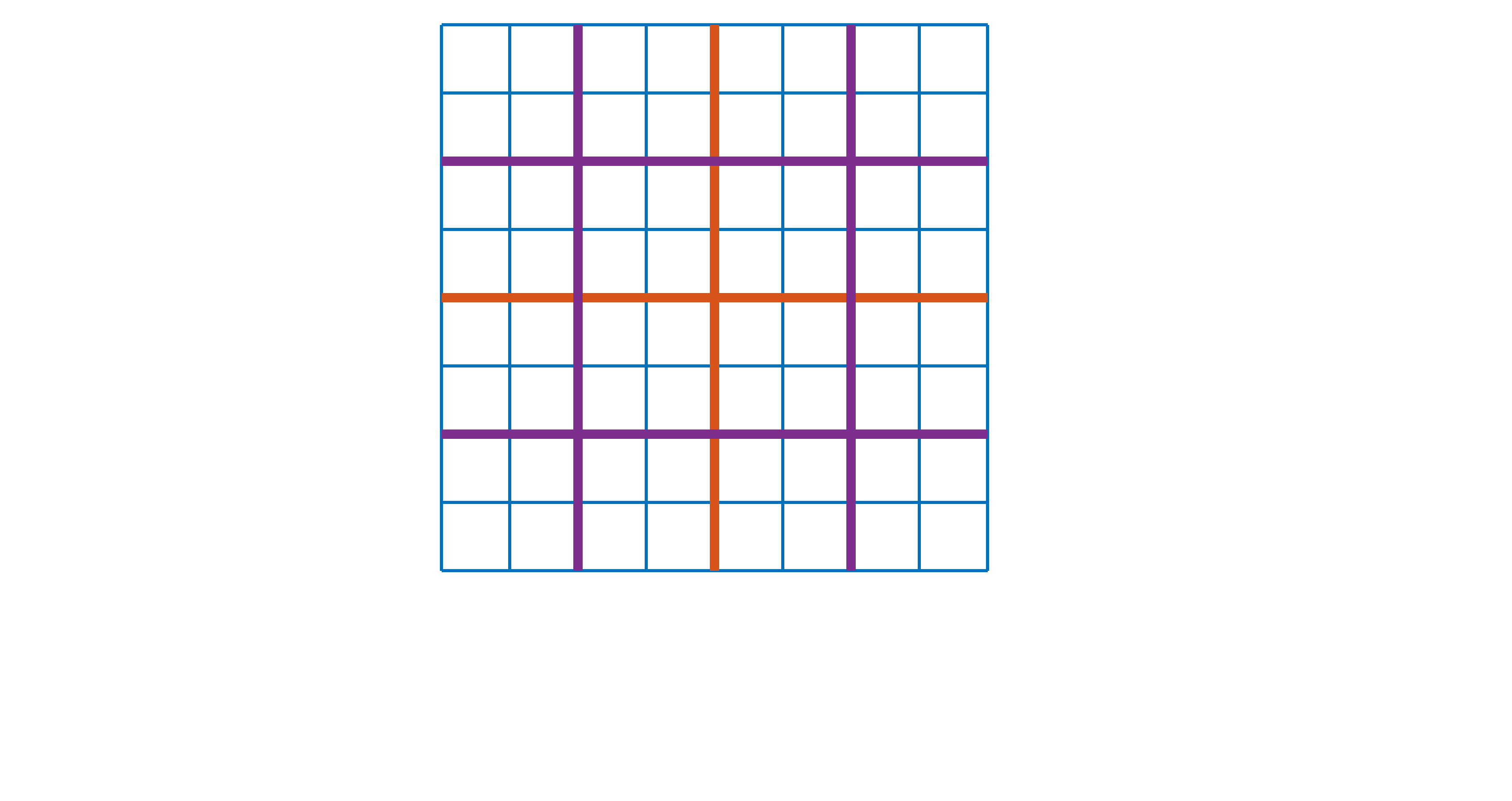}
			\put(20,0){\begin{overpic}[width=0.77\linewidth,trim={3cm 1.5cm 2cm 8cm},clip]{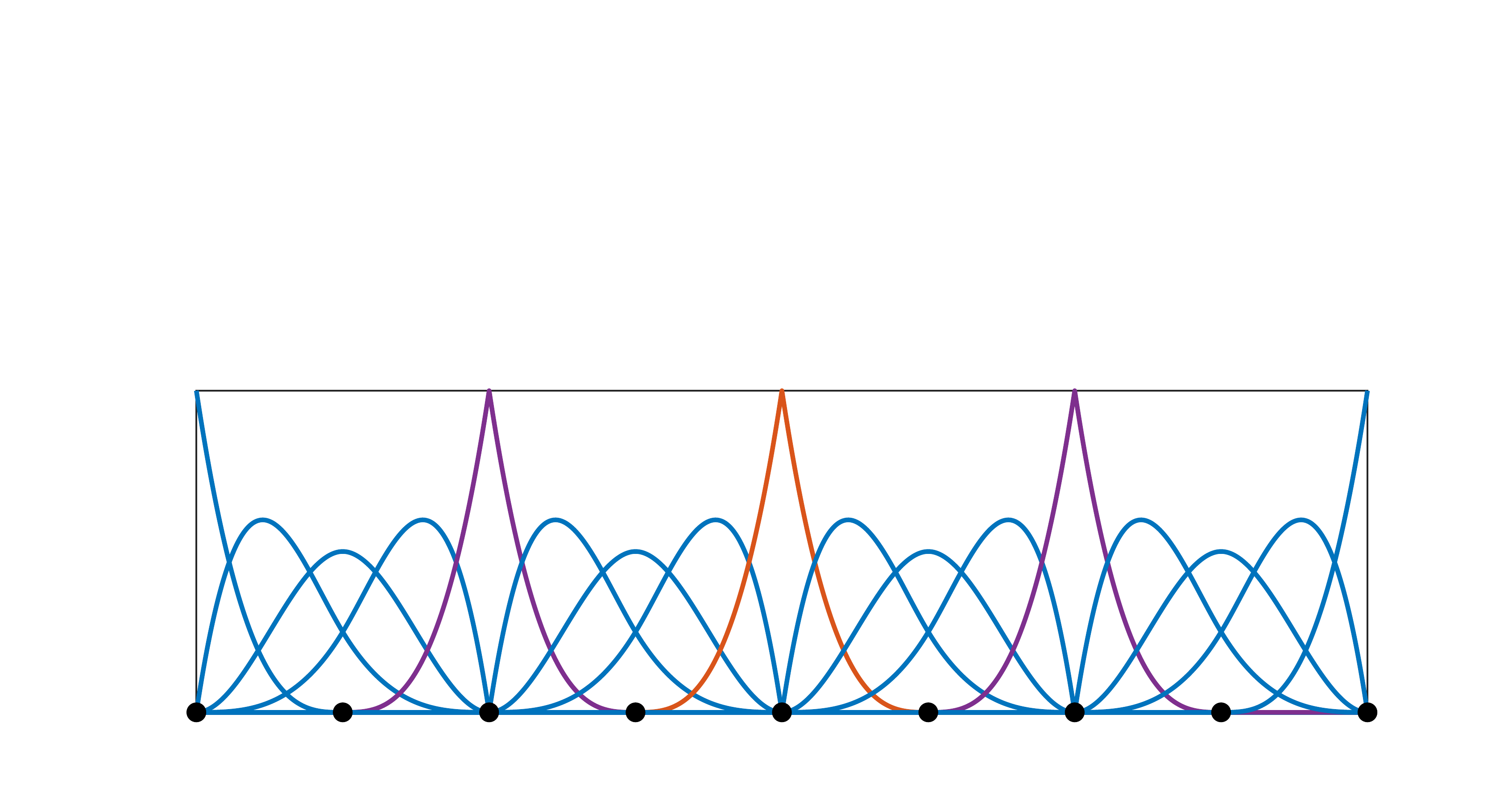}\end{overpic}}
			\begin{turn}{-90}
				\put(-100,0){\begin{overpic}[width=0.77\linewidth,trim={3cm 1.5cm 2cm 8cm},clip]{figs/BasisN8L2-eps-converted-to}\end{overpic}}
			\end{turn}
		\end{overpic}
	\caption*{\qquad\qquad (b) $\l=2~(\BS=2)$}
	\end{subfigure}
	\begin{subfigure}{.325\textwidth}\centering
		\begin{overpic}[width=.952\linewidth,trim={5.9cm 1cm 11.5cm -0.18cm}, clip]{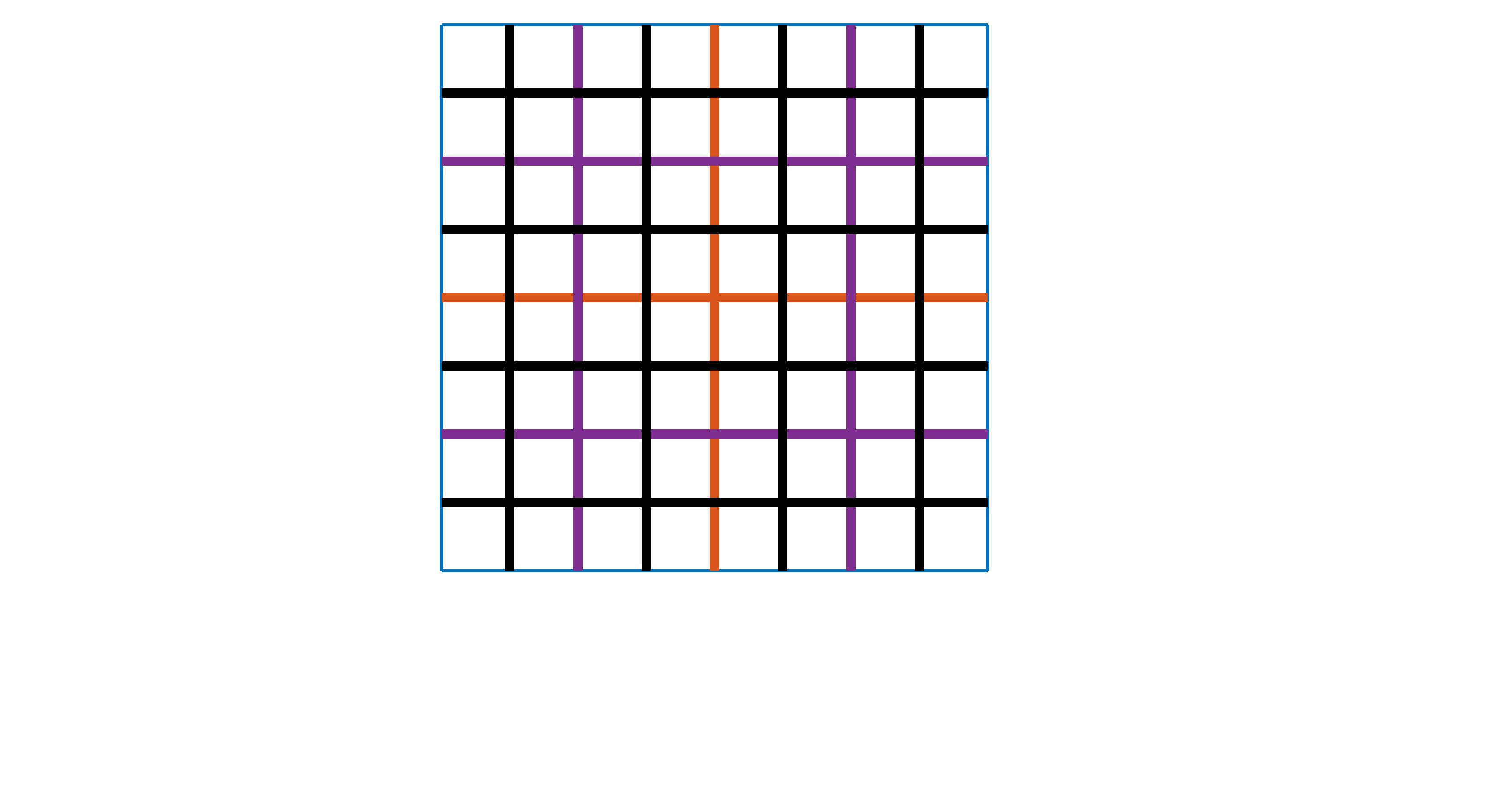}
			\put(20,0){\begin{overpic}[width=0.77\linewidth,trim={3cm 1.5cm 2cm 8cm},clip]{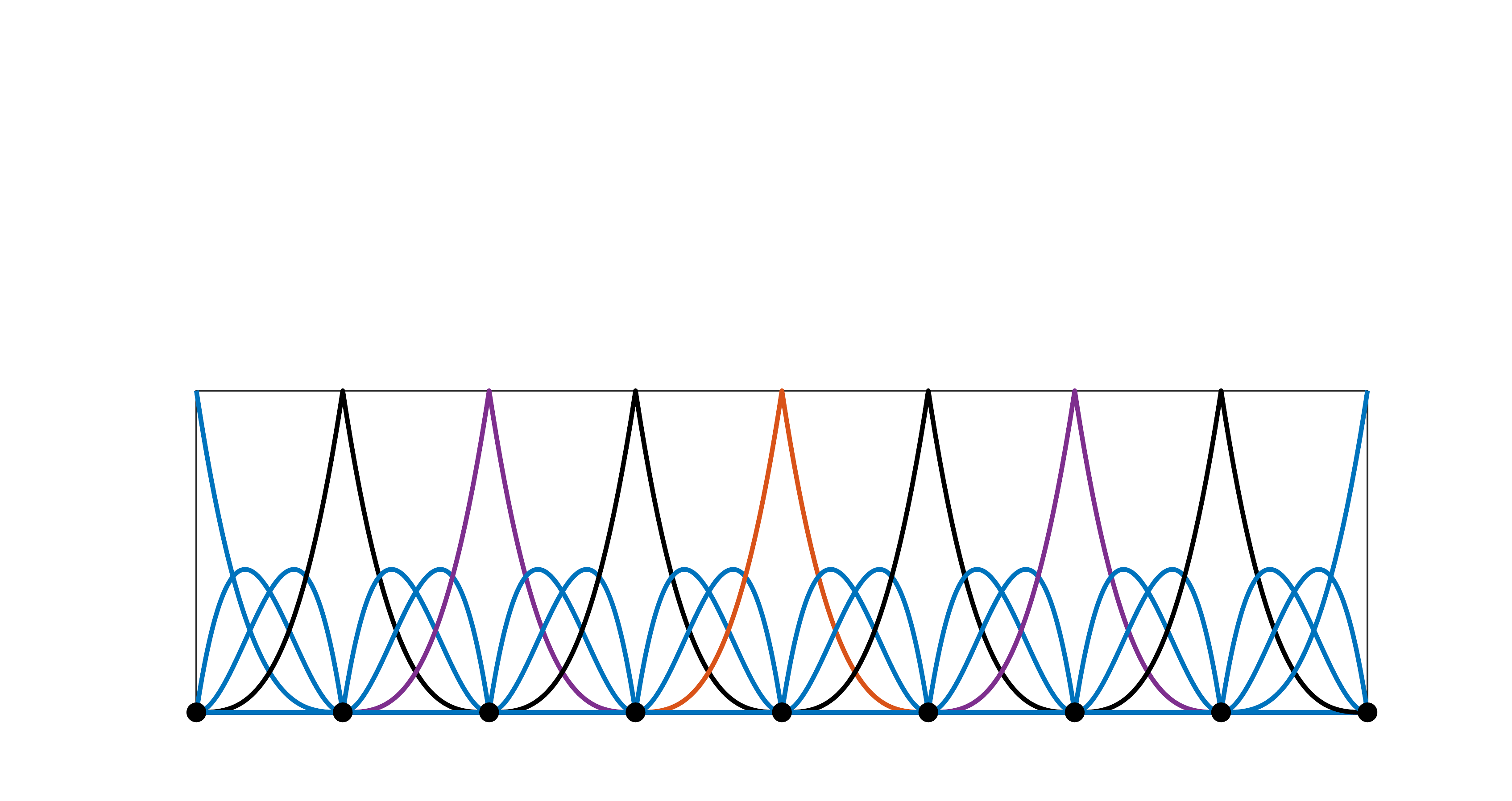}\end{overpic}}
			\begin{turn}{-90}
				\put(-100,0){\begin{overpic}[width=0.77\linewidth,trim={3cm 1.5cm 2cm 8cm},clip]{figs/BasisN8L3-eps-converted-to}\end{overpic}}
			\end{turn}
		\end{overpic}
	\caption*{\qquad\qquad (c) $\l=3~(\BS=1)$}
	\end{subfigure}
	\caption{Different partitioning levels ${(\l=1,2,3)}$ and respective $C^{\.0}$ separators (in red, purple and black) for the bicubic ${8\times8}$ grid of \fig{\ref{fig.N8L0}}.}
	\label{fig.N8L123}	
\end{figure}

In refined IGA, we improve the approximation space to reduce the computational cost of the solution as well as to approximate better the solution. That is, rIGA reduces the continuity of certain basis functions to reduce the interconnection between degrees of freedom of the mesh~\cite{ Garcia2017}.  By increasing the multiplicity of some existing knots up to the degree of B-spline bases in the \mbox{$h$-refinement} sense, the continuity and support size of the basis functions decreases without adding new elements.  The resulting zero-continuity basis functions partition the computational space into interconnected blocks separated by $C^{\.0}$ hyperplanes.  This connectivity reduction significantly reduces the solver cost by reducing the cost of the matrix factorization as well as the forward/backward elimination (see~\cite{ Garcia2017, Garcia2019}).  The knot insertion steps add new control variables and, therefore, enrich the Galerkin space, modifying the spectral approximation properties of the IGA approach.  \fig{\ref{fig.N8L123}} describes three symmetric partitioning levels with the respective blocksizes of 4, 2, and 1 for the bicubic ${8\times8}$ mesh of \fig{\ref{fig.N8L0}} ($\BS$ is the number of elements of $C^{\.0}$ blocks in each direction).

\begin{rem}
  For simplicity, we assume the mesh size in each direction is a power of two (i.e., $\ne=2^s$). This assumption allows us to split the mesh symmetrically and obtain $2^{\.\l}$ blocks (i.e., macroelements) in each direction with $\BS$ $2^{s-\.\l}$, where ${\l=0,1,...,s}$ is the partitioning level.  Here, ${\l=0}$ refers to the maximum-continuity IGA with $C^{\p-1}$ continuity everywhere, while ${\l=s}$ is equivalent to FEA with $C^{\.0}$ continuity at all element interfaces (knot lines).
\end{rem}

\begin{figure}[h]
  \centering
  \begin{subfigure}{.245\textwidth}\centering
    \begin{overpic}[width=1\linewidth,trim={2.5cm 1cm 2cm 0.5cm},clip]{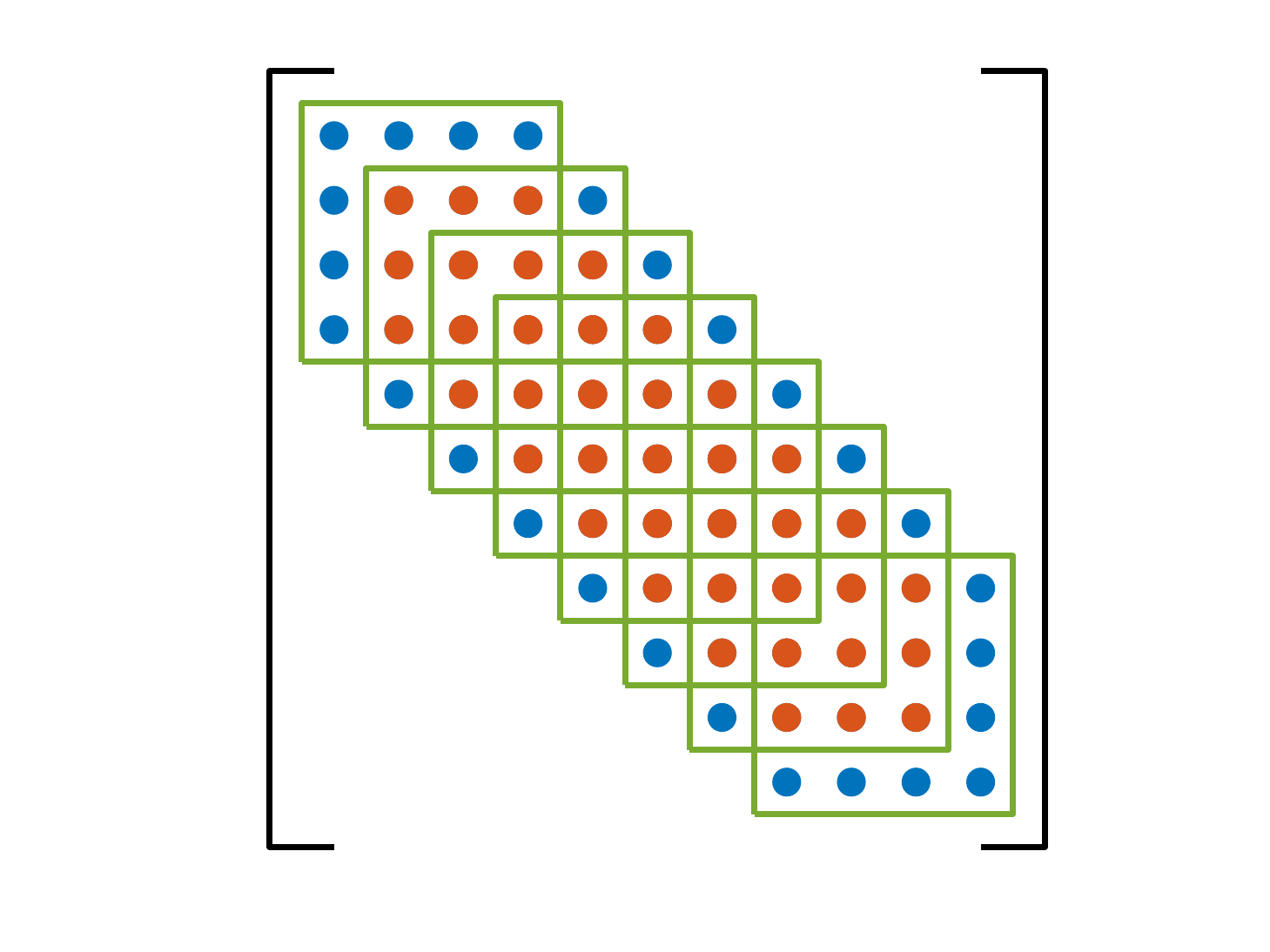}\end{overpic}
    \caption{$C^{\p-1}$ IGA}
  \end{subfigure}
  \begin{subfigure}{.245\textwidth}\centering
    \begin{overpic}[width=1\linewidth,trim={2.5cm 1cm 2cm 0.5cm},clip]{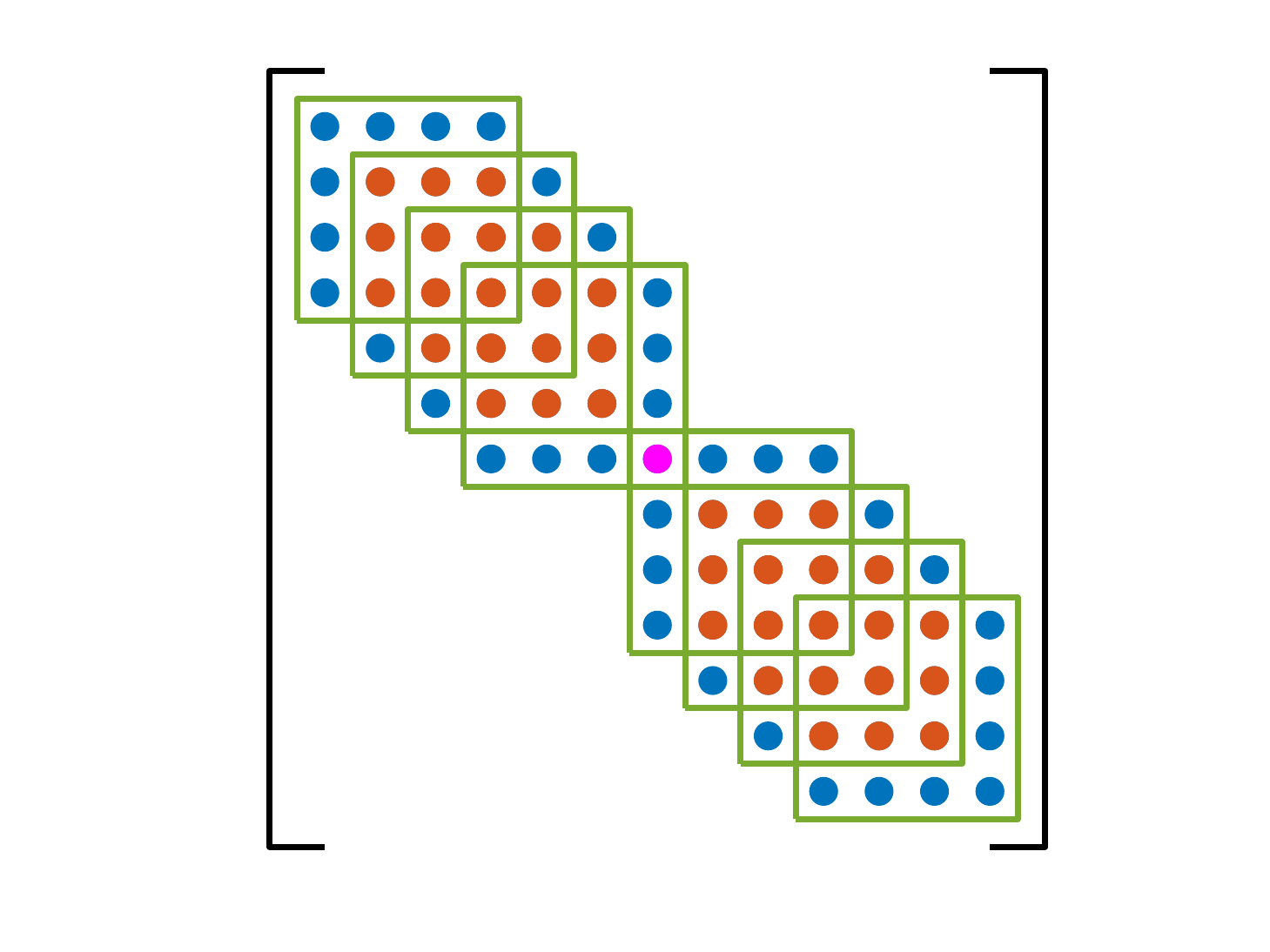}\end{overpic}
    \caption{rIGA with $\l=1$}
  \end{subfigure}
  \begin{subfigure}{.245\textwidth}\centering
    \begin{overpic}[width=1\linewidth,trim={2.5cm 1cm 2cm 0.5cm},clip]{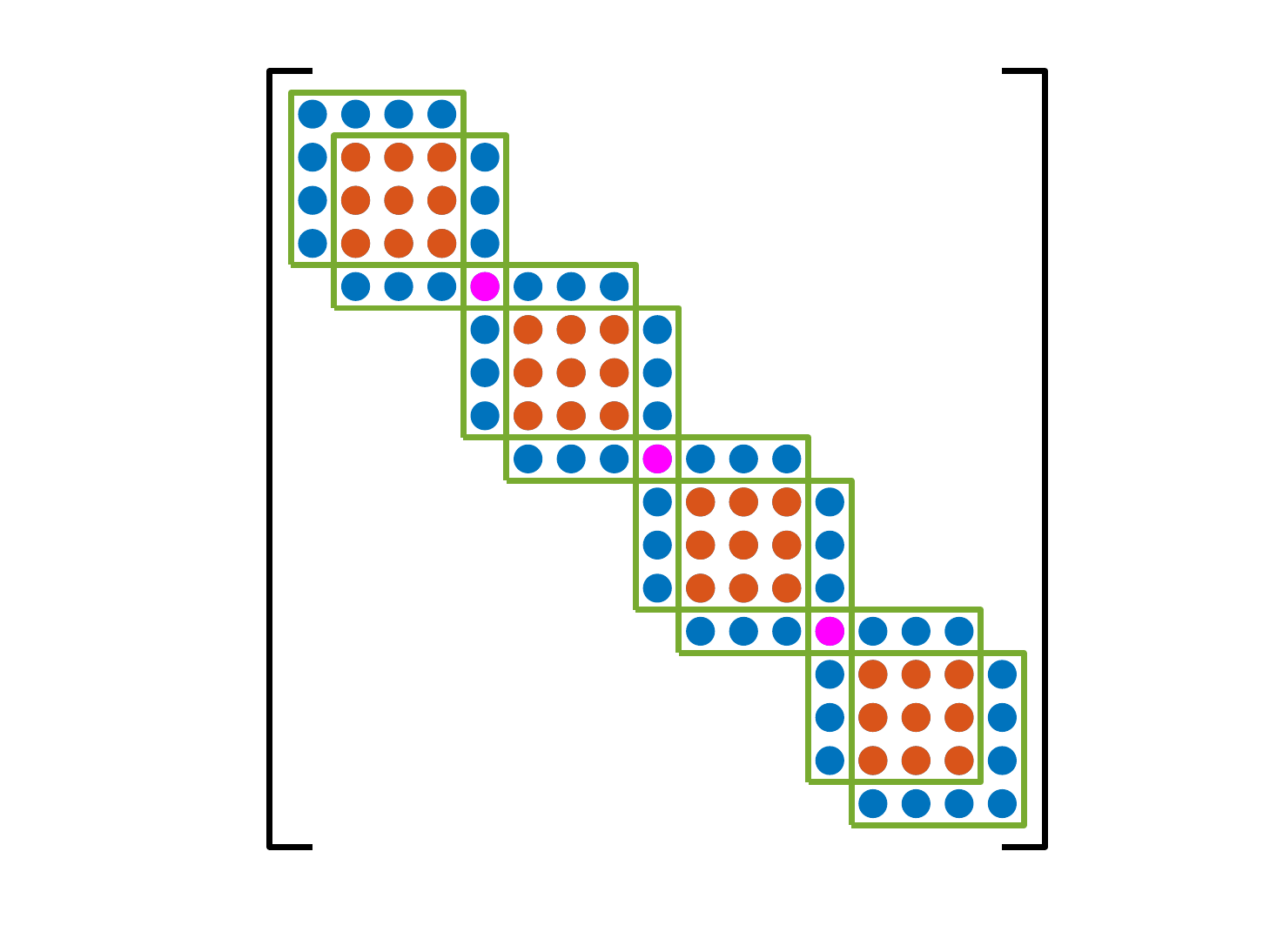}\end{overpic}
    \caption{rIGA with $\l=2$}
  \end{subfigure}
  \begin{subfigure}{.245\textwidth}\centering
    \begin{overpic}[width=1\linewidth,trim={2.5cm 1cm 2cm 0.5cm},clip]{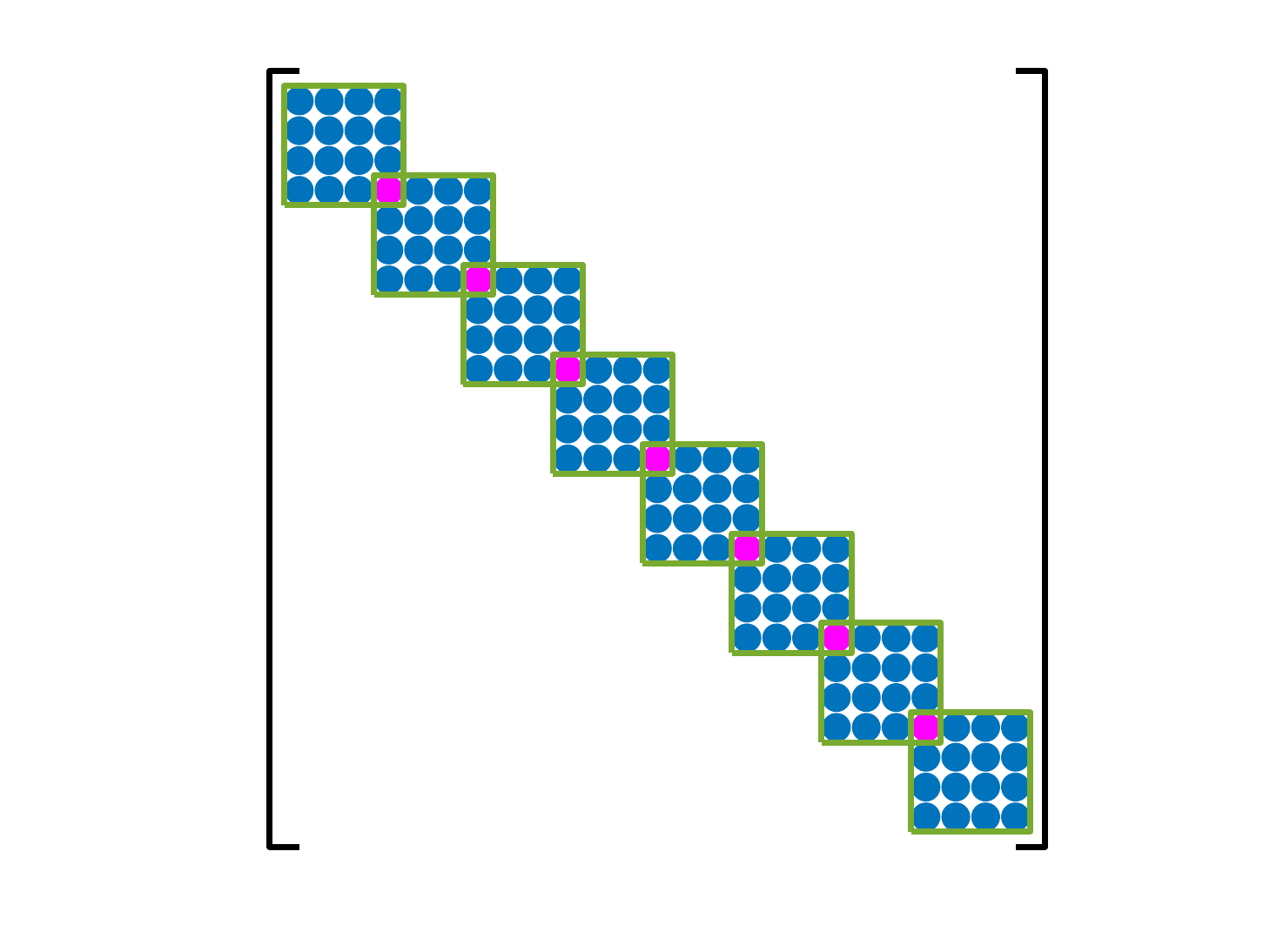}\end{overpic}
    \caption{rIGA with $\l=3$ (FEA)}
  \end{subfigure}
  \caption{
    Matrix patterns of a cubic eight-element domain in 1D under different discretizations (green squares delimit elemental matrices). (a)~maximum-continuity $C^{\p-1}$ IGA, (b)~rIGA with ${\BS=4}$, (c)~rIGA with ${\BS=2}$, and (d) rIGA with ${\BS=1}$, which is equal to the cubic FEA discretization. Uniform high-continuity implies strong interconnection between degrees of freedom (red dots) in (a). The interconnection weakens for the interior elements of each macroelement under rIGA discretizations in (b) and (c).  Magenta dots denote the ${C^{\.0}}$ interconnection between elements.
      }
	\label{fig.matrixN8}	
\end{figure}

\fig{\ref{fig.matrixN8}} depicts the matrix patterns of a 1D domain under different discretizations with ${\ne=8}$ and ${\p=3}$.  The figure shows the strong interconnectivity between degrees of freedom for maximum-continuity IGA. The figure also shows that rIGA partitioning weakens this connectivity accelerating the solution and reducing the memory required to solve.  For the maximum partitioning level ${(\l=3)}$, each macroelement consists of one element, and the interconnection reduces to a single degree of freedom.
 
Increasing the multiplicity of an existing knot up to $\p$ (i.e., the degree of bases) adds ${\p-1}$ control variables to each direction.
Thus, the ${\l\.}$th partitioning level ${(\l>0)}$ adds ${[2^{\l-1}(\p-1)]^{\.d}}$ control variables.
Consequently, the total number of degrees of freedom to discretize \eq{\eqref{eq.helmholtz}} is 

\begin{equation}
 	\begin{aligned}
 		\rm{IGA} &: & \N &= \pra{\ne+\p-2}^{\.d} \,,\\
 		\rm{rIGA}&: & \N &= \bra{\ne+2^\l(\p-1)-1}^{\.d} \,.
 	\end{aligned}
 	\label{eq.dofs}
 \end{equation}


\section{Solving generalized Hermitian eigenproblems} 
\label{sec.Eigensolution}

\subsection{Shift-and-invert algorithm} 
\label{sub:shifinvert}

\begin{figure}[h]
	\centering
	\begin{overpic}[width=0.5\linewidth,trim={12cm 2.25cm 2cm 6.5cm},clip]{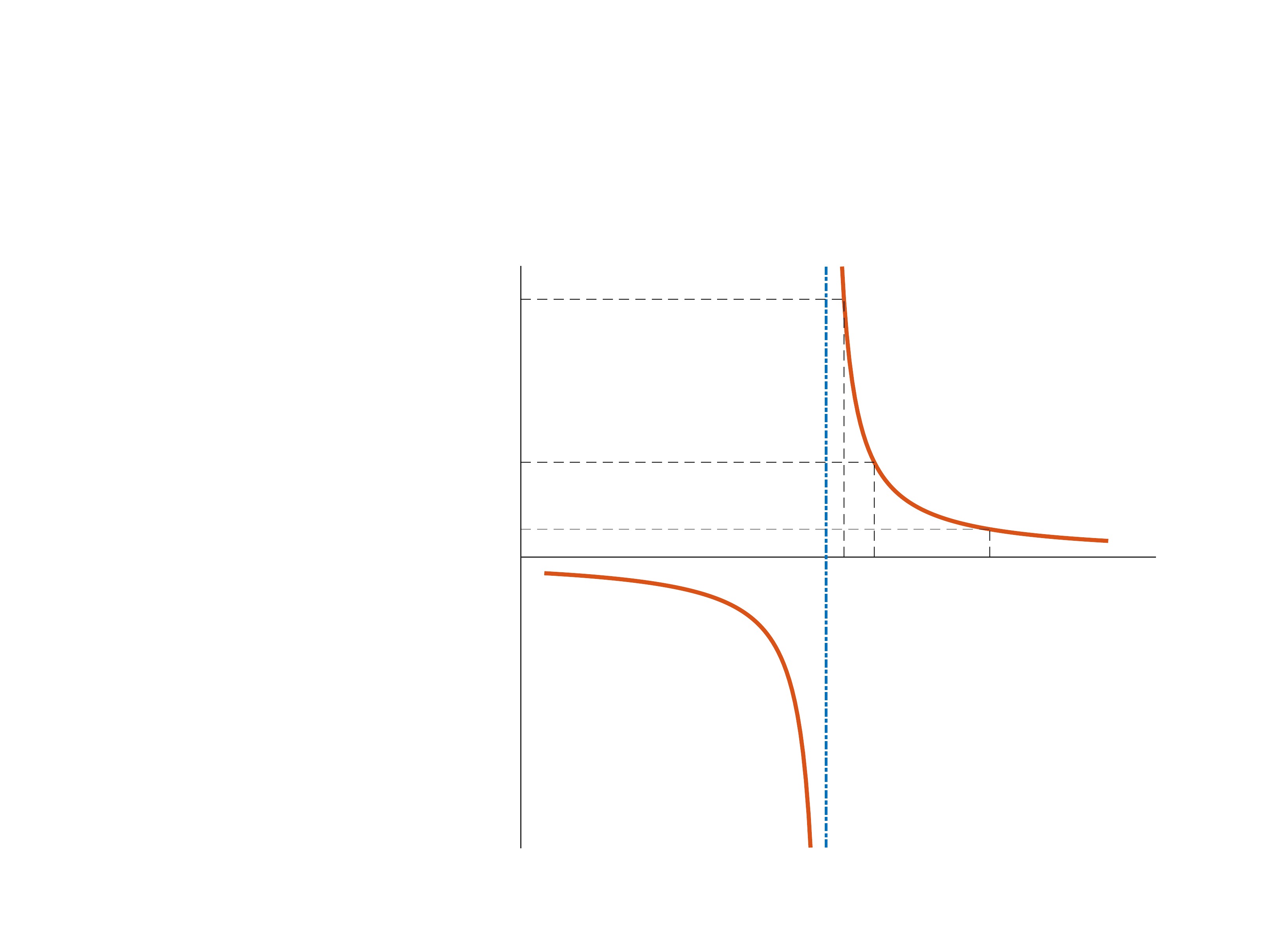}
		\put(93,38){\footnotesize $\lmh$}
		\put(50,38){\footnotesize $\lmh_1$}
		\put(55,38){\footnotesize $\lmh_2$}		
		\put(71,38){\footnotesize $\lmh_3$}
		\put(3,45){\footnotesize $\theta_3$}
		\put(3,54){\footnotesize $\theta_2$}
		\put(3,76){\footnotesize $\theta_1$}		
		\put(9,80){\footnotesize $\theta$}
		\put(51,20){\footnotesize $\lmh=\sg$}
		\put(56,65){\footnotesize $\theta=\dfrac{1}{\lmh-\sg} \quad (\sg<\lmh)$}
		\put(13,7){\footnotesize $\theta=\dfrac{1}{\lmh-\sg} \quad (\sg>\lmh)$}
	\end{overpic}
	\caption{The shift-and-invert spectral transformation: $\lmh$ values near the shift $\sg$ are well-separated on the $\theta$ axis.}
	\label{fig.ST}	
\end{figure}
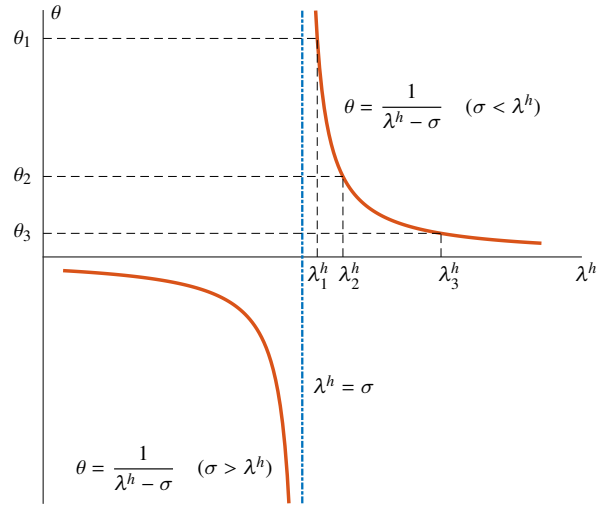

The eigenproblem~\eqref{eq.GHEP} defines a large sparse system of matrices with eigenvalues that may have arbitrary multiplicity.  Numerically, we seek to compute the the eigenpairs $\lmh_i$ and $\uuh_i$ within the given interval ${\lmh_i\in[\lm_s,\lm_e]}$, where either $\lm_s$ or $\lm_e$ can be infinite.  The conversion of a generalized eigenproblem to a standard one is fraught. The transformation can factor $\M$ (or $\K$) into its Cholesky decomposition as ${\M=\LL\LL^T}$ (or ${\K=\LL\LL^T}$) and solve ${\LL^{-1}\K\LL^{-T}\w=\lmh\.\w}$ (or ${\LL^{-1}\M\LL^{-T}\w=\w/\lmh}$), where ${\uuh=\LL^{-T}\w}$ (superscript $T$ refers to the transpose or conjugate transpose of a real symmetric or a Hermitian matrix, respectively). Either transformation may numerically fail.  This fragility of the computation process can have many causes. For example, $\K$ could be semidefinite, the eigenvalues may be insufficiently separated, or $\LL$ may be poorly conditioned. Any of these failings affect the extraction of eigenvectors from ${\uuh=\LL^{-T}\w}$ in the backward elimination~\cite{ Grimes1994}.   Thus, the eigensolvers generally perform a spectral transformation (ST) and solve the following shifted problem
\begin{equation} 
  (\K-\sg\M)\.\uuh=(\lmh-\sg)\.\M\uuh\,,
  \label{eq.shifted}
\end{equation}
to obtain an accurate approximation of eigenpairs (see, e.g.,~\cite{ NourOmid1987,Grimes1994,Demmel2000,XUE2011}).  Then, most eigensolvers use the shift-and-invert eigenproblem by solving the following system:
\begin{equation}
  (\K-\sg\M)^{-1}\M\uuh=\theta\.\uuh\,,
  \label{eq.ST_GHEP}
\end{equation}
where ${\theta=1/(\lmh-\sg)}$.  The resulting operator matrix ${(\K-\sg\M)^{-1}\M}$ is not symmetric, but self-adjoint with respect to $\M$.  Since matrices $\K$ and $\M$ have different null spaces, the ST matrix ${\K-\sg \M}$ is non-singular unless $\sg$ is an eigenvalue.  Hence, \eq{\eqref{eq.ST_GHEP}} resolves the (potential) issue of dealing with semidefinite system matrices.  Another main advantage behind using the shift-and-invert problem is to transform the eigenvalues $\lmh_i$ nearest the shift $\sg$ into the largest and well separated eigenvalues $\theta_i$ of the reciprocal eigenproblem of \eq{\eqref{eq.ST_GHEP}} (see \fig{\ref{fig.ST}}).  A well-selected shift enables the eigensolver to compute many eigenpairs in a single iteration.

In practical cases where the requested eigenvalue interval is large, we select additional shifts to prevent deterioration of the convergence rate of the eigensolution when the desired eigenvalues are far from $\sg$. An example is to find all critical speeds of a turbine shaft in a given working interval.  For this purpose, we select $\sg_k$s using a spectrum slicing technique (see, e.g.,~\cite{ Campos2012}). This method finds the requested eigenpairs with the true multiplicities while minimizing the computational effort (see \Sec{\ref{sub:Lanczos}}).

Apart from slicing the spectrum, Krylov eigensolvers also incorporate effective restarting mechanisms.  Restarting prevents an increase of the computational work needed for each shift in systems with significant numbers of degrees of freedom.  Well-known restarting techniques are the Sorensen's implicit restart~\cite{ Sorensen1992}, employed in the context of Lanczos methods, the thick-restart Lanczos method of Wu and Simon~\cite{ Wu2000}, and its unsymmetrical equivalent for the Arnoldi case referred to as the Krylov--Schur method of Stewart~\cite{ Stewart2002, Stewart2002a}.


\subsection{Lanczos method formulation for the generalized Hermitian eigenproblem} 
\label{sub:Lanczos}

Herein we focus on the implementation of the Lanczos method for solving the generalized Hermitian eigenproblems.  For each individual shift $\sg_k$, given the factorization ${(\K-\sg_k\M)=\LL\LL^T}$, we define the operator matrix ${\H:=\LL^{-T}\LL^{-1}\M}$, so that the eigenproblem of \eq{\eqref{eq.ST_GHEP}} becomes ${\H\uuh=\theta\.\uuh}$.  The $m$-step Lanczos decomposition consists of reducing the ${\N\times \N}$ matrix $\H$ to a symmetric tridiagonal matrix $\T_m$ ${(m\ll N)}$ as follows (see, e.g.,~\cite{ Demmel2000, Stewart2001}),
\begin{equation}
  \H\.\V_m=\V_m\T_m+\beta_m \vv_{m+1}\e_m^T\,,
  \label{eq.Lanczos}
\end{equation}
where $\e_m$ is the ${m\.}$th coordinate vector, and the term ${\beta_m \vv_{m+1}\e_m^T}$ is the residual of the $m$-step Lanczos decomposition.  In the above equation,
\begin{equation}
  \T_m=\left[
    \begin{array}{ccccc} \alpha_1 & \beta_2 \\
      \beta_2 & \alpha_2 & \beta_3 \\
                                  & \ddots & \ddots & \ddots \\
                                  & & \beta_{m-1} & \alpha_{m-1} & \beta_m \\
                                  & & & \beta_m & \alpha_m
    \end{array}
  \right],
  \label{eq.LanczosTm}
\end{equation}
and ${\V_m:=[\vv_1,\vv_2,...,\vv_m]}$ is the matrix of Lanczos vectors.  Assuming ${\beta_1=0}$ and $\vv_1$ is an initial \textit{generalized} unit vector of length $N$, i.e., ${\normM{\vv_1}=(\vv_1^T\M\vv_1)^{1/2}=1}$, the components of $\T_m$ and vectors $\vv_{j+1}$ $(j=1,2,...,m)$ are obtained by the following recurrence formulae,
\begin{equation}
  \begin{aligned}
    \alpha_j    &= \vv_j^T\M\.\H\vv_j \,,\\
    \beta_{j+1} &= \normM{\.\H\vv_j-\alpha_j\vv_j-\beta_j\vv_{j-1}} \,,\\
    \vv_{j+1}   &= \dfrac{1}{\beta_{j+1}} \pra{\H\vv_j-\alpha_j\vv_j-\beta_j\vv_{j-1}} \,.
  \end{aligned}
  \label{eq.LanczosAlphaBeta}
\end{equation}

In this way, the Lanczos vector $\vv_{m+1}$ is $\M$-orthogonal with respect to the columns of $\V_m$ in the Gram--Schmidt sense, resulting in ${\V_m^T\.\M\vv_{m+1}=0}$.  Hence, the $\M$-inner product of $\V_m$ premultiplied in \eq{\eqref{eq.Lanczos}} leads to the following equation, noting that $\V_m$ is an $\M$-orthogonal matrix (i.e., ${\V_m^T\.\M\.\V_m}=\textbf{I}$\.),
\begin{equation}
	\V_m^T\.\M\.\H\.\V_m=\T_m\,.
	\label{eq.LanczosRitz}
\end{equation}

The above equation reveals that $\T_m$ is the $\M$-orthogonal projection of $\H$ onto the ${m\.}$th order Krylov subspace $\kr_m(\H,\vv_1)$.
Therefore, if $\om_j$ and $\w_j$ are the eigenpairs of $\T_m$ (commonly referred to as Ritz values and Ritz vectors), the Rayleigh--Ritz approximation of the eigenpairs of $\H$ can be computed as
\begin{equation}
  \begin{aligned}
    \theta_j &= \om_j \,,\\
    \uuh_j   &= \V_m\w_j \,.
  \end{aligned}
  \label{eq.RitzValues}
\end{equation}

The eigenvalues of the original GHEP of \eq{\eqref{eq.GHEP}} for each Lanczos iteration is then obtained as ${\lmh_j=\sg_k+1/\theta_j}$ for each shift $\sg_k$.  Since $\T_m$ is a symmetric tridiagonal matrix, there exist multiple methods for computing its eigenpairs (see, e.g.,~\cite{ COAKLEY2013379}).

\begin{rem}
  By employing an LDL Cholesky factorization of the ST matrix, i.e., ${(\K-\sg_k\M)=\LL\DD\LL^T}$, where $\DD$ is a diagonal matrix, based on the Sylvester's law of inertia~\cite{ Parlett1998}, the number of eigenvalues smaller than $\sg_k$ is equal to the number of negative eigenvalues of $\DD$. Therefore, by defining ${\nu_{\.k}:=\nu\.(\K-\sg_k\M)}$ as the number of eigenvalues smaller than $\sg_k$, for two consecutive shifts $\sg_k$ and $\sg_{k+1}$, the interval ${[\sg_k,\sg_{k+1}]}$ has ${\nu_{\.k+1}-\nu_{\.k}}$ eigenvalues including their multiplicities. This rule drives the spectrum slicing technique when determining the required number of shifts.
  \label{rem.inertia}
\end{rem}

For each shift, only $m$ Lanczos steps are carried out followed by a restarting algorithm until computing all eigenpairs corresponding to the interval ${[\sg_k,\sg_{k+1}]}$.  After computing $\V_m$, a new Lanczos process starts, which benefits from the previously obtained spectral approximation.  The thick-restart approach presented in~\cite{ Wu2000} is an effective restarting technique in the case of Hermitian eigenproblems.  When restarting, the eigensolver keeps an appropriate number of Lanczos vectors, let say ${c<m}$.  The Lanczos recurrence of \eq{\eqref{eq.LanczosAlphaBeta}} continues after restarting with the following initial values:
\begin{equation}
	\begin{aligned}
		\rr_{c+1}	&= \H\vv_c-\alpha_c\vv_c-\sum_{i\.=1}^c \beta_i\vv_{i-1} \,,\\
		\beta_{c+1} &= \normM{\.\rr_{c+1}} \,,\\
		\vv_{c+1}   &= \dfrac{\rr_{c+1}}{\beta_{c+1}} \,,
	\end{aligned}
	\label{eq.LanczosRestart}
\end{equation}
where the orthogonalization is with respect to all previously stored Lanczos vectors.  Another points to consider are the deflation and spectrum recycling algorithms. The former is for maintaining the orthogonality of eigenvectors associated to a cluster of eigenvalues obtained from different shifts. The latter is to transform some Lanczos basis from a previous shift to the current one in case they create the same Krylov subspace.  More descriptions about the mathematical details of the above-mentioned algorithms can be found in, e.g.,~\cite{ Demmel2000, Campos2012, Wu2000, Olsson2006}.


\section{Cost estimation of the eigensolution} 
\label{sec.CostEstimate}

Estimating the cost of an eigensolution is challenging because it contains several numerical algorithms. 
In addition, different eigensolvers like, e.g., SLEPc~\cite{ slepc}, ANASAZI~\cite{ Hernandez2009}, ARPACK~\cite{ Lehoucq1998}, have their own methodologies for the eigensolution. 
In here, we focus on the formulation of Lanczos decomposition in \Sec{\ref{sub:Lanczos}}, 
for which we estimate its computational cost (measured in time)
based on the most expensive operations.


\subsection{Most expensive numerical operations} 
\label{sub:ExpensiveProcedures}

According to \eq{\eqref{eq.LanczosTm}}, constructing $\T_m$ involves $m$ computations of $\alpha_j$ and ${m-1}$ of $\beta_j$.
\eq{\eqref{eq.LanczosAlphaBeta}} shows that computing each $\alpha_j$ requires one forward/backward elimination of the Cholesky factors of the ST matrix, two multiplications of the mass matrix $\M$ by the respective vectors, and one vector product.
In the following, we refer to these operations as ``\fb~elimination", ``mat--vec" and ``vec--vec", respectively.
On the other hand, the computation of $\beta_j$ only needs one mat--vec followed by one vec--vec for the $\M$-norm calculation.
Additionally, extracting eigenvectors by \eq{\eqref{eq.RitzValues}} requires $\N_0$ multiplications of the ${\N\times m}$ Lanczos matrix $\V_m$ by the Ritz vectors.

To determine the cost of matrix factorization and \fb~elimination, we follow the theoretical estimates of 2D and 3D systems in terms of floating point operations (FLOPs) described in~\cite{ Garcia2017,Collier2012}. 
We have

\begin{equation}
	\begin{aligned}
 		\rm{IGA} &: & \FL_{\,\rm{fact}} &= \Or{\N^{\.(d+1)/2}\p^{\.3}} \,,\\
 		\rm{rIGA}&: & \FL_{\,\rm{fact}} &= 2^{\.d\l}\.\Or{(2^{-\.d\l}\N)^{\.(d+1)/2}\p^{\.3}}+\Or{\N^{\.(d+1)/2}} \,,
 	\end{aligned}
	\label{eq.flopsFact}
\end{equation}

\begin{equation}
	\begin{aligned}
 		\rm{IGA} &: & \FL_{\,\rm{\fb}} &= \Or{\N^{\.(d+1)/3}\p^{\.2}} \,,\\
 		\rm{rIGA}&: & \FL_{\,\rm{\fb}} &= 2^{\.d\l}\Or{(2^{-\.d\l}\N)^{\.(d+1)/3}\p^{\.2}}+\Or{\N^{\.(d+1)/3}} \,.
 	\end{aligned}
	\label{eq.flopsFb}
\end{equation}

Garcia el al.~\cite{ Garcia2017} show the factorization and \fb~elimination costs in large systems
reduces by up to $\O(\p^2)$ and $\O(\p)$, respectively, when using rIGA with $\BS$ of $2^{\.4}$ elements.
The cost of vec--vec products is proportional to the length of the vectors, $\N$, 
and the cost of mat--vecs is proportional to the number of nonzeros of the mass matrix, ${\nz\.(\M)}$.
In particular, the number of nonzeros of either mass or stiffness matrix is related to the sum of interactions that each basis function has with all other bases~\cite{ Collier2013}. 
As a result, referring to the matrix layouts of 1D systems in \fig{\ref{fig.matrixN8}} and the tensor-product property described in \Lem{\ref{lem.tensorproduct}}, 
one obtains the number of nonzeros of $\M$ as 

\begin{equation}
	\begin{aligned}
 		\rm{IGA} &: & \nz\.(\M) &= \bra{\ne(2\p+1)+\p^2}^{\.d} \,,\\
 		\rm{rIGA}&: & \nz\.(\M) &= 2^{\.d\l}\bra{2^{-\l}\ne(2\p+1)+\p^2-1}^{\.d} \,.
 	\end{aligned}
	\label{eq.nonzeros}
\end{equation}

Therefore, the cost of mat--vecs with IGA discretization is close to $\O(\N\p^{\.d})$ when $\N$ is sufficiently large, while that of rIGA discretizations is slightly higher.
Considering the optimal $\BS$ of $2^{\.4}$,
the ratio ${\nz(\.\M_{\rm{rIGA}})/\nz(\.\M_{\rm{IGA}})}$ is equal to $(1.03)^{\.d}$ and $(1.14)^{\.d}$,
for ${\p=2~\rm{and}~5}$, respectively, when $\ne$ is in the range of ${2^6 \sim 2^{11}}$.
The degradation incurred by rIGA, however, is not comparable to the improvements we obtain in the factorization and \fb~elimination steps by using rIGA.

For large systems, we identify the three most expensive operations, respectively, as
the Cholesky factorizations of the ST matrix, \fb~eliminations, and mat--vecs.
The costs of vec--vec products and extracting eigenvectors are significantly lower than that of mat--vecs, since ${\N\ll\nz\.(\M)}$ and, henceforth, we can exclude them from our cost estimates.
In addition, other numerical procedures in the eigenanalysis (e.g., system integration) are assumed to be of a lower order compared to the most expensive operations.
Hence, 
the total cost of the eigensolution is governed by the number of factorizations for each shift, $\Nfa$, the number of \fb~eliminations, $\Nfb$, and the number of mat--vecs, $\Nmv$.

\tab{\ref{tab:matvecs}} expresses these numbers in terms of number of shifts, $\Nsh$, and the total number of iterations, $\Nit$, 
carried out by the eigensolver.
The table also compares how an rIGA discretization improves or degrades the performance of each operation with respect to that of an IGA discretization.
To build this table, we assume the following:

\begin{itemize}
	\item IGA and rIGA discretizations use the same number of shifts. 
	This is derived from \Rem{\ref{rem.inertia}} and confirmed by numerical results (see \Sec{\ref{sec.Results}}).
	\item IGA and rIGA discretizations require the same number of iterations. 
	We show this numerically in \Sec{\ref{sec.Results}} for a sufficiently large number of eigenpairs (${\N_0\geq2^{10}}$).
	\item The number of Lanczos steps $m$ has the same average per shift through all iterations under both IGA and rIGA discretizations.
	Numerical results of \Sec{\ref{sec.Results}} confirm this assumption.
	\item The number of degrees of freedom is sufficiently large, so the cost improvements due to the use of rIGA described in~\cite{ Garcia2017} hold.
\end{itemize}

\begin{table}[!h]
	\centering
	\caption{Main numerical operations required by the Lanczos eigensolution algorithm
	and the improvement/degradation we obtain by using rIGA versus IGA discretization.
	We express the number of times we call each operation in terms of
	number of shifts, $\Nsh$, and the total number of iterations, $\Nit$.}
	\label{tab:matvecs}
	\small
	\begin{tabular}{@{}llll@{}}
		\toprule
		Numerical operation & 
		Matrix factorization & 
		\Fb~elimination & 
		Matrix--vector product  
		\\ \midrule
		\begin{tabular}[c]{@{}l@{}} Number of times \\ the operation is called \end{tabular} & 
		$\Nfa=\Nsh$ & 
		$\Nfb \approx m\.\Nit$ & 
		$\Nmv \approx 3m\.\Nit$  
		\\[10pt]
		\begin{tabular}[c]{@{}l@{}} Improvement/degradation of\\ performing one operation in rIGA\end{tabular} & 
		\begin{tabular}[c]{@{}l@{}} Improving by \\[2pt] $\O(\p^{\.2})$\end{tabular} & 
		\begin{tabular}[c]{@{}l@{}} Improving by \\[2pt] $\O(\p)$\end{tabular} & 
		\begin{tabular}[c]{@{}l@{}} Degrading by \\[2pt] $\nz(\.\M_{\rm{rIGA}})/\nz(\.\M_{\rm{IGA}})$\end{tabular}   
		\\ \bottomrule
	\end{tabular}
\end{table}


\subsection{Theoretical time estimates} 
\label{sub:TimeEstimates}

Referring to \eqz{\eqref{eq.flopsFact}}{\eqref{eq.nonzeros}}, we express the number of FLOPs of algorithms described in \tab{\ref{tab:matvecs}}
as $\O(\N^{\.a}\p^{\.b})$, where factors $a$ and $b$ vary for different operations and space dimensions as displayed in \tab{\ref{tab:constants}}.
Herein, we are interested in measuring the computational time. 
Since time and FLOPs are correlated for the type of operations considered here (as already shown in, e.g.,~\cite{ Garcia2017}), we estimate the time $\ti$ required to perform each operation as

\begin{equation}
	\ti \approx A \N^{\.a}\p^{\.b} \,,
	\label{eq.time}
\end{equation}

\noindent and in the logarithmic form as

\begin{equation}
	\log\ti \approx \log A + a\log\N + b\log\p \,.
	\label{eq.timelog}
\end{equation}

\begin{table}[!h]
	\centering
	\caption{Constants $a$ and $b$ in \eq{\eqref{eq.time}} for the theoretical estimation of computational time of the most expensive numerical operations of the Lanczos eigensolution algorithm. Constant $a$ is equal for both IGA and rIGA discretizations.}
	\label{tab:constants}
	\small
	\begin{tabular}{@{}lllll@{}}
		\toprule
		Constants & Discretization method & Matrix factorization & \Fb~elimination & Mat--vec product $^{*}$ \\ \midrule
		$a$ & IGA and rIGA & $(d+1)/\.2$ & $(d+1)/\.3$ & $1~$ \\
		\multirow{2}{*}{$b$} & IGA & 3 & 2 & $d$ \\
		 & rIGA & 1 & 1 & $d$ \\ \bottomrule
		\multicolumn{5}{l}{\footnotesize ${^*}$ The time of mat--vecs under rIGA increases by a factor in the range of $(1.03)^{\.d}$ and $(1.14)^{\.d}$ with ${\p=2\sim 5}$.}
	\end{tabular}
\end{table}

\begin{rem}
	When seeking for a sufficiently large number of eigenpairs, $\N_0$, 
	we assume the computational time grows linearly with respect to $\N_0$.
	Numerical results of \Sec{\ref{sec.Results}} confirm this assumption.
	\label{rem.time}
\end{rem}


\section{Implementation details} 
\label{sec.Implementation}

We discretize the model problem using PetIGA~\cite{ petiga}, a high-performance isogeometric analysis implementation based on PETSc (portable extensible toolkit for scientific computation)~\cite{ petsc}.
PetIGA has been utilized in many scientific and engineering applications (see, e.g.,~\cite{ Garcia2017,Garcia2018,Garcia2019,Collier2012,Collier2013,Vignal2015,Crtes2015,Espath2016,Espath2017}).
It allows us to investigate both IGA and rIGA discretizations on test cases with different numbers of elements in 2D and 3D, different polynomial degrees of the B-spline spaces, and different partitioning levels of the mesh.

We also use SLEPc, the scalable library for eigenvalue problem computations~\cite{ slepc}, for performing the eigenanalysis, allowing us to apply the shift-and-invert spectral transformation.
SLEPc, which has been used in solving different eigenproblems (see, e.g.,~\cite{ Campos2012,Romero2014,Campos2016,Faber2018,Keeli2018,AraujoC2020}),
employs the Krylov--Schur algorithm by default, which is equal to the thick-restart Lanczos algorithm in the case of generalized Hermitian eigenproblems.
SLEPc computes almost the same number of eigenpairs are computed for each shift, allowing us to efficiently estimate the computational costs.

We use multifrontal direct solver MUMPS~\cite{ mumps} to construct the LDL Cholesky factors, compute the required inertia for shift selections, and perform the forward/backward eliminations.
We employ the sequential version of MUMPS, which runs on a single thread (core).
We also use the automatic orderings provided by METIS~\cite{ Karypis1998}. 
For each test case, since all ST matrices have the same sparsity pattern, we only perform one {\em symbolic} factorization, followed by a certain number of numerical factorizations depending on the number of required shifts.
We executed all tests on a workstation
equipped with an Intel Xeon Gold 6230 CPU at 2.10 GHz
with 256 GB of RAM.


\section{Numerical results} 
\label{sec.Results}

We report the computational times (in seconds) required for finding the eigenpairs of the Laplace operator described in \Sec{\ref{sec.Preliminaries}}. 
We test different mesh sizes with ${\ne=2^s}$ elements in each direction, different partitioning levels of the rIGA discretization, namely ${\l=0~\mbox{(IGA)},1,2,...,s~\mbox{(FEA)}}$, and different polynomial degrees $\p$ of B-spline bases.
For 2D problems, we consider uniform meshes with ${s=8,9,10,11}$
and degree of ${\p=2,3,4,5}$. 
For the 3D case, we test on ${s=5,6,7}$ and the same degrees as in 2D.
To investigate the computational savings of the rIGA discretization compared to its IGA counterpart, we find the first $\N_0$ eigenpairs of the PDE system with different partitioning levels $\l$, where $\N_0$ can be as large as the number of eigenmodes of the IGA discretization $(\N_{\.\rm{IGA}})$.
We report the elapsed time for finding all $\N_0$ eigenpairs, $\tiN$,
the average required time per eigenmode, ${\tia:=\tiN/\N_0}$, 
and a normalization of time given by ${\tih:=\tiN/\N_0\.\N}$.

\begin{figure}[!h]
	\centering
	\begin{subfigure}{0.49\textwidth}\centering
		\begin{overpic}[width=1.0\textwidth,trim={0cm 0cm 0cm 0cm} ,clip]{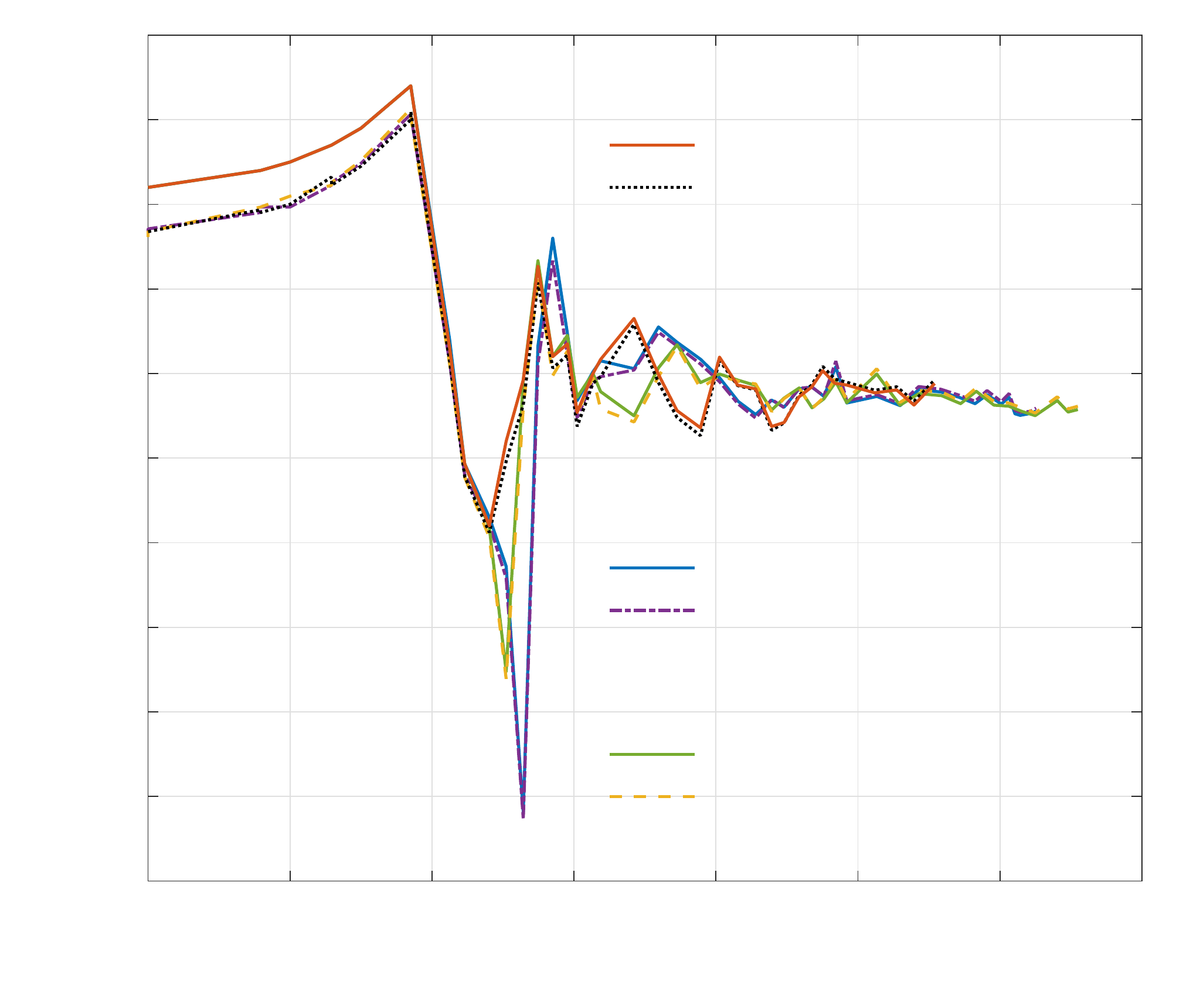}
			\put(32,2){\footnotesize No. of computed eigenpairs, $\N_0$}
			\put(11.00,6.50){\footnotesize $2^{0}$}
			\put(23.00,6.50){\footnotesize $2^{2}$}
			\put(35.00,6.50){\footnotesize $2^{4}$}
			\put(47.00,6.50){\footnotesize $2^{6}$}
			\put(59.00,6.50){\footnotesize $2^{8}$}
			\put(71.00,6.50){\footnotesize $2^{10}$}
			\put(83.00,6.50){\footnotesize $2^{12}$}
			\put(95.00,6.50){\footnotesize $2^{14}$}
			\put(7.75,10.00){\footnotesize 0}
			\put(6.75,17.10){\footnotesize 10}
			\put(6.75,24.20){\footnotesize 20}
			\put(6.75,31.30){\footnotesize 30}
			\put(6.75,38.40){\footnotesize 40}
			\put(6.75,45.50){\footnotesize 50}
			\put(6.75,52.60){\footnotesize 60}
			\put(6.75,59.70){\footnotesize 70}
			\put(6.75,66.80){\footnotesize 80}
			\put(6.75,73.90){\footnotesize 90}
			\put(5.75,81.00){\footnotesize 100}
			\put(60.00,36.50){\scriptsize  $\Nfb/\Nit$}
			\put(60.00,32.50){\scriptsize  $\Nmv/\.3\.\Nit$}
			\put(60.00,20.75){\scriptsize  $\Nfb/\Nit$}
			\put(60.00,16.75){\scriptsize  $\Nmv/\.3\.\Nit$}
			\put(60.00,72.25){\scriptsize  $\Nfb/\Nit$}
			\put(60.00,68.00){\scriptsize  $\Nmv/\.3\.\Nit$}
			\put(60.00,41.00){\scriptsize  $\text{mesh:~}512^{\.2},\,\p=2,\,\l=0$}
			\put(60.00,25.00){\scriptsize  $\text{mesh:~}256^{\.2},\,\p=4,\,\l=3$}
			\put(60.00,76.50){\scriptsize  $\text{mesh:~} 64^{\.3},\,\p=2,\,\l=2$}
			\begin{turn}{90}
			    \put(24,-3.5){\footnotesize Average No. of Lanczos steps, $m$}
			\end{turn}
		\end{overpic}
		\caption{}
		\label{fig.mconverge.a}
	\end{subfigure}
	\begin{subfigure}{0.49\textwidth}\centering
		\begin{overpic}[width=1.0\textwidth,trim={0cm 0cm 0cm 0cm} ,clip]{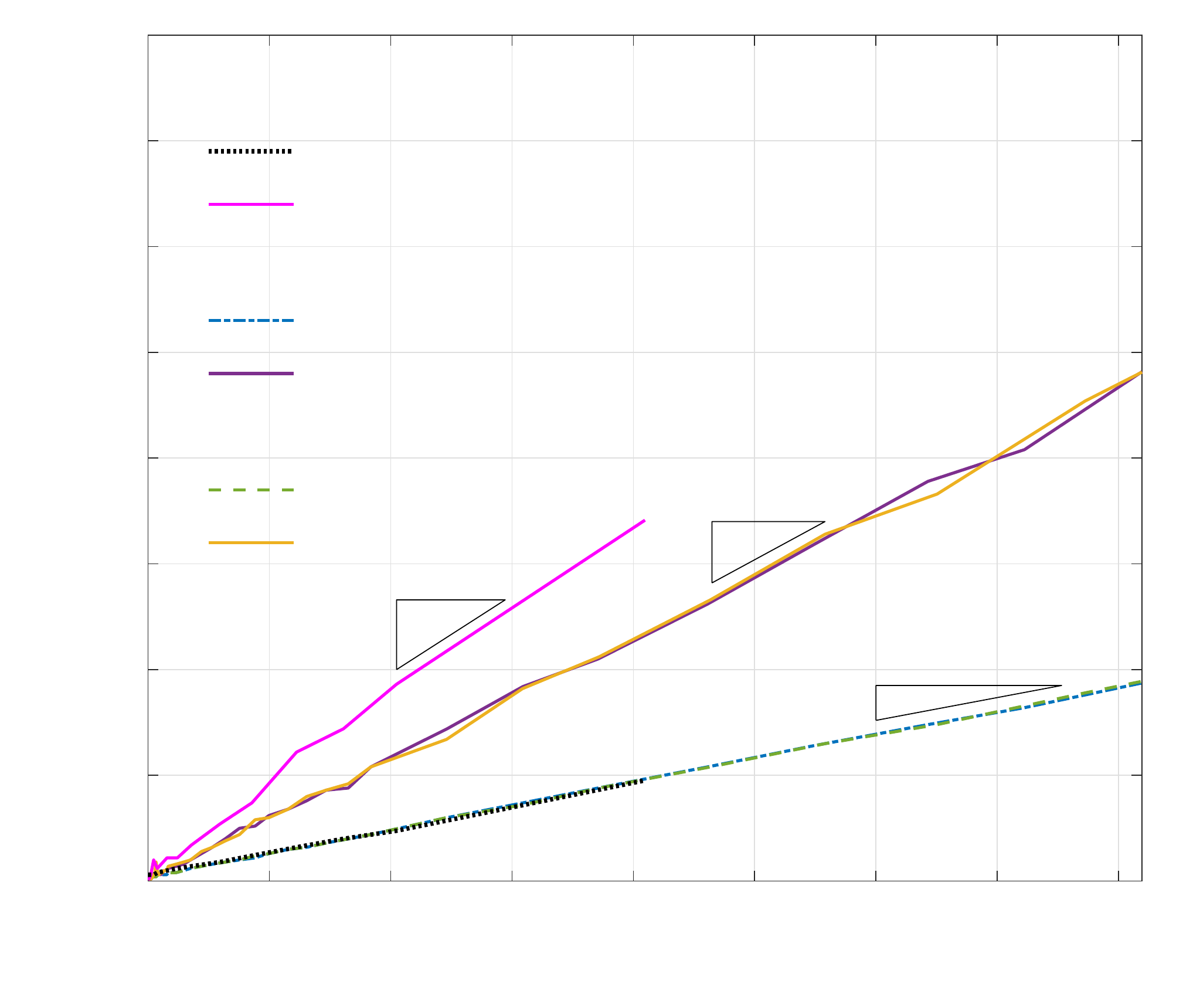}
			\put(32,2){\footnotesize No. of computed eigenpairs, $\N_0$}
			\put(11.50,6.50){\footnotesize $0$}
			\put(20.34,6.50){\footnotesize $500$}
			\put(29.00,6.50){\footnotesize $1000$}
			\put(39.53,6.50){\footnotesize $1500$}
			\put(50.00,6.50){\footnotesize $2000$}
			\put(60.72,6.50){\footnotesize $2500$}
			\put(70.81,6.50){\footnotesize $3000$}
			\put(80.91,6.50){\footnotesize $3500$}
			\put(91.00,6.50){\footnotesize $4000$}
			\put(7.75,10.00){\footnotesize $0$}
			\put(6.75,18.88){\footnotesize $50$}
			\put(5.75,27.75){\footnotesize $100$}
			\put(5.75,36.62){\footnotesize $150$}
			\put(5.75,45.50){\footnotesize $200$}
			\put(5.75,54.38){\footnotesize $250$}
			\put(5.75,63.25){\footnotesize $300$}
			\put(5.75,72.12){\footnotesize $350$}
			\put(5.75,81.00){\footnotesize $400$}
			\put(26.00,72.25){\scriptsize  $\Nsh=\Nfa$}
			\put(26.00,67.50){\scriptsize  $\Nit$}
			\put(26.00,58.00){\scriptsize  $\Nsh=\Nfa$}
			\put(26.00,53.50){\scriptsize  $\Nit$}
			\put(26.00,43.50){\scriptsize  $\Nsh=\Nfa$}
			\put(26.00,39.00){\scriptsize  $\Nit$}			
			\put(26.00,76.50){\scriptsize  $\text{mesh:~} 64^{\.3},\,\p=2,\,\l=2$}		
			\put(26.00,62.00){\scriptsize  $\text{mesh:~}512^{\.2},\,\p=2,\,\l=0$}
			\put(26.00,47.50){\scriptsize  $\text{mesh:~}256^{\.2},\,\p=4,\,\l=3$}
			\put(71.75,25){\scriptsize 1}
			\put(80,28.25){\scriptsize 43}
			\put(58,38){\scriptsize 1}
			\put(63,42){\scriptsize 16}
			\put(31,31){\scriptsize 1}
			\put(36,35.25){\scriptsize 12}
			\begin{turn}{90}
			    \put(35,-3.5){\footnotesize $\Nsh~~\rm{and}~~\Nit$}
			\end{turn}
		\end{overpic}
		\caption{}
		\label{fig.mconverge.b}
	\end{subfigure}
	\caption{(a) Independence of the average number of Lanczos steps, ${m\approx\Nfb/\Nit \approx\Nmv/3\.\Nit}$ from number of computed eigenpairs when $\N_0$ is large.
	(b) Linear increase of the numbers of shifts and iterations in terms of the number of computed eigenpairs. 
	We test on three different systems in terms of space dimension, domain size, polynomial degree and partitioning level.}
	\label{fig.mconverge}
\end{figure}
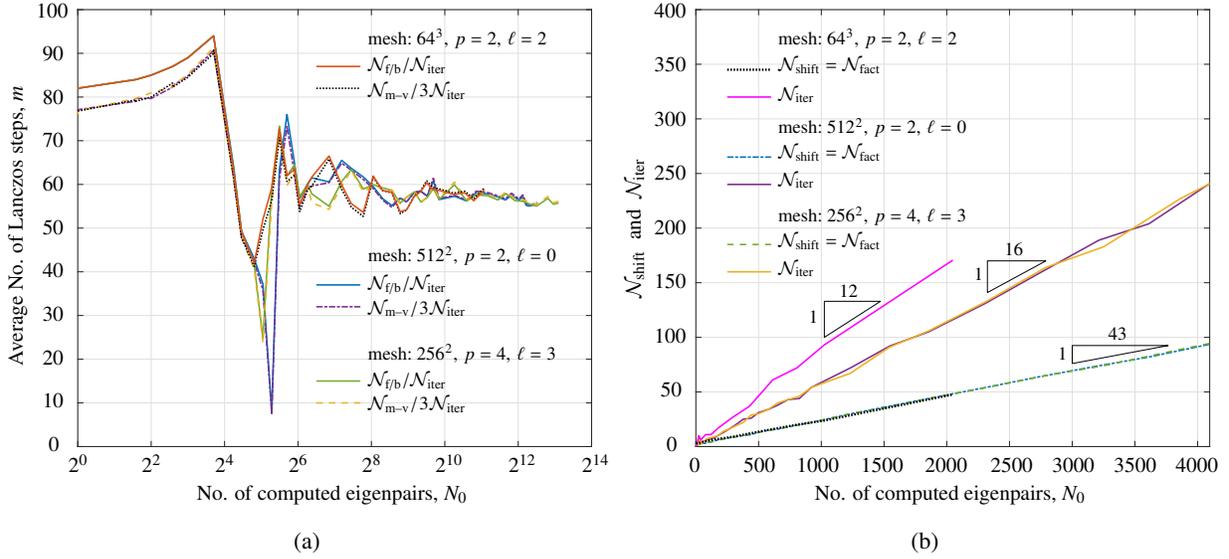

Before proceeding with time performance of IGA and rIGa discretizations in eigenanalysis,
we show in \fig{\ref{fig.mconverge.a}} that the number of Lanczos steps $m$ for a sufficiently large $\N_0$ becomes constant, independently of the mesh size and dimension. This confirms the assumptions of \Sec{\ref{sec.CostEstimate}}. 
Furthermore, \fig{\ref{fig.mconverge.b}} demonstrates that the number of shifts and the number of iterations increase linearly with $\N_0$,
indicating the proportional relationship of $\tiN$ and $\tia$ (see \Rem{\ref{rem.time}}).
The figure also shows that 3D systems need more iterations than 2D systems for finding the same number of eigenpairs.
These observations allow to predict expected times for large systems by solving only a small portion of the spectrum.

\subsection{Time performance of 2D eigenproblems} 
\label{sub:CostEvaluation2D}

\begin{figure}[!h]
	\centering
	\begin{subfigure}{0.49\textwidth}\centering
		\begin{overpic}[width=1.0\textwidth,trim={0cm 0cm 0cm 0cm} ,clip]{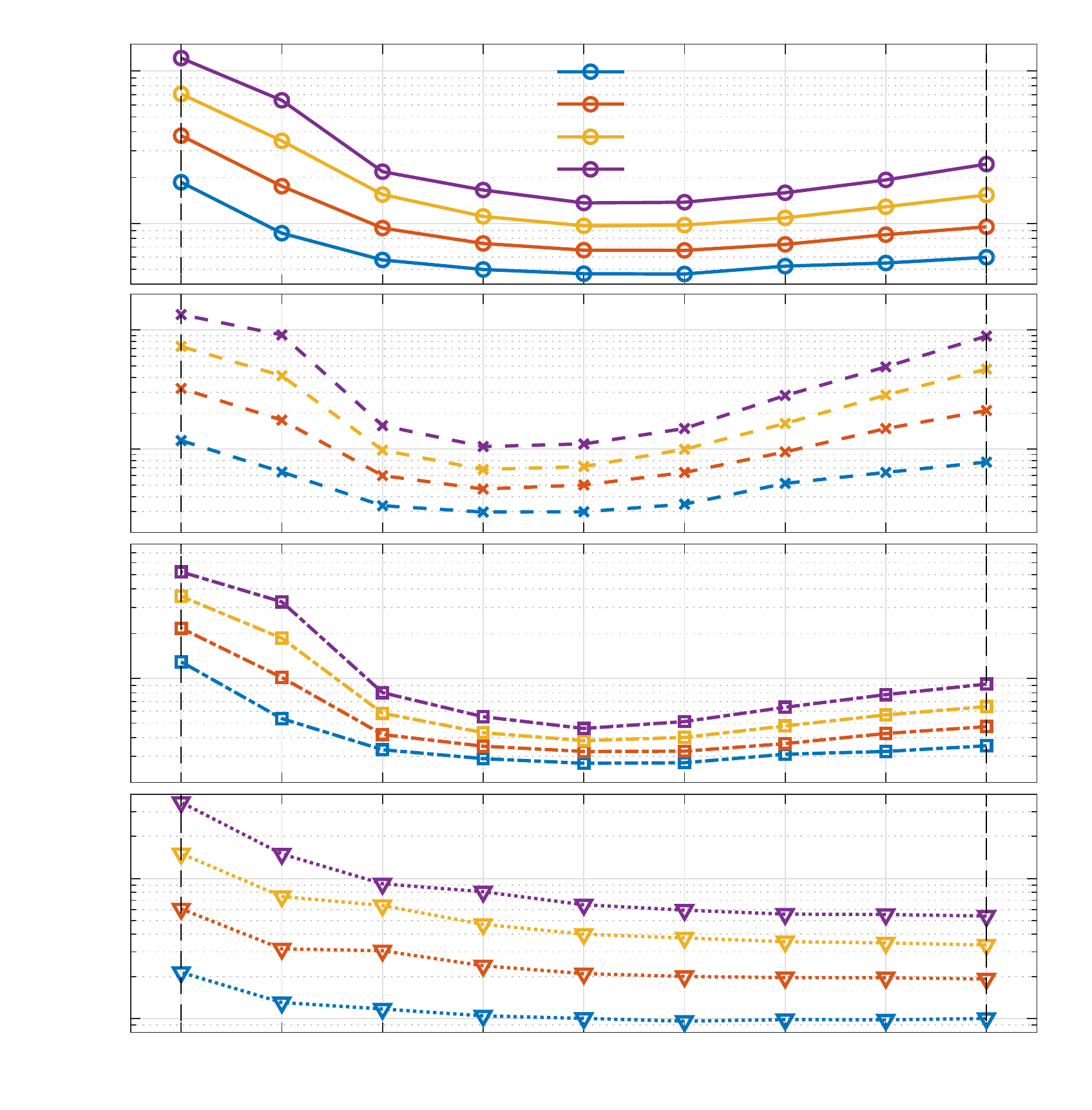}
			\put(46.5,0){\scriptsize  Blocksize}
			\put(57,84){\scriptsize  $\p=5$}
			\put(57,87){\scriptsize  $\p=4$}
			\put(57,90){\scriptsize  $\p=3$}
			\put(57,93){\scriptsize  $\p=2$}
			\put(17,93){\line(1,2){2.25}}
			\put(88.5,93){\line(-1,2){2.25}}
			\put(17,98){\scriptsize  $C^{\.0}\,(\rm{FEA})$}
			\put(81,98){\scriptsize  $C^{\p-1}\,(\rm{IGA})$}
			\put(14.50,3.50){\scriptsize $2^{\.0}$}
			\put(23.69,3.50){\scriptsize $2^{\.1}$}
			\put(32.88,3.50){\scriptsize $2^{\.2}$}
			\put(42.06,3.50){\scriptsize $2^{\.3}$}
			\put(51.25,3.50){\scriptsize $2^{\.4}$}
			\put(60.44,3.50){\scriptsize $2^{\.5}$}
			\put(69.62,3.50){\scriptsize $2^{\.6}$}
			\put(78.81,3.50){\scriptsize $2^{\.7}$}
			\put(88.00,3.50){\scriptsize $2^{\.8}$}
			\put(6.5,92.00){\scriptsize $10^{\,4}$}
			\put(6.5,78.00){\scriptsize $10^{\,3}$}
			\put(6.5,69.00){\scriptsize $10^{\,3}$}
			\put(6.5,58.00){\scriptsize $10^{\,2}$}
			\put(6.5,38.00){\scriptsize $10^{\,3}$}
			\put(6.5,19.50){\scriptsize $10^{\,3}$}
			\put(6.5,7.00){\scriptsize $10^{\,2}$}
			\begin{turn}{90}
				\put(9,-0.5){\scriptsize Mat--vec time}
			    \put(15,-3.5){\scriptsize (sec)}
			    \put(30.5,-0.5){\scriptsize \Fb~elimination}
			    \put(34,-3.5){\scriptsize time (sec)}
			    \put(55,-0.5){\scriptsize Factorization}
			    \put(57,-3.5){\scriptsize time (sec)}
			    \put(79,-0.5){\scriptsize Total time}
			    \put(82,-3.5){\scriptsize (sec)}
			\end{turn}
		\end{overpic}
		\caption{Mesh: $256^{\,2}$}
	\end{subfigure}
	\begin{subfigure}{0.49\textwidth}\centering
		\begin{overpic}[width=1.0\textwidth,trim={0cm 0cm 0cm 0cm} ,clip]{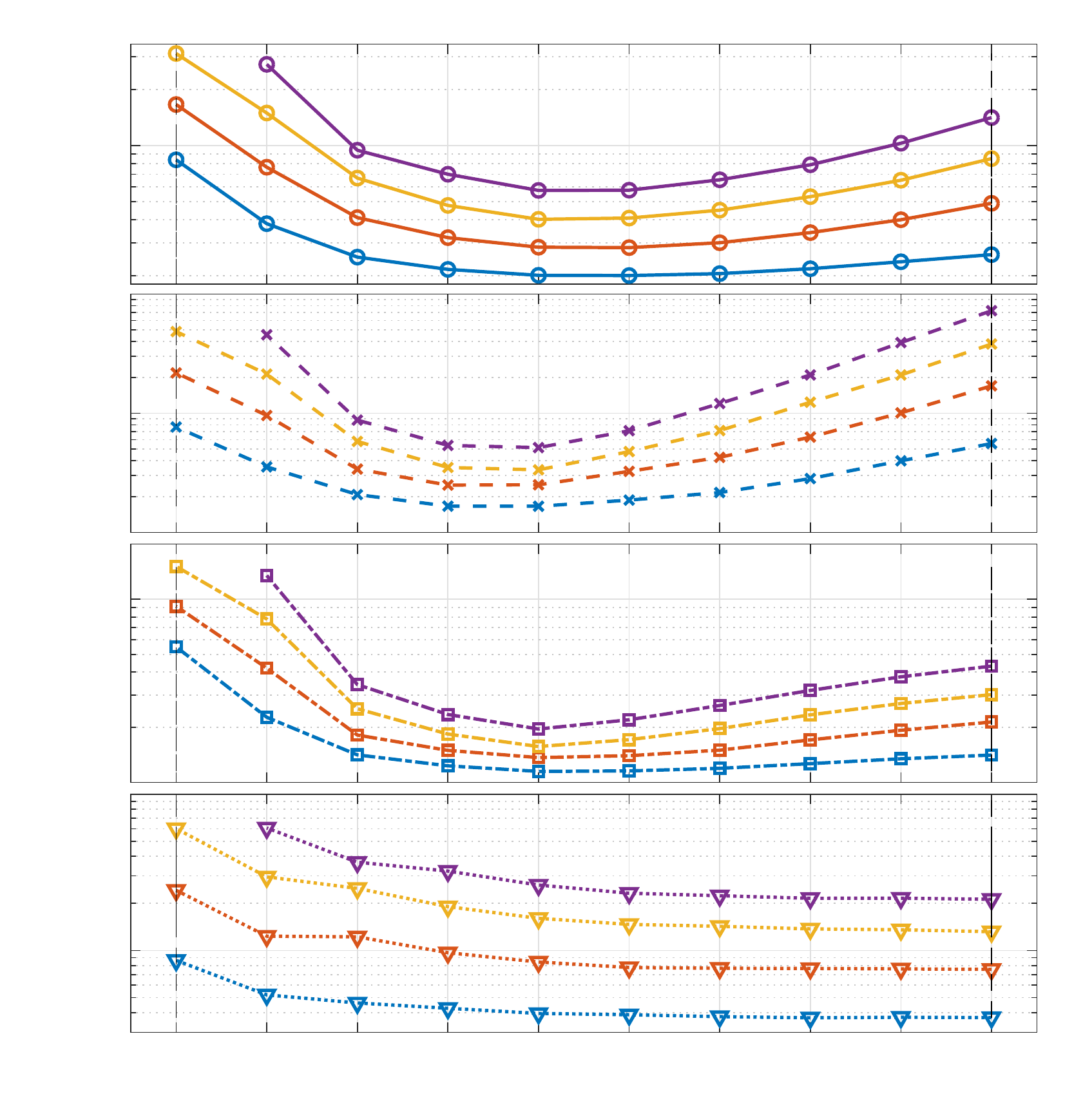}
			\put(46.5,0){\scriptsize  Blocksize}
			\put(17,93){\line(1,2){2.25}}
			\put(88.5,93){\line(-1,2){2.25}}
			\put(17,98){\scriptsize  $C^{\.0}\,(\rm{FEA})$}
			\put(81,98){\scriptsize  $C^{\p-1}\,(\rm{IGA})$}
			\put(14.50,3.50){\scriptsize $2^{\.0}$}
			\put(22.67,3.50){\scriptsize $2^{\.1}$}
			\put(30.83,3.50){\scriptsize $2^{\.2}$}
			\put(39.00,3.50){\scriptsize $2^{\.3}$}
			\put(47.17,3.50){\scriptsize $2^{\.4}$}
			\put(55.33,3.50){\scriptsize $2^{\.5}$}
			\put(63.50,3.50){\scriptsize $2^{\.6}$}
			\put(71.67,3.50){\scriptsize $2^{\.7}$}
			\put(79.83,3.50){\scriptsize $2^{\.8}$}
			\put(88.00,3.50){\scriptsize $2^{\.9}$}
			\put(6.5,86.00){\scriptsize $10^{\,4}$}
			\put(6.5,71.50){\scriptsize $10^{\,4}$}
			\put(6.5,61.00){\scriptsize $10^{\,3}$}
			\put(6.5,52.00){\scriptsize $10^{\,2}$}
			\put(6.5,45.00){\scriptsize $10^{\,4}$}
			\put(6.5,29.00){\scriptsize $10^{\,3}$}
			\put(6.5,13.00){\scriptsize $10^{\,3}$}
			\begin{turn}{90}
			    \put(9,-0.5){\scriptsize Mat--vec time}
			    \put(15,-3.5){\scriptsize (sec)}
			    \put(30.5,-0.5){\scriptsize \Fb~elimination}
			    \put(34,-3.5){\scriptsize time (sec)}
			    \put(55,-0.5){\scriptsize Factorization}
			    \put(57,-3.5){\scriptsize time (sec)}
			    \put(79,-0.5){\scriptsize Total time}
			    \put(82,-3.5){\scriptsize (sec)}
			\end{turn}
		\end{overpic}
		\caption{Mesh: $512^{\,2}$}
	\end{subfigure}
	\begin{subfigure}{0.49\textwidth}\centering
		\begin{overpic}[width=1.0\textwidth,trim={0cm 0cm 0cm 0cm} ,clip]{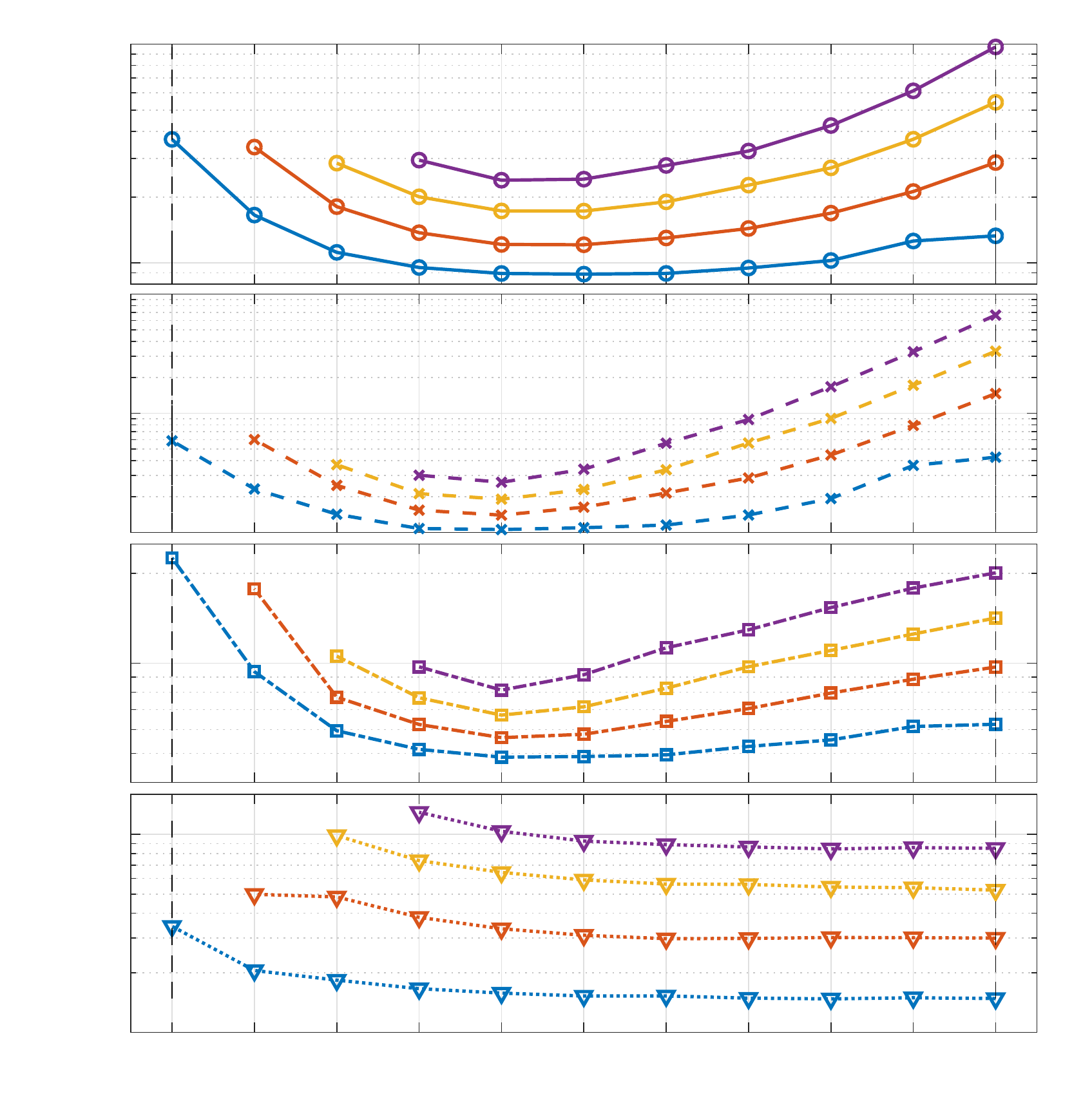}
			\put(46.5,0){\scriptsize  Blocksize}
			\put(16,90){\line(2,1){5}}
			\put(88.5,93){\line(-1,2){2.25}}
			\put(22,92){\scriptsize  $C^{\.0}\,(\rm{FEA})$}
			\put(81,98){\scriptsize  $C^{\p-1}\,(\rm{IGA})$}
			\put(14.50,3.50){\scriptsize $2^{\.0}$}
			\put(21.85,3.50){\scriptsize $2^{\.1}$}
			\put(29.20,3.50){\scriptsize $2^{\.2}$}
			\put(36.55,3.50){\scriptsize $2^{\.3}$}
			\put(43.90,3.50){\scriptsize $2^{\.4}$}
			\put(51.25,3.50){\scriptsize $2^{\.5}$}
			\put(58.60,3.50){\scriptsize $2^{\.6}$}
			\put(65.95,3.50){\scriptsize $2^{\.7}$}
			\put(73.30,3.50){\scriptsize $2^{\.8}$}
			\put(80.65,3.50){\scriptsize $2^{\.9}$}
			\put(88.00,3.50){\scriptsize $2^{\.10}$}
			\put(6.5,94.00){\scriptsize $10^{\,5}$}
			\put(6.5,76.00){\scriptsize $10^{\,4}$}
			\put(6.5,71.50){\scriptsize $10^{\,5}$}
			\put(6.5,61.50){\scriptsize $10^{\,4}$}
			\put(6.5,52.00){\scriptsize $10^{\,3}$}
			\put(6.5,39.00){\scriptsize $10^{\,4}$}
			\put(6.5,24.00){\scriptsize $10^{\,4}$}
			\put(6.5,7.00){\scriptsize $10^{\,3}$}
			\begin{turn}{90}
			    \put(9,-0.5){\scriptsize Mat--vec time}
			    \put(15,-3.5){\scriptsize (sec)}
			    \put(30.5,-0.5){\scriptsize \Fb~elimination}
			    \put(34,-3.5){\scriptsize time (sec)}
			    \put(55,-0.5){\scriptsize Factorization}
			    \put(57,-3.5){\scriptsize time (sec)}
			    \put(79,-0.5){\scriptsize Total time}
			    \put(82,-3.5){\scriptsize (sec)}
			\end{turn}
		\end{overpic}
		\caption{Mesh: $1024^{\,2}$}
	\end{subfigure}
	\begin{subfigure}{0.49\textwidth}\centering
		\begin{overpic}[width=1.0\textwidth,trim={0cm 0cm 0cm 0cm} ,clip]{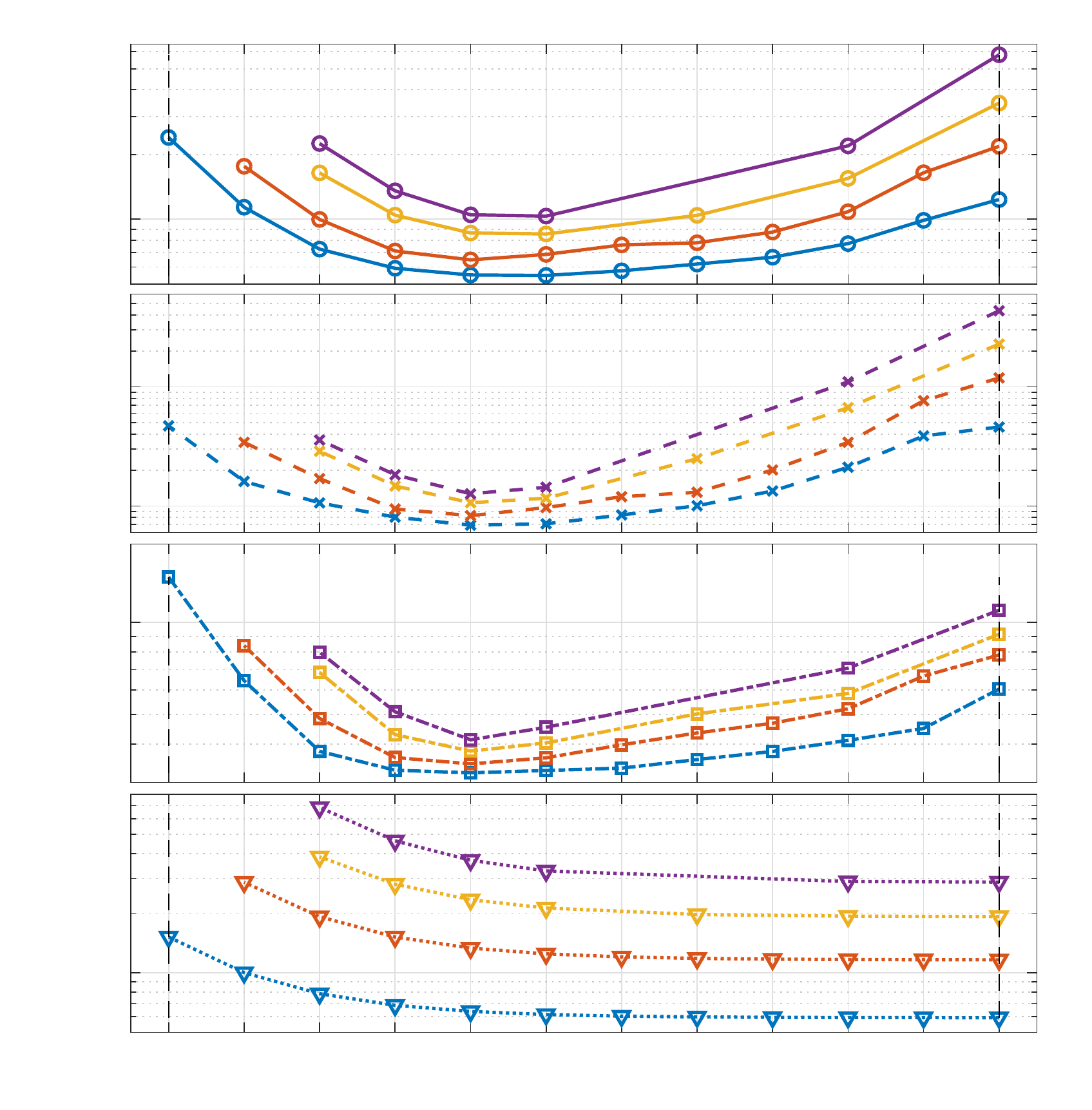}
			\put(46.5,0){\scriptsize  Blocksize}
			\put(16,90){\line(2,1){5}}
			\put(88.5,93){\line(-1,2){2.25}}
			\put(22,92){\scriptsize  $C^{\.0}\,(\rm{FEA})$}
			\put(81,98){\scriptsize  $C^{\p-1}\,(\rm{IGA})$}
			\put(14.50,3.50){\scriptsize $2^{\.0}$}
			\put(21.18,3.50){\scriptsize $2^{\.1}$}
			\put(27.86,3.50){\scriptsize $2^{\.2}$}
			\put(34.55,3.50){\scriptsize $2^{\.3}$}
			\put(41.23,3.50){\scriptsize $2^{\.4}$}
			\put(47.91,3.50){\scriptsize $2^{\.5}$}
			\put(54.59,3.50){\scriptsize $2^{\.6}$}
			\put(61.27,3.50){\scriptsize $2^{\.7}$}
			\put(67.95,3.50){\scriptsize $2^{\.8}$}
			\put(74.64,3.50){\scriptsize $2^{\.9}$}
			\put(81.32,3.50){\scriptsize $2^{\.10}$}
			\put(88.00,3.50){\scriptsize $2^{\.11}$}
			\put(6.5,79.00){\scriptsize $10^{\,5}$}
			\put(6.5,64.00){\scriptsize $10^{\,5}$}
			\put(6.5,53.50){\scriptsize $10^{\,4}$}
			\put(6.5,43.00){\scriptsize $10^{\,5}$}
			\put(6.5,11.50){\scriptsize $10^{\,4}$}
			\begin{turn}{90}
			    \put(9,-0.5){\scriptsize Mat--vec time}
			    \put(15,-3.5){\scriptsize (sec)}
			    \put(30.5,-0.5){\scriptsize \Fb~elimination}
			    \put(34,-3.5){\scriptsize time (sec)}
			    \put(55,-0.5){\scriptsize Factorization}
			    \put(57,-3.5){\scriptsize time (sec)}
			    \put(79,-0.5){\scriptsize Total time}
			    \put(82,-3.5){\scriptsize (sec)}
			\end{turn}
		\end{overpic}
		\caption{Mesh: $2048^{\,2}$}
		\label{fig.timevsBS.2048}
	\end{subfigure}
	\caption{The total times and those of the most expensive numerical operations for finding the first ${\N_0=4096}$ eigenpairs of the 2D Laplace operator with different mesh sizes of ${\ne=2^s}$ ${(s=8,9,10,11)}$ and degrees ${\p=2,3,4,5}$.
	We test rIGA discretizations with different $\BS$s $(2^{s-\.\l})$ obtained using different levels of partitioning ${\l=0,1,...,s}$.}
	\label{fig.timevsBS}
\end{figure}
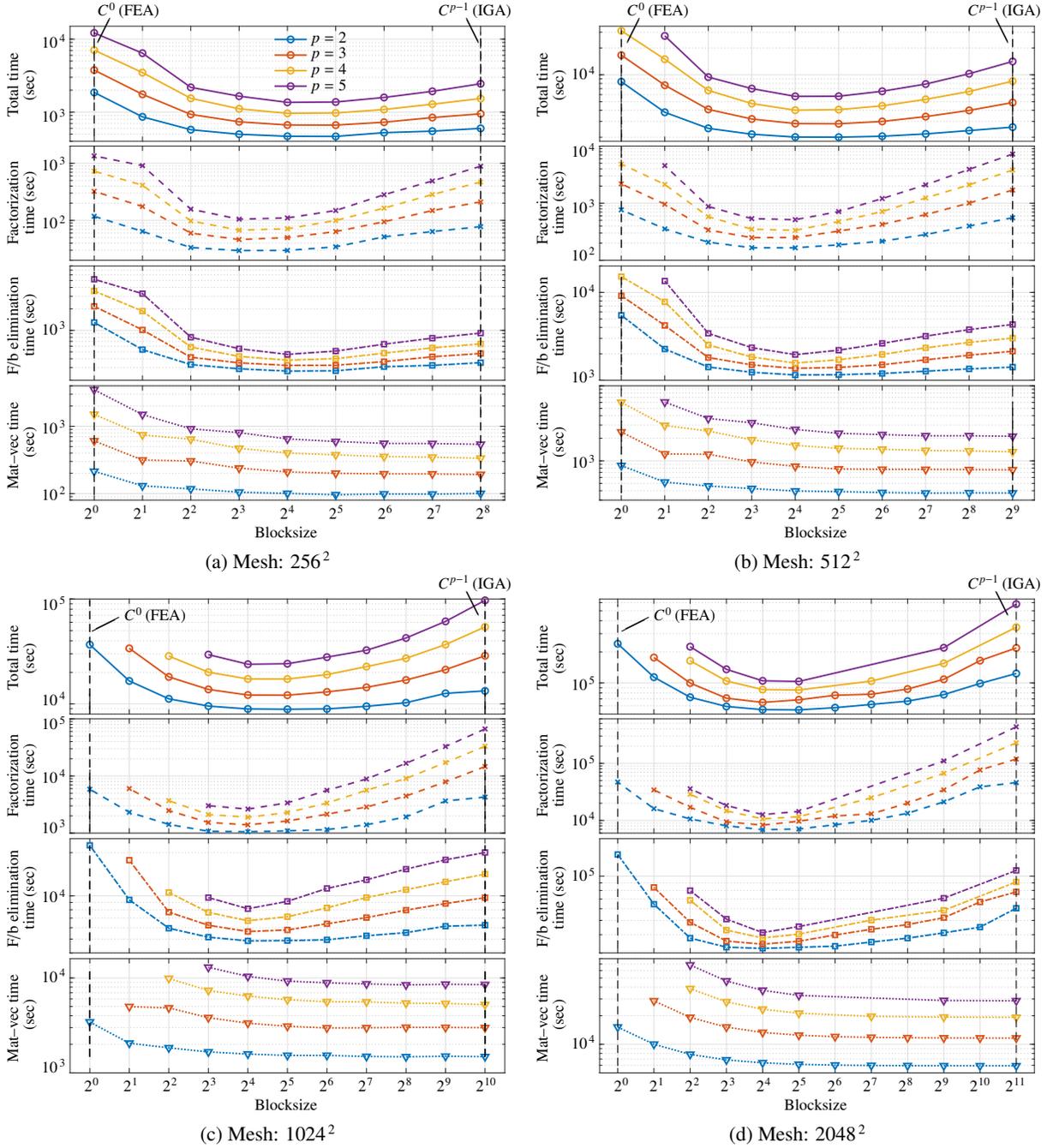

The partitioning level $\l$ of rIGA affects the computational time of the eigenanalysis. 
To see this, we focus on the most expensive numerical operations described in \Sec{\ref{sub:ExpensiveProcedures}}. 
We monitor the total elapsed time for different $\BS$s considering the partitioning scheme presented in \Sec{\ref{sub:rIGApartitioning}}. 
There is an optimal $\BS$ of ${2^{s-\.\l}}$ at which the factorization and \fb~elimination times are minimum; however,
the time of mat--vec products always increases with $\l$.
\fig{\ref{fig.timevsBS}} shows the computational times of finding ${\N_0=2^{12}}$ eigenpairs (i.e., $\tiN$) as a function of the $\BS$.
The rIGA factorization time reaches a minimum at a ${\BS}$ of 16 elements almost in all cases, which coincides with the optimal $\BS$ for the f/b~elimination time. 
Hence, 
we obtain the maximum savings for the total elapsed time (considering all numerical operations) by employing macroelements of size 16.

\begin{table}[!h]
\centering
\caption{The average computational times per eigenvalue, ${\tia=\tiN/\N_0}$, and the number of executions of the most expensive operations for finding ${\N_0=4096}$ eigenpairs of the 2D test cases illustrated in \fig{\ref{fig.timevsBS}}.
For rIGA discretizations, we select a $\BS$ of 16 elements.}
\label{tab:T1_2D}
\footnotesize
\begin{tabular}{@{}llllllllllllll@{}}
\toprule
\multicolumn{3}{l}{Discretization}  & \multicolumn{3}{l}{Factorization} & \multicolumn{3}{l}{\Fb~elimination} & \multicolumn{3}{l}{Mat--vec product} & \multicolumn{2}{l}{Total} \\ 
\cmidrule(r){1-3} \cmidrule(lr){4-6} \cmidrule(lr){7-9} \cmidrule(lr){10-12} \cmidrule(l){13-14} 
Mesh & $\p$ & Method & $\Nfa$ &
\begin{tabular}[c]{@{}l@{}}$\tiax{\text{fact}}$\\(sec)\end{tabular} & 
\begin{tabular}[c]{@{}l@{}}Improved\\by\end{tabular} & $\Nfb$ & \begin{tabular}[c]{@{}l@{}}$\tiax{\fb}$\\(sec)\end{tabular} & 
\begin{tabular}[c]{@{}l@{}}Improved\\by\end{tabular} & $\Nmv$ & \begin{tabular}[c]{@{}l@{}}$\tiax{\text{m--v}}$\\(sec)\end{tabular} &
\begin{tabular}[c]{@{}l@{}}Degraded\\by\end{tabular} & \begin{tabular}[c]{@{}l@{}}$\tia$\\(sec)\end{tabular} & 
\begin{tabular}[c]{@{}l@{}}Improved\\by\end{tabular} \\ \midrule
\multirow{8}{*}{$256^{\,2}$} 
 & \multirow{2}{*}{2} & IGA & 96 & 0.019 & \multirow{2}{*}{2.612} & 13672 & 0.086 & \multirow{2}{*}{1.311} & 41143 & 0.024 & \multirow{2}{*}{0.998} & 0.146 & \multirow{2}{*}{1.282} \\ 
				  &  & rIGA & 96 & 0.007 &						  & 13538 & 0.065 &						   & 40705 & 0.024 &  & 0.114 &  \\[2pt] 
 & \multirow{2}{*}{3} & IGA & 96 & 0.051 & \multirow{2}{*}{4.218} & 14401 & 0.115 & \multirow{2}{*}{1.469} & 43459 & 0.046 & \multirow{2}{*}{0.917} & 0.232 & \multirow{2}{*}{1.424} \\ 
  				  &  & rIGA & 96 & 0.012 &						  & 14026 & 0.078 &						   & 42378 & 0.051 &  & 0.163 &  \\[2pt] 
 & \multirow{2}{*}{4} & IGA & 92 & 0.114 & \multirow{2}{*}{6.581} & 14637 & 0.157 & \multirow{2}{*}{1.696} & 44109 & 0.081 & \multirow{2}{*}{0.837} & 0.375 & \multirow{2}{*}{1.593} \\ 
  				  &  & rIGA & 96 & 0.017 &						  & 14905 & 0.093 &						   & 45318 & 0.097 &  & 0.235 &  \\[2pt] 
 & \multirow{2}{*}{5} & IGA & 96 & 0.217 & \multirow{2}{*}{8.048} & 15050 & 0.223 & \multirow{2}{*}{1.989} & 45511 & 0.131 & \multirow{2}{*}{0.830} & 0.598 & \multirow{2}{*}{1.797} \\ 
 				  &  & rIGA & 94 & 0.026 &						  & 14878 & 0.112 &						   & 44958 & 0.158 &  & 0.332 &  \\ \cmidrule(){1-14}
\multirow{8}{*}{$512^{\,2}$} 
 & \multirow{2}{*}{2} & IGA & 96 & 0.135 & \multirow{2}{*}{3.361} & 14127 & 0.344 & \multirow{2}{*}{1.228} & 42693 & 0.091 & \multirow{2}{*}{0.940} & 0.632 & \multirow{2}{*}{1.291} \\ 
 			 	  &  & rIGA & 96 & 0.040 &						  & 13495 & 0.280 &						   & 40570 & 0.096 &  & 0.489 &  \\[2pt] 
 & \multirow{2}{*}{3} & IGA & 93 & 0.414 & \multirow{2}{*}{6.771} & 14980 & 0.521 & \multirow{2}{*}{1.564} & 45344 & 0.185 & \multirow{2}{*}{0.897} & 1.193 & \multirow{2}{*}{1.721} \\ 
 			 	  &  & rIGA & 96 & 0.061 &						  & 13881 & 0.333 &						   & 41715 & 0.206 &  & 0.693 &  \\[2pt] 
 & \multirow{2}{*}{4} & IGA & 96 & 0.931 & \multirow{2}{*}{11.39} & 14647 & 0.735 & \multirow{2}{*}{1.920} & 44237 & 0.322 & \multirow{2}{*}{0.821} & 2.073 & \multirow{2}{*}{2.117} \\ 
    		      &  & rIGA & 92 & 0.081 &						  & 15162 & 0.383 &						   & 45762 & 0.392 &  & 0.979 &  \\[2pt] 
 & \multirow{2}{*}{5} & IGA & 96 & 1.772 & \multirow{2}{*}{14.09} & 14449 & 1.053 & \multirow{2}{*}{2.205} & 43401 & 0.519 & \multirow{2}{*}{0.810} & 3.447 & \multirow{2}{*}{2.462} \\ 
  				  &  & rIGA & 95 & 0.125 &						  & 14465 & 0.477 &						   & 43465 & 0.640 &  & 1.399 &  \\ \cmidrule(){1-14}
\multirow{8}{*}{$1024^{\,2}$}
 & \multirow{2}{*}{2} & IGA & 96 & 1.044 & \multirow{2}{*}{4.050} & 14642 & 1.526 & \multirow{2}{*}{1.287} & 44475 & 0.361 & \multirow{2}{*}{0.938} & 3.246 & \multirow{2}{*}{1.486} \\ 
			  	  &  & rIGA & 96 & 0.257 &						  & 14717 & 1.185 &						   & 44688 & 0.385 &  & 2.183 &  \\[2pt] 
 & \multirow{2}{*}{3} & IGA & 96 & 3.572 & \multirow{2}{*}{10.48} & 14503 & 2.370 & \multirow{2}{*}{1.721} & 43768 & 0.730 & \multirow{2}{*}{0.897} & 7.022 & \multirow{2}{*}{2.369} \\ 
  				  &  & rIGA & 94 & 0.340 &						  & 14551 & 1.377 &						   & 43742 & 0.814 &  & 2.964 &  \\[2pt] 
 & \multirow{2}{*}{4} & IGA & 96 & 8.101 & \multirow{2}{*}{17.46} & 14942 & 3.457 & \multirow{2}{*}{2.110} & 45267 & 1.277 & \multirow{2}{*}{0.813} & 13.23 & \multirow{2}{*}{3.137} \\ 
  				  &  & rIGA & 95 & 0.463 &						  & 15259 & 1.638 &						   & 46212 & 1.570 &  & 4.218 &  \\[2pt] 
 & \multirow{2}{*}{5} & IGA & 96 & 16.26 & \multirow{2}{*}{25.28} & 14245 & 4.894 & \multirow{2}{*}{2.464} & 42963 & 2.078 & \multirow{2}{*}{0.820} & 23.71 & \multirow{2}{*}{4.067} \\ 
  				  &  & rIGA & 96 & 0.643 &						  & 14294 & 1.985 &						   & 43047 & 2.534 &  & 5.830 &  \\ \cmidrule(){1-14}
\multirow{8}{*}{$2048^{\,2}$}
 & \multirow{2}{*}{2} & IGA & 96 & 11.23 & \multirow{2}{*}{6.671} & 14526 & 14.77 & \multirow{2}{*}{1.878} & 43701 & 1.447 & \multirow{2}{*}{0.930} & 30.14 & \multirow{2}{*}{2.240} \\ 
  				  &  & rIGA & 95 & 1.683 &						  & 14563 & 7.868 &						   & 43954 & 1.557 &  & 13.45 &  \\[2pt] 
 & \multirow{2}{*}{3} & IGA & 96 & 29.02 & \multirow{2}{*}{14.33} & 14351 & 19.09 & \multirow{2}{*}{2.269} & 42979 & 2.838 & \multirow{2}{*}{0.872} & 53.27 & \multirow{2}{*}{3.369} \\ 
  				  &  & rIGA & 96 & 2.024 &						  & 14566 & 8.415 &						   & 43976 & 3.254 &  & 15.81 &  \\[2pt] 
 & \multirow{2}{*}{4} & IGA & 94 & 55.73 & \multirow{2}{*}{21.43} & 14453 & 22.31 & \multirow{2}{*}{2.408} & 43549 & 4.694 & \multirow{2}{*}{0.822} & 84.54 & \multirow{2}{*}{4.011} \\ 
  				  &  & rIGA & 96 & 2.599 &						  & 14240 & 9.263 &						   & 42728 & 5.713 &  & 21.08 &  \\[2pt] 
 & \multirow{2}{*}{5} & IGA & 96 & 106.1 & \multirow{2}{*}{34.38} & 14312 & 26.73 & \multirow{2}{*}{2.650} & 43178 & 7.015 & \multirow{2}{*}{0.776} & 141.7 & \multirow{2}{*}{5.536} \\ 
  				  &  & rIGA & 96 & 3.088 &						  & 14419 & 10.08 &						   & 43526 & 9.040 &  & 25.59 &  \\ \bottomrule
\end{tabular}
\end{table}

When a matrix is sufficiently large, its Cholesky factorization is more expensive than \fb~elimination and mat--vec multiplication.
However, Krylov eigensolvers perform multiple \fb~eliminations and mat--vecs per Cholesky factorization. 
In our 2D case, ${\Nfb\approx 150~\Nfa}$ and ${\Nmv\approx 450~\Nfa}$ (see \tab{\ref{tab:matvecs}} and \fig{\ref{fig.mconverge}}).
It brings the costs of these two operations in a comparable range with the matrix factorization.
To compare results, \tab{\ref{tab:T1_2D}} reports the number of executions of each operation for finding ${\N_0=2^{12}}$ eigenpairs.
We also report the average computational times per eigenvalue 
obtained by dividing the time of each numerical procedure by $\N_0$ (i.e., ${\tia=\tiN/\N_0}$). 
Results indicate an improvement in the cost of matrix factorization close to $\O(\p^{\.2})$ for large problems, and of almost $\O(\p)$ for \fb~eliminations.
We observe a slight degradation in cost of mat--vec multiplications due to the increase of nonzero terms of the mass matrix under rIGA.
In summary, the total observed time saving for the entire eigensolution is up to $\O(\p)$ for large domains.
The time of finding ${\N_0=2^{\.12}}$ discrete eigenpairs with ${\ne=2048}$ and quintic basis functions reduces from 161 hours to 29 hours using the optimal rIGA discretization (see \fig{\ref{fig.timevsBS.2048}}).
If the domain is very large, the total computational cost is governed only by matrix factorization.
Therefore, we predict the total time improvements of up to $\O(\p^2)$.
To observe them, we would need larger computational resources.
On the other hand, for small problems (e.g., 2D systems with ${\ne\leq 256}$ and ${\p\leq 3}$),
the cost of IGA is comparable (or smaller) to that of rIGA. This occurs because the cost of matrix factorization becomes a small fraction of the total cost.

\begin{figure}[!h]
	\centering
	\begin{subfigure}{0.49\textwidth}\centering
		\begin{overpic}[width=1.0\textwidth,trim={0cm 0cm 0cm 0cm} ,clip]{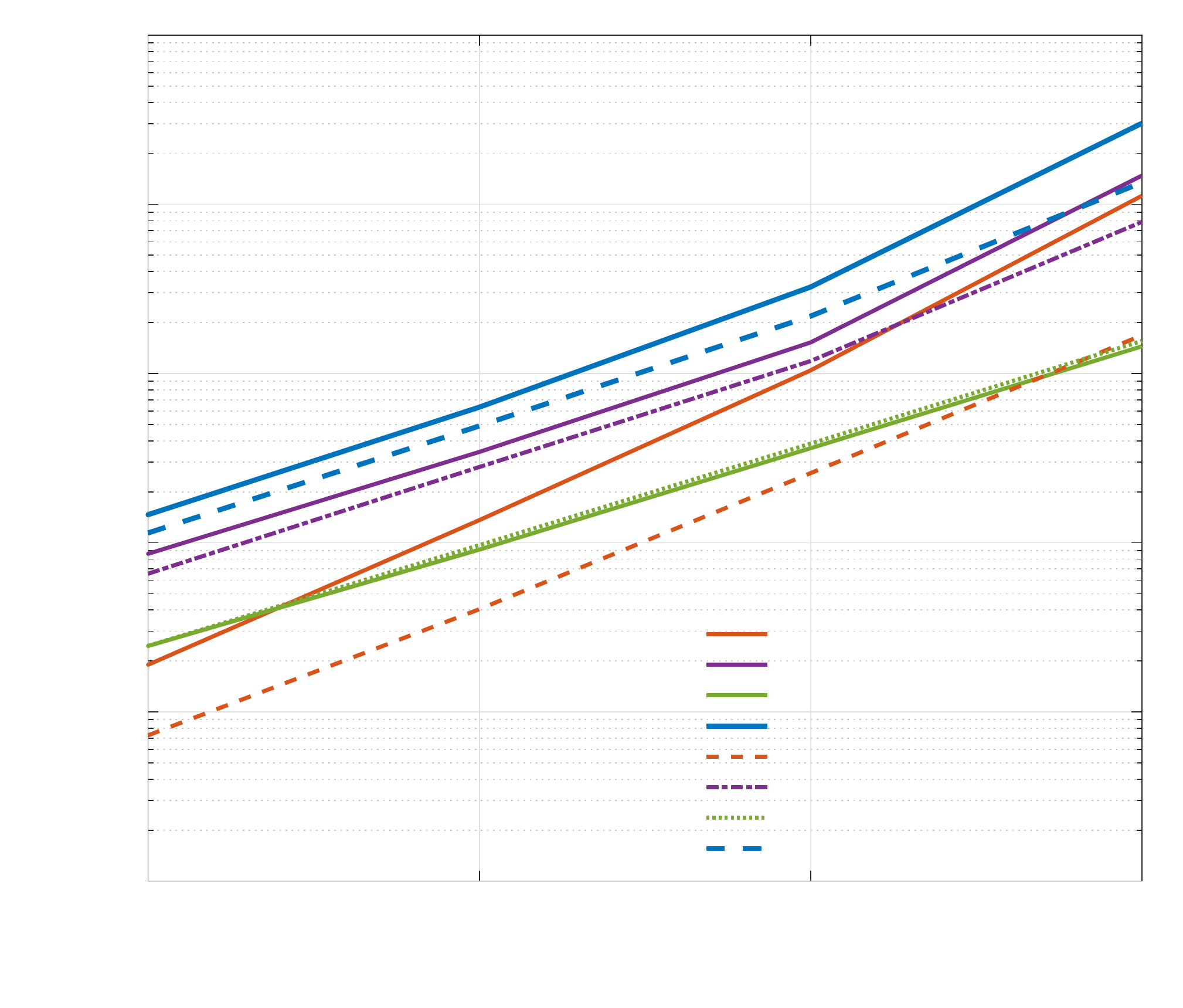}
			\put(30,2){\footnotesize No. of elements in each direction, $\ne$}
			\put(11.00,6.50){\footnotesize $2^{8}$}
			\put(38.67,6.50){\footnotesize $2^{9}$}
			\put(66.33,6.50){\footnotesize $2^{10}$}
			\put(94.00,6.50){\footnotesize $2^{11}$}
			\put(5.50,10.00){\footnotesize $10^{-3}$}
			\put(5.50,24.20){\footnotesize $10^{-2}$}
			\put(5.50,38.40){\footnotesize $10^{-1}$}
			\put(5.75,52.60){\footnotesize $10^{\,0}$}
			\put(5.75,66.80){\footnotesize $10^{\,1}$}
			\put(5.75,81.00){\footnotesize $10^{\,2}$}
			\put(66.00,31.50){\scriptsize  IGA, factorization}
			\put(66.00,28.83){\scriptsize  IGA, \fb~elimination}
			\put(66.00,26.16){\scriptsize  IGA, mat--vec product}
			\put(66.00,23.49){\scriptsize  IGA, total time}
			\put(66.00,20.81){\scriptsize rIGA, factorization}
			\put(66.00,18.14){\scriptsize rIGA, \fb~elimination}
			\put(66.00,15.47){\scriptsize rIGA, mat--vec product}
			\put(66.00,12.80){\scriptsize rIGA, total time}
			\begin{turn}{90}
			    \put(26,-3.5){\footnotesize Time per eigenvalue, $\tia$ (sec)}
			\end{turn}
		\end{overpic}
		\caption{$\p=2$}
	\end{subfigure}
	\begin{subfigure}{0.49\textwidth}\centering
		\begin{overpic}[width=1.0\textwidth,trim={0cm 0cm 0cm 0cm} ,clip]{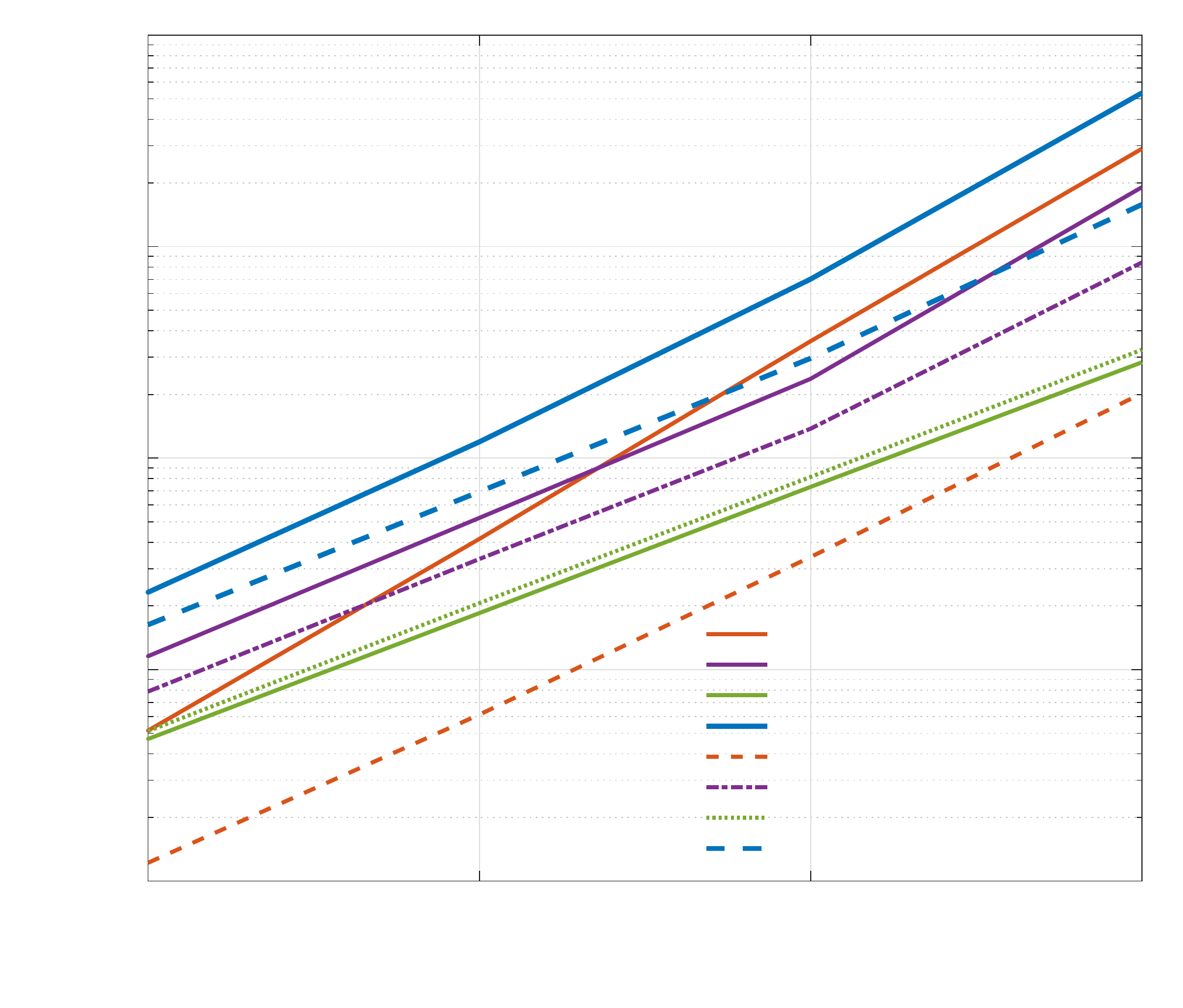}
			\put(30,2){\footnotesize No. of elements in each direction, $\ne$}
			\put(11.00,6.50){\footnotesize $2^{8}$}
			\put(38.67,6.50){\footnotesize $2^{9}$}
			\put(66.33,6.50){\footnotesize $2^{10}$}
			\put(94.00,6.50){\footnotesize $2^{11}$}
			\put(5.50,10.00){\footnotesize $10^{-2}$}
			\put(5.50,27.75){\footnotesize $10^{-1}$}
			\put(5.75,45.50){\footnotesize $10^{\,0}$}
			\put(5.75,63.25){\footnotesize $10^{\,1}$}
			\put(5.75,81.00){\footnotesize $10^{\,2}$}
			\put(66.00,31.50){\scriptsize  IGA, factorization}
			\put(66.00,28.83){\scriptsize  IGA, \fb~elimination}
			\put(66.00,26.16){\scriptsize  IGA, mat--vec product}
			\put(66.00,23.49){\scriptsize  IGA, total time}
			\put(66.00,20.81){\scriptsize rIGA, factorization}
			\put(66.00,18.14){\scriptsize rIGA, \fb~elimination}
			\put(66.00,15.47){\scriptsize rIGA, mat--vec product}
			\put(66.00,12.80){\scriptsize rIGA, total time}
			\begin{turn}{90}
			    \put(26,-3.5){\footnotesize Time per eigenvalue, $\tia$ (sec)}
			\end{turn}
		\end{overpic}
		\caption{$\p=3$}
	\end{subfigure}
	\begin{subfigure}{0.49\textwidth}\centering
		\begin{overpic}[width=1.0\textwidth,trim={0cm 0cm 0cm 0cm} ,clip]{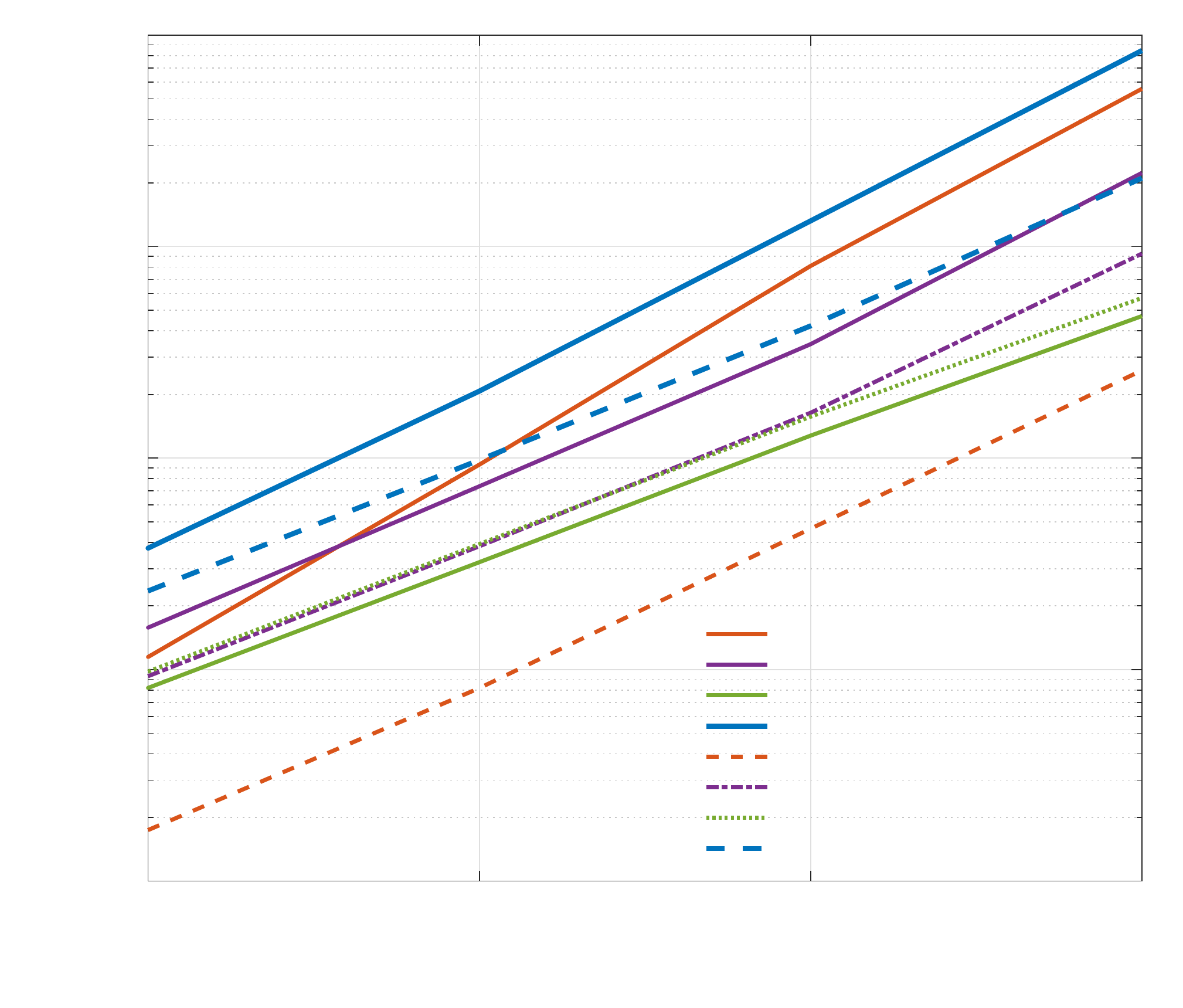}
			\put(30,2){\footnotesize No. of elements in each direction, $\ne$}
			\put(11.00,6.50){\footnotesize $2^{8}$}
			\put(38.67,6.50){\footnotesize $2^{9}$}
			\put(66.33,6.50){\footnotesize $2^{10}$}
			\put(94.00,6.50){\footnotesize $2^{11}$}
			\put(5.50,10.00){\footnotesize $10^{-2}$}
			\put(5.50,27.75){\footnotesize $10^{-1}$}
			\put(5.75,45.50){\footnotesize $10^{\,0}$}
			\put(5.75,63.25){\footnotesize $10^{\,1}$}
			\put(5.75,81.00){\footnotesize $10^{\,2}$}
			\put(66.00,31.50){\scriptsize  IGA, factorization}
			\put(66.00,28.83){\scriptsize  IGA, \fb~elimination}
			\put(66.00,26.16){\scriptsize  IGA, mat--vec product}
			\put(66.00,23.49){\scriptsize  IGA, total time}
			\put(66.00,20.81){\scriptsize rIGA, factorization}
			\put(66.00,18.14){\scriptsize rIGA, \fb~elimination}
			\put(66.00,15.47){\scriptsize rIGA, mat--vec product}
			\put(66.00,12.80){\scriptsize rIGA, total time}
			\begin{turn}{90}
			    \put(26,-3.5){\footnotesize Time per eigenvalue, $\tia$ (sec)}
			\end{turn}
		\end{overpic}
		\caption{$\p=4$}
	\end{subfigure}
	\begin{subfigure}{0.49\textwidth}\centering
		\begin{overpic}[width=1.0\textwidth,trim={0cm 0cm 0cm 0cm} ,clip]{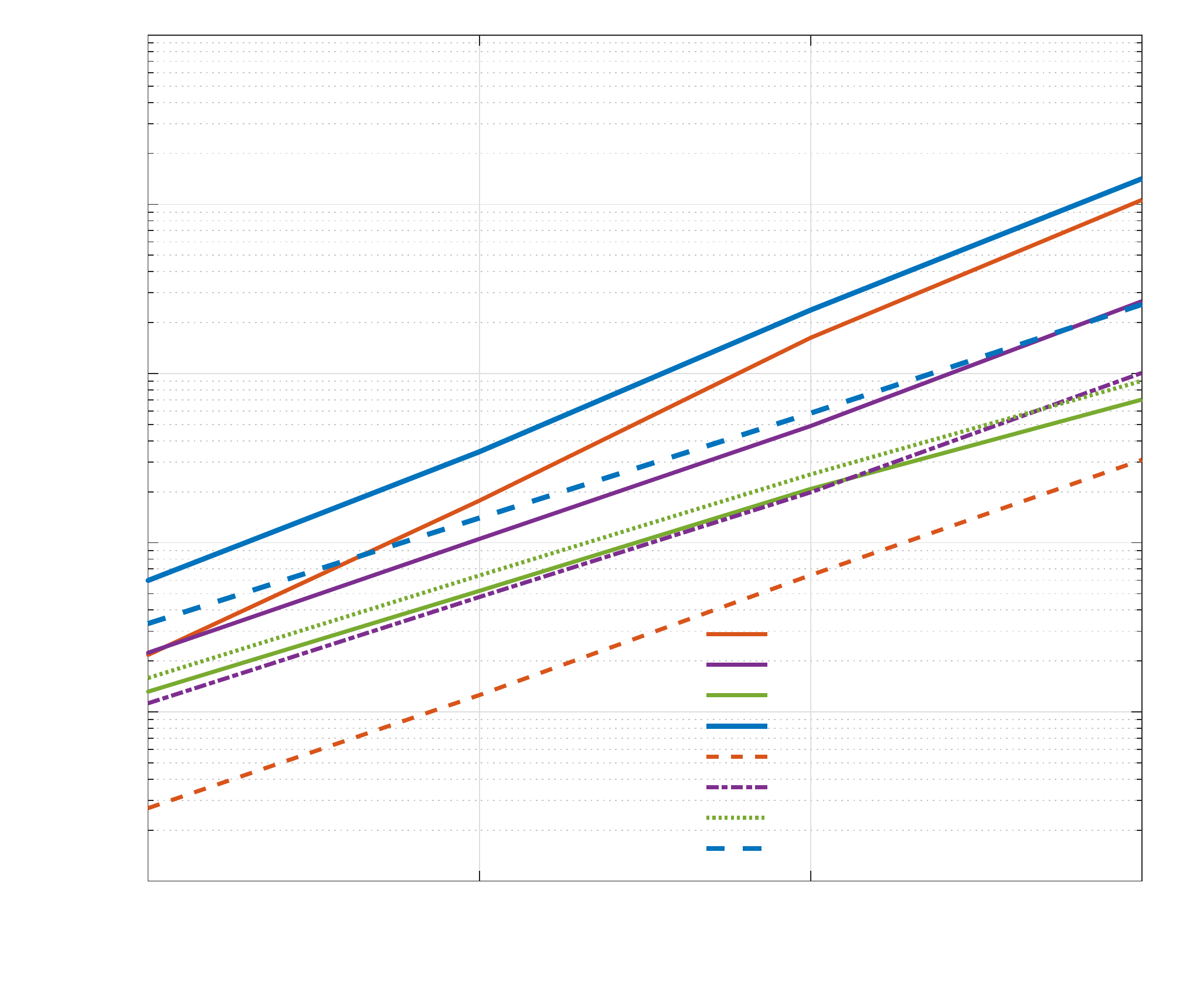}
			\put(30,2){\footnotesize No. of elements in each direction, $\ne$}
			\put(11.00,6.50){\footnotesize $2^{8}$}
			\put(38.67,6.50){\footnotesize $2^{9}$}
			\put(66.33,6.50){\footnotesize $2^{10}$}
			\put(94.00,6.50){\footnotesize $2^{11}$}
			\put(5.50,10.00){\footnotesize $10^{-2}$}
			\put(5.50,24.20){\footnotesize $10^{-1}$}
			\put(5.75,38.40){\footnotesize $10^{\,0}$}
			\put(5.75,52.60){\footnotesize $10^{\,1}$}
			\put(5.75,66.80){\footnotesize $10^{\,2}$}
			\put(5.75,81.00){\footnotesize $10^{\,3}$}
			\put(66.00,31.50){\scriptsize  IGA, factorization}
			\put(66.00,28.83){\scriptsize  IGA, \fb~elimination}
			\put(66.00,26.16){\scriptsize  IGA, mat--vec product}
			\put(66.00,23.49){\scriptsize  IGA, total time}
			\put(66.00,20.81){\scriptsize rIGA, factorization}
			\put(66.00,18.14){\scriptsize rIGA, \fb~elimination}
			\put(66.00,15.47){\scriptsize rIGA, mat--vec product}
			\put(66.00,12.80){\scriptsize rIGA, total time}
			\begin{turn}{90}
			    \put(26,-3.5){\footnotesize Time per eigenvalue, $\tia$ (sec)}
			\end{turn}
		\end{overpic}
		\caption{$\p=5$}
	\end{subfigure}
	\caption{Average time per eigenvalue, ${\tia=\tiN/\N_0}$, versus the number of elements in each direction, $\ne$. 
	2D eigenproblems under maximum-continuity IGA and optimal rIGA discretizations with $\BS$ of 16 elements in each direction.}
	\label{fig.T1}
\end{figure}

\begin{figure}[!h]
	\centering
	\begin{subfigure}{0.49\textwidth}\centering
		\begin{overpic}[width=1.0\textwidth,trim={0cm 0cm 0cm 0cm} ,clip]{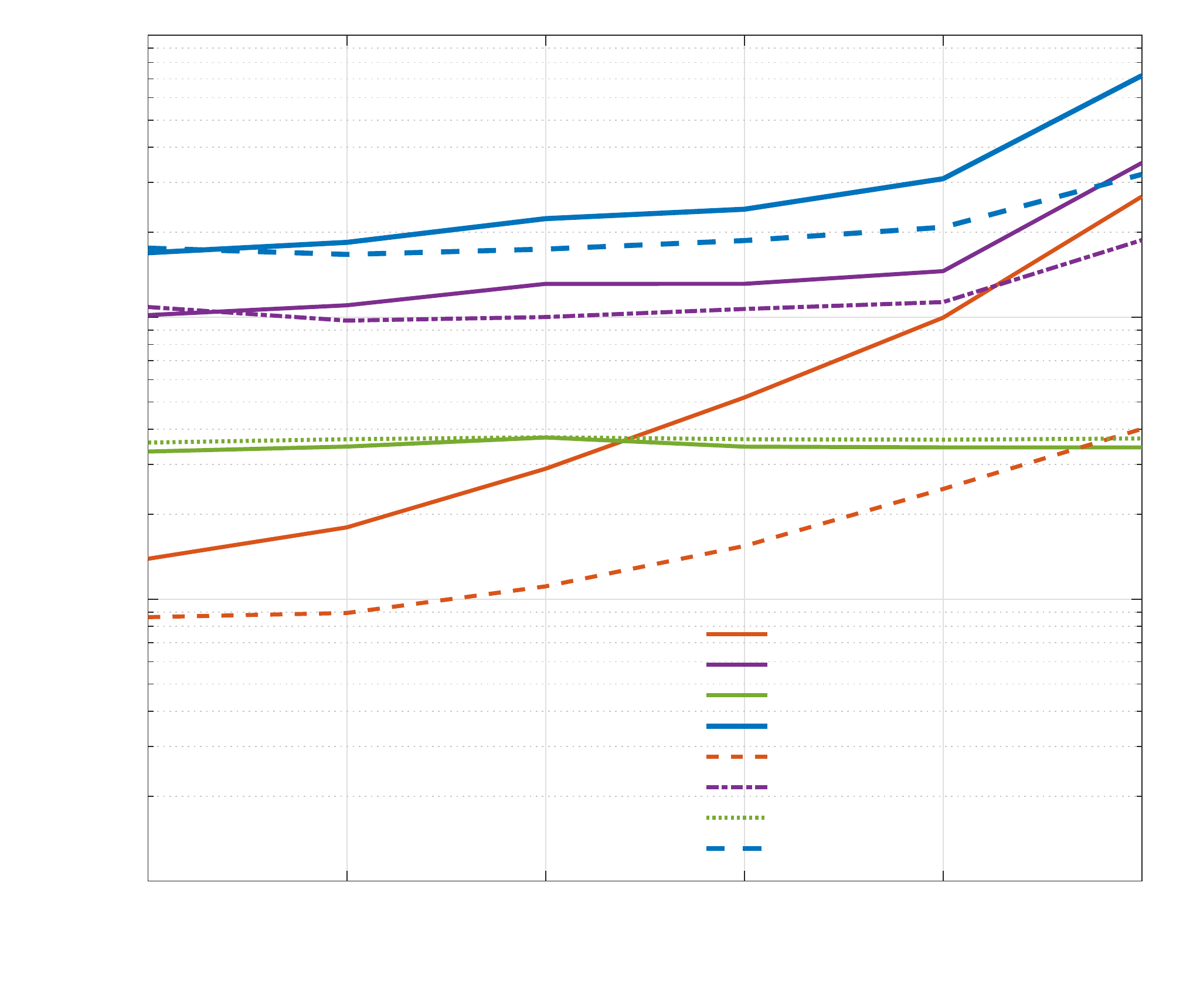}
			\put(30,2){\footnotesize No. of elements in each direction, $\ne$}
			\put(11.00,6.50){\footnotesize $2^{6}$}
			\put(27.60,6.50){\footnotesize $2^{7}$}
			\put(44.20,6.50){\footnotesize $2^{8}$}
			\put(60.80,6.50){\footnotesize $2^{9}$}
			\put(77.40,6.50){\footnotesize $2^{10}$}
			\put(94.00,6.50){\footnotesize $2^{11}$}
			\put(5.50,10.00){\footnotesize $10^{-8}$}
			\put(5.50,33.67){\footnotesize $10^{-7}$}
			\put(5.50,57.33){\footnotesize $10^{-6}$}
			\put(5.50,81.00){\footnotesize $10^{-5}$}
			\put(66.00,31.50){\scriptsize  IGA, factorization}
			\put(66.00,28.83){\scriptsize  IGA, \fb~elimination}
			\put(66.00,26.16){\scriptsize  IGA, mat--vec product}
			\put(66.00,23.49){\scriptsize  IGA, total time}
			\put(66.00,20.81){\scriptsize rIGA, factorization}
			\put(66.00,18.14){\scriptsize rIGA, \fb~elimination}
			\put(66.00,15.47){\scriptsize rIGA, mat--vec product}
			\put(66.00,12.80){\scriptsize rIGA, total time}
			\begin{turn}{90}
				\put(29,-3.5){\footnotesize Normalized time, $\tih$ (sec)}
			\end{turn}
		\end{overpic}
		\caption{$\p=2$}
	\end{subfigure}
	\begin{subfigure}{0.49\textwidth}\centering
		\begin{overpic}[width=1.0\textwidth,trim={0cm 0cm 0cm 0cm} ,clip]{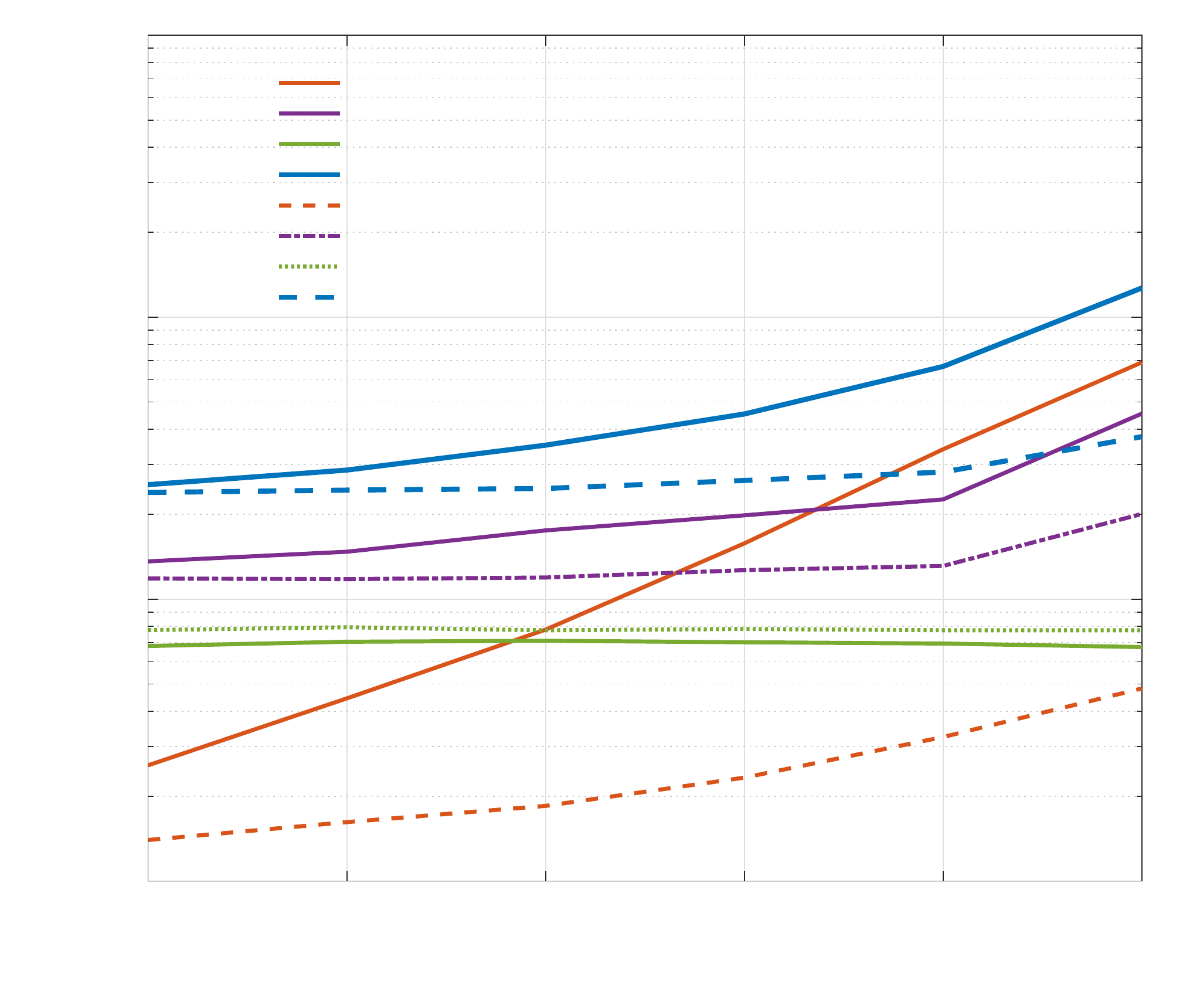}
			\put(30,2){\footnotesize No. of elements in each direction, $\ne$}
			\put(11.00,6.50){\footnotesize $2^{6}$}
			\put(27.60,6.50){\footnotesize $2^{7}$}
			\put(44.20,6.50){\footnotesize $2^{8}$}
			\put(60.80,6.50){\footnotesize $2^{9}$}
			\put(77.40,6.50){\footnotesize $2^{10}$}
			\put(94.00,6.50){\footnotesize $2^{11}$}
			\put(5.50,10.00){\footnotesize $10^{-7}$}
			\put(5.50,33.67){\footnotesize $10^{-6}$}
			\put(5.50,57.33){\footnotesize $10^{-5}$}
			\put(5.50,81.00){\footnotesize $10^{-4}$}
			\put(30.00,78.00){\scriptsize  IGA, factorization}
			\put(30.00,75.33){\scriptsize  IGA, \fb~elimination}
			\put(30.00,72.66){\scriptsize  IGA, mat--vec product}
			\put(30.00,69.99){\scriptsize  IGA, total time}
			\put(30.00,67.31){\scriptsize rIGA, factorization}
			\put(30.00,64.64){\scriptsize rIGA, \fb~elimination}
			\put(30.00,61.97){\scriptsize rIGA, mat--vec product}
			\put(30.00,59.30){\scriptsize rIGA, total time}
			\begin{turn}{90}
			    \put(29.00,-3.5){\footnotesize Normalized time, $\tih$ (sec)}
			\end{turn}
		\end{overpic}
		\caption{$\p=3$}
	\end{subfigure}
	\begin{subfigure}{0.49\textwidth}\centering
		\begin{overpic}[width=1.0\textwidth,trim={0cm 0cm 0cm 0cm} ,clip]{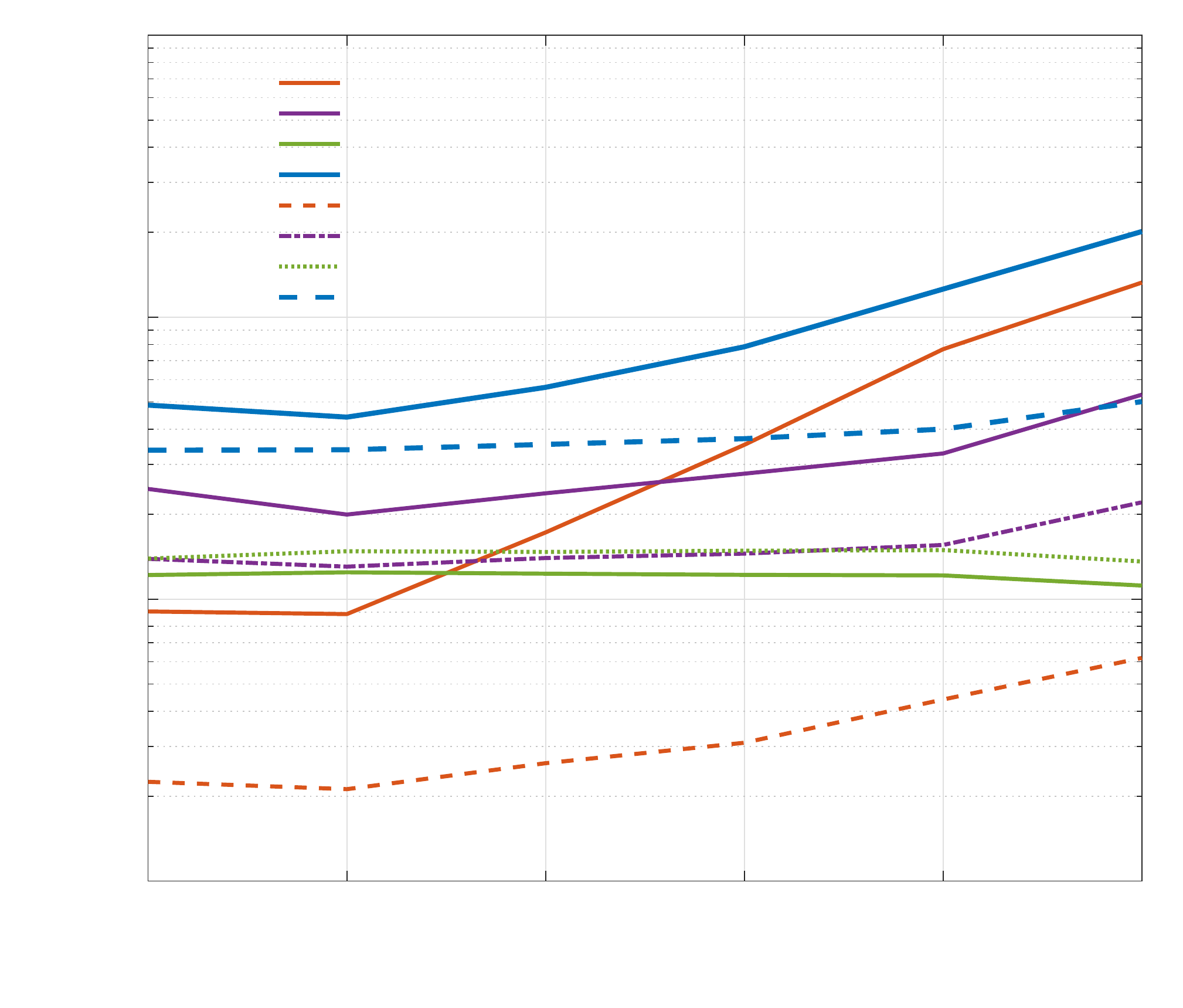}
			\put(30,2){\footnotesize No. of elements in each direction, $\ne$}
			\put(11.00,6.50){\footnotesize $2^{6}$}
			\put(27.60,6.50){\footnotesize $2^{7}$}
			\put(44.20,6.50){\footnotesize $2^{8}$}
			\put(60.80,6.50){\footnotesize $2^{9}$}
			\put(77.40,6.50){\footnotesize $2^{10}$}
			\put(94.00,6.50){\footnotesize $2^{11}$}
			\put(5.50,10.00){\footnotesize $10^{-7}$}
			\put(5.50,33.67){\footnotesize $10^{-6}$}
			\put(5.50,57.33){\footnotesize $10^{-5}$}
			\put(5.50,81.00){\footnotesize $10^{-4}$}
			\put(30.00,78.00){\scriptsize  IGA, factorization}
			\put(30.00,75.33){\scriptsize  IGA, \fb~elimination}
			\put(30.00,72.66){\scriptsize  IGA, mat--vec product}
			\put(30.00,69.99){\scriptsize  IGA, total time}
			\put(30.00,67.31){\scriptsize rIGA, factorization}
			\put(30.00,64.64){\scriptsize rIGA, \fb~elimination}
			\put(30.00,61.97){\scriptsize rIGA, mat--vec product}
			\put(30.00,59.30){\scriptsize rIGA, total time}
			\begin{turn}{90}
			    \put(29.00,-3.5){\footnotesize Normalized time, $\tih$ (sec)}
			\end{turn}
		\end{overpic}
		\caption{$\p=4$}
	\end{subfigure}
	\begin{subfigure}{0.49\textwidth}\centering
		\begin{overpic}[width=1.0\textwidth,trim={0cm 0cm 0cm 0cm} ,clip]{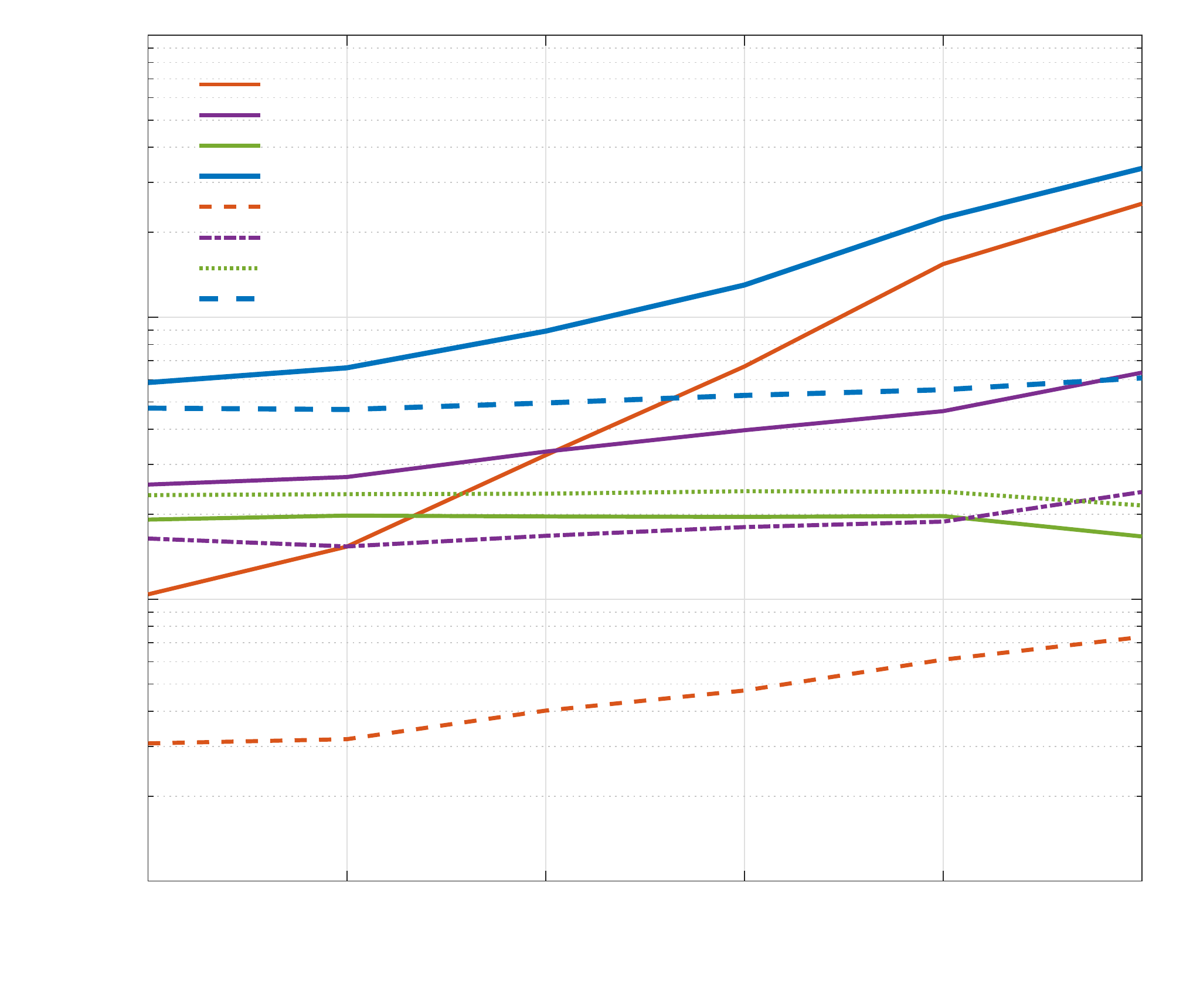}
			\put(30,2){\footnotesize No. of elements in each direction, $\ne$}
			\put(11.00,6.50){\footnotesize $2^{6}$}
			\put(27.60,6.50){\footnotesize $2^{7}$}
			\put(44.20,6.50){\footnotesize $2^{8}$}
			\put(60.80,6.50){\footnotesize $2^{9}$}
			\put(77.40,6.50){\footnotesize $2^{10}$}
			\put(94.00,6.50){\footnotesize $2^{11}$}
			\put(5.50,10.00){\footnotesize $10^{-7}$}
			\put(5.50,33.67){\footnotesize $10^{-6}$}
			\put(5.50,57.33){\footnotesize $10^{-5}$}
			\put(5.50,81.00){\footnotesize $10^{-4}$}
			\put(23.00,78.00){\scriptsize  IGA, factorization}
			\put(23.00,75.33){\scriptsize  IGA, \fb~elimination}
			\put(23.00,72.66){\scriptsize  IGA, mat--vec product}
			\put(23.00,69.99){\scriptsize  IGA, total time}
			\put(23.00,67.31){\scriptsize rIGA, factorization}
			\put(23.00,64.64){\scriptsize rIGA, \fb~elimination}
			\put(23.00,61.97){\scriptsize rIGA, mat--vec product}
			\put(23.00,59.30){\scriptsize rIGA, total time}
			\begin{turn}{90}
			    \put(29.00,-3.5){\footnotesize Normalized time, $\tih$ (sec)}
			\end{turn}
		\end{overpic}
		\caption{$\p=5$}
	\end{subfigure}
	\caption{Normalized average time per eigenvalue, ${\tih=\tiN/\N_0\.\N}$, versus the number of elements in each direction, $\ne$.
	2D eigenproblems under maximum-continuity IGA and optimal rIGA discretizations with $\BS$ of 16 elements in each direction.}
	\label{fig.That}
\end{figure}

\fig{\ref{fig.T1}} shows the average time required to compute an eigenvalue, ${\tia=\tiN/\N_0}$.
It also describes the contribution of each of the most expensive operations to the total time.
\fig{\ref{fig.That}} demonstrates the same results after normalizing the time per eigenvalue by the number of degrees of freedom as ${\tih=\tiN/\N_0\.\N}$.
Both figures confirm that matrix factorization is the most decisive numerical operation of the eigenanalysis for large problems.
rIGA with \BS~of 16 reduces the cost of the matrix factorization by a factor of up to ${\O(\p^{\.2}}$). 
It also
improves the cost of \fb~elimination, which is the governing cost in dealing with small problems. 
On the other hand,
the cost of mat--vec multiplication 
slightly increases when using rIGA,
which is not as decisive as the improvements of the other two operations.

\subsection{Time performance of 3D eigenproblems} 
\label{sub:CostEvaluation3D}

\begin{figure}[!h]
	\centering
	\begin{subfigure}{0.49\textwidth}\centering
		\begin{overpic}[width=1.0\textwidth,trim={0cm 0cm 0cm 0cm} ,clip]{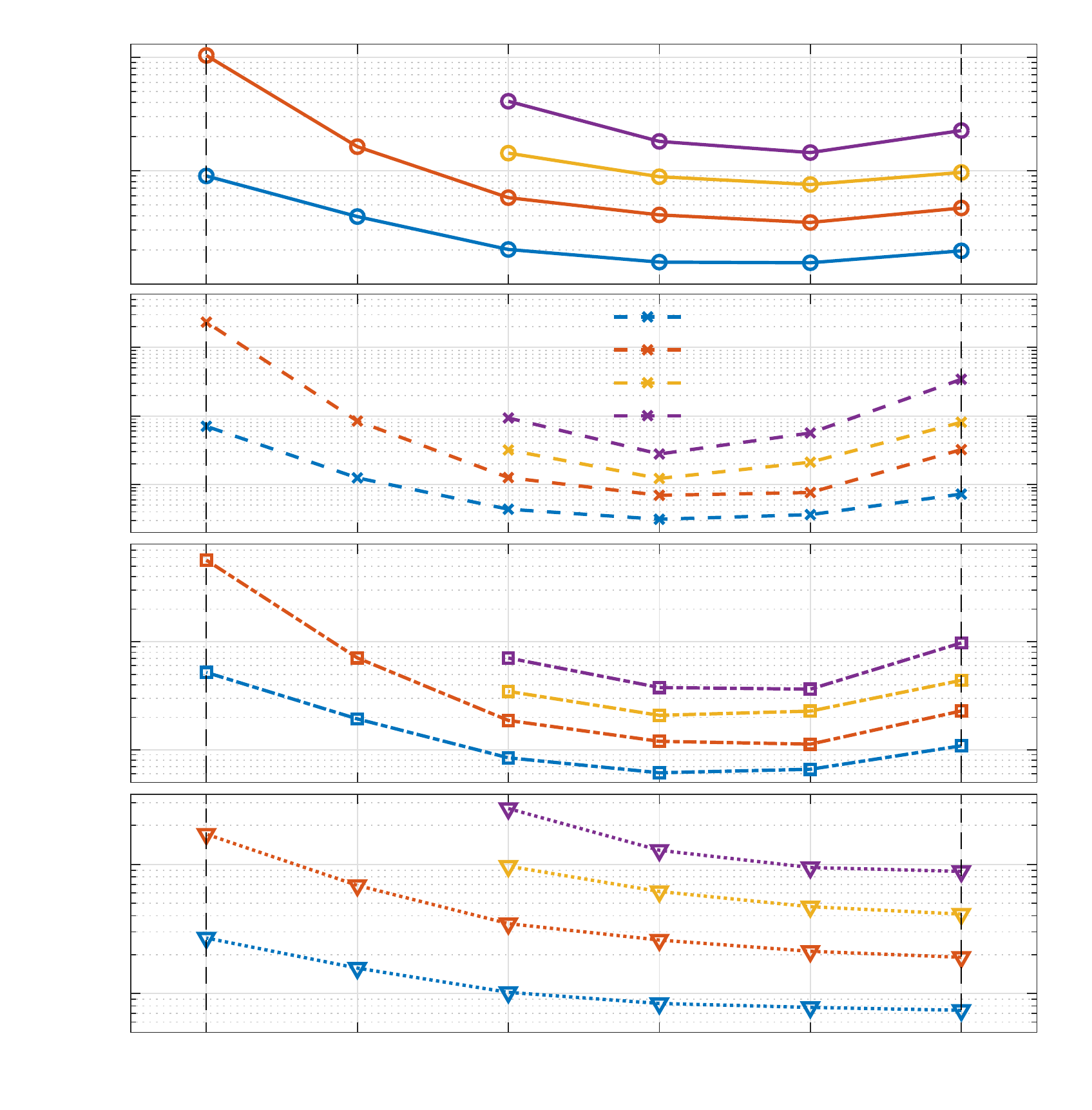}
			\put(46.5,0){\scriptsize  Blocksize}
			\put(62,62){\scriptsize  $\p=5$}
			\put(62,65){\scriptsize  $\p=4$}
			\put(62,68){\scriptsize  $\p=3$}
			\put(62,71){\scriptsize  $\p=2$}
			\put(20,93){\line(1,2){2.25}}
			\put(85.5,93){\line(-1,2){2.25}}
			\put(20,98){\scriptsize  $C^{\.0}\,(\rm{FEA})$}
			\put(78,98){\scriptsize  $C^{\p-1}\,(\rm{IGA})$}
			\put(16,3.50){\scriptsize $2^{\.0}$}
			\put(30,3.50){\scriptsize $2^{\.1}$}
			\put(43.90,3.50){\scriptsize $2^{\.2}$}
			\put(58.60,3.50){\scriptsize $2^{\.3}$}
			\put(72.30,3.50){\scriptsize $2^{\.4}$}
			\put(86.00,3.50){\scriptsize $2^{\.5}$}
			\put(6.5,94.00){\scriptsize $10^{\,4}$}
			\put(6.5,84.00){\scriptsize $10^{\,3}$}
			\put(6.5,74.00){\scriptsize $10^{\,2}$}
			\put(6.5,68.00){\scriptsize $10^{\,3}$}
			\put(6.5,61.00){\scriptsize $10^{\,2}$}
			\put(6.5,55.50){\scriptsize $10^{\,1}$}
			\put(6.5,42.00){\scriptsize $10^{\,3}$}
			\put(6.5,32.00){\scriptsize $10^{\,2}$}
			\put(6.5,21.50){\scriptsize $10^{\,3}$}
			\put(6.5,10.00){\scriptsize $10^{\,2}$}
			\begin{turn}{90}
				\put(9,-0.5){\scriptsize Mat--vec time}
			    \put(15,-3.5){\scriptsize (sec)}
			    \put(30.5,-0.5){\scriptsize \Fb~elimination}
			    \put(34,-3.5){\scriptsize time (sec)}
			    \put(55,-0.5){\scriptsize Factorization}
			    \put(57,-3.5){\scriptsize time (sec)}
			    \put(79,-0.5){\scriptsize Total time}
			    \put(82,-3.5){\scriptsize (sec)}
			\end{turn}
		\end{overpic}
		\caption{Mesh: $32^{\,3}$}
		\label{fig.timevsBS3D.32}
	\end{subfigure}
	\begin{subfigure}{0.49\textwidth}\centering
		\begin{overpic}[width=1.0\textwidth,trim={0cm 0cm 0cm 0cm} ,clip]{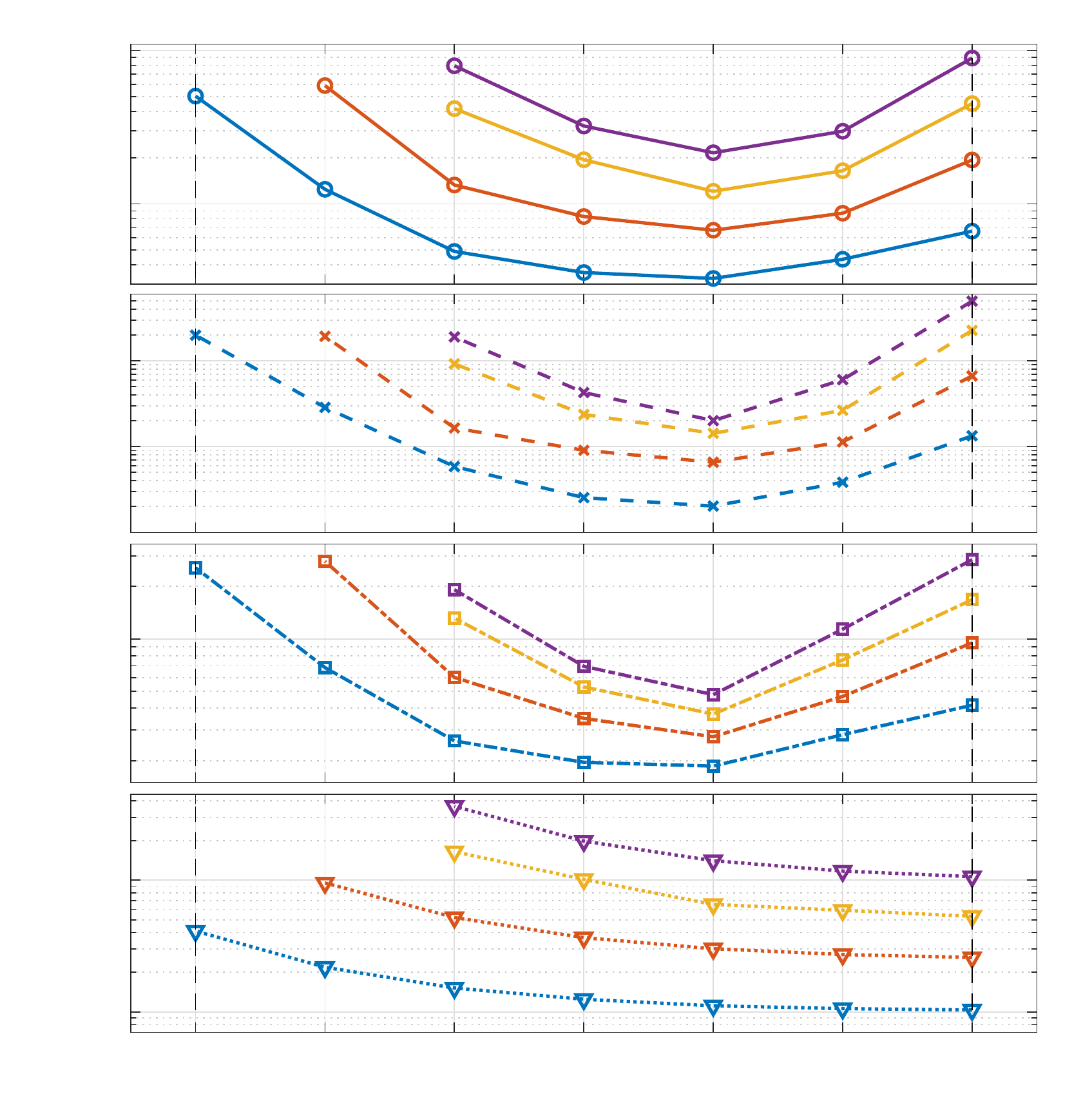}
			\put(46.5,0){\scriptsize  Blocksize}
			\put(18.5,93){\line(1,2){2.25}}
			\put(87,93){\line(-1,2){2.25}}
			\put(18.5,98){\scriptsize  $C^{\.0}\,(\rm{FEA})$}
			\put(80,98){\scriptsize  $C^{\p-1}\,(\rm{IGA})$}
			\put(16.00,3.50){\scriptsize $2^{\.0}$}
			\put(27.25,3.50){\scriptsize $2^{\.1}$}
			\put(39.00,3.50){\scriptsize $2^{\.2}$}
			\put(51.25,3.50){\scriptsize $2^{\.3}$}
			\put(63.25,3.50){\scriptsize $2^{\.4}$}
			\put(75,3.50){\scriptsize $2^{\.5}$}
			\put(86.00,3.50){\scriptsize $2^{\.6}$}
			\put(6.5,94.00){\scriptsize $10^{\,5}$}
			\put(6.5,81.00){\scriptsize $10^{\,4}$}
			\put(6.5,67.00){\scriptsize $10^{\,4}$}
			\put(6.5,59.00){\scriptsize $10^{\,3}$}
			\put(6.5,51.50){\scriptsize $10^{\,2}$}
			\put(6.5,42.00){\scriptsize $10^{\,4}$}
			\put(6.5,20.00){\scriptsize $10^{\,4}$}
			\put(6.5,8.00){\scriptsize $10^{\,3}$}
			\begin{turn}{90}
			    \put(9,-0.5){\scriptsize Mat--vec time}
			    \put(15,-3.5){\scriptsize (sec)}
			    \put(30.5,-0.5){\scriptsize \Fb~elimination}
			    \put(34,-3.5){\scriptsize time (sec)}
			    \put(55,-0.5){\scriptsize Factorization}
			    \put(57,-3.5){\scriptsize time (sec)}
			    \put(79,-0.5){\scriptsize Total time}
			    \put(82,-3.5){\scriptsize (sec)}
			\end{turn}
		\end{overpic}
		\caption{Mesh: $64^{\,3}$}
	\end{subfigure}
	\begin{subfigure}{0.49\textwidth}\centering
		\begin{overpic}[width=1.0\textwidth,trim={0cm 0cm 0cm 0cm} ,clip]{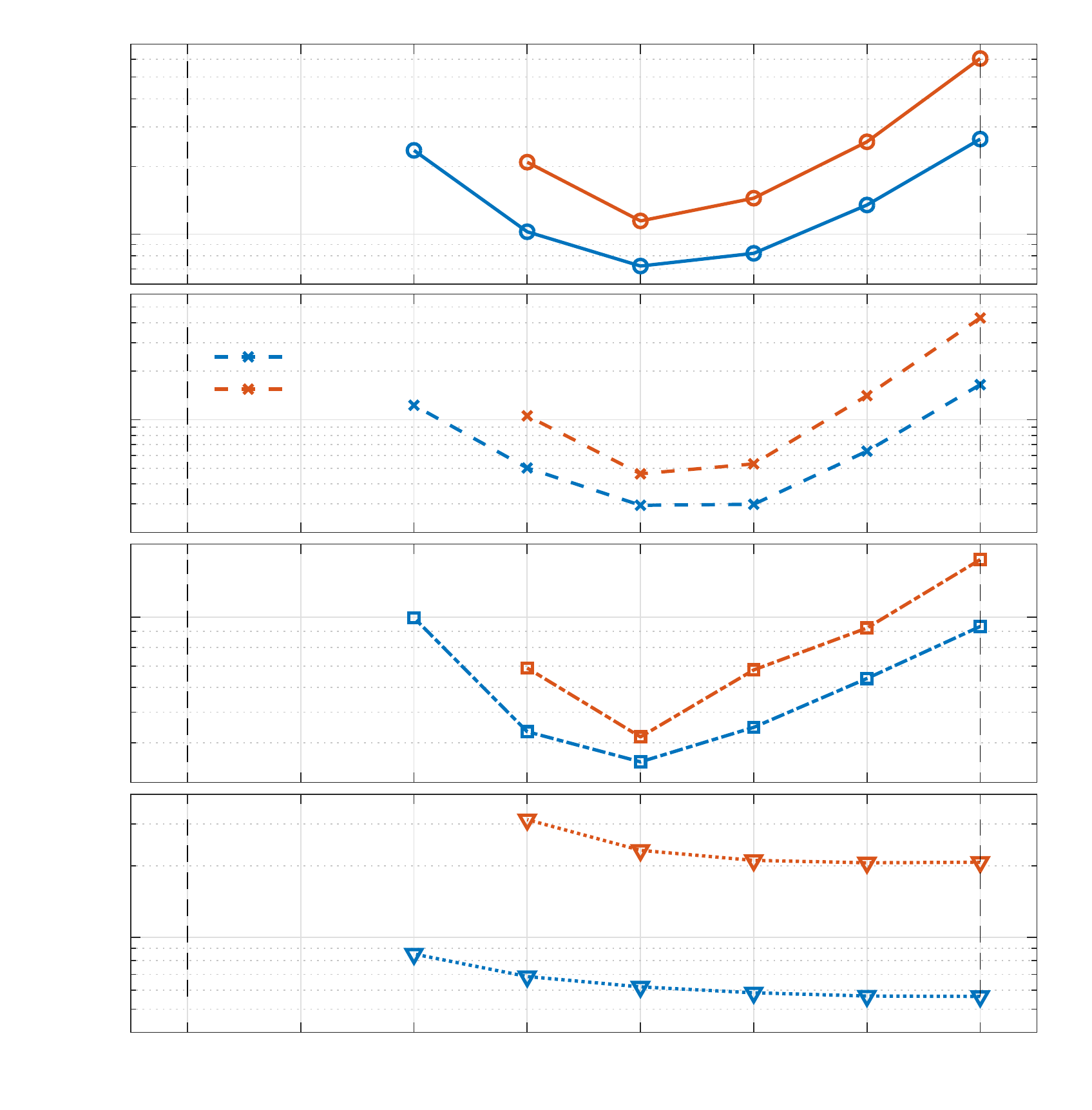}
			\put(46.5,0){\scriptsize  Blocksize}
			\put(26,64,5){\scriptsize  $\p=3$}
			\put(26,67.5){\scriptsize  $\p=2$}
			\put(18,90){\line(2,1){5}}
			\put(87.5,89){\line(-2,1){7}}
			\put(24,92){\scriptsize  $C^{\.0}\,(\rm{FEA})$}
			\put(65.5,92){\scriptsize  $C^{\p-1}\,(\rm{IGA})$}
			\put(15.50,3.50){\scriptsize $2^{\.0}$}
			\put(25.50,3.50){\scriptsize $2^{\.1}$}
			\put(35.50,3.50){\scriptsize $2^{\.2}$}
			\put(46.25,3.50){\scriptsize $2^{\.3}$}
			\put(56.50,3.50){\scriptsize $2^{\.4}$}
			\put(67.00,3.50){\scriptsize $2^{\.5}$}
			\put(77.20,3.50){\scriptsize $2^{\.6}$}
			\put(87.50,3.50){\scriptsize $2^{\.7}$}
			\put(6.5,78.00){\scriptsize $10^{\,5}$}
			\put(6.5,61.00){\scriptsize $10^{\,5}$}
			\put(6.5,43.50){\scriptsize $10^{\,5}$}
			\put(6.5,15.00){\scriptsize $10^{\,4}$}
			\begin{turn}{90}
			    \put(9,-0.5){\scriptsize Mat--vec time}
			    \put(15,-3.5){\scriptsize (sec)}
			    \put(30.5,-0.5){\scriptsize \Fb~elimination}
			    \put(34,-3.5){\scriptsize time (sec)}
			    \put(55,-0.5){\scriptsize Factorization}
			    \put(57,-3.5){\scriptsize time (sec)}
			    \put(79,-0.5){\scriptsize Total time}
			    \put(82,-3.5){\scriptsize (sec)}
			\end{turn}
		\end{overpic}
		\caption{Mesh: $128^{\,3}$}
		\label{fig.timevsBS3D.128}
	\end{subfigure}
	\caption{The total computational times and those of the most expensive numerical operations for finding the first ${\N_0=1024}$ eigenpairs of the 3D Laplace operator with different mesh sizes of ${\ne=2^s}$ ${(s=5,6,7)}$ and degrees ${\p=2,3,4,5}$.
	We test rIGA discretizations with different $\BS$s $(2^{s-\.\l})$ obtained using different levels of partitioning ${\l=0,1,...,s}$. Since the eigenanalysis in 3D case is highly demanding, for ${\ne=128}$, we only test up to ${\p=3}$.}
	\label{fig.timevsBS3D}
\end{figure}

\fig{\ref{fig.timevsBS3D}} depicts the times required to compute ${\N_0=2^{10}}$ eigenpairs of the Laplace operator in 3D versus different $\BS$s of the rIGA discretizations.
Similar to the 2D test cases, the factorization time of rIGA reaches a minimum at ${\BS}$ of 16 elements expect for ${\ne=32}$, where referring to \fig{\ref{fig.timevsBS3D.32}}, the optimal $\BS$ for matrix factorization is 8 elements
(see the same inference in~\cite{ Garcia2017} for optimal $\BS$ of small and large domains). 
However, the maximum computational saving for the total elapsed time, considering all numerical operations of the eigenanalysis, is achieved by employing an rIGA discretization with $\BS$ of 16 elements in each direction.
\tab{\ref{tab:T1_3D}} reports the number of times we perform each operation when finding ${\N_0}$ eigenpairs
as well as 
the average computational times per eigenvalue, ${\tia=\tiN/\N_0}$. 
In 3D systems, 
we observe improvements of $\O(\p^2)$ in matrix factorization and $\O(\p)$ in the total eigenanalysis
when using rIGA with ${\ne\geq64}$. 
We find ${\N_0=2^{\.10}}$ eigenpairs with a ${128^3}$ mesh and cubic bases in 168 hours using IGA and 32 hours using rIGA. (see \fig{\ref{fig.timevsBS3D.128}}).
However, we expect an improvement of ${\O(\p^{\.2})}$ 
for sufficiently large problems.

\begin{table}[!h]
\centering
\caption{The average computational times per eigenvalue, ${\tia=\tiN/\N_0}$, and the number of executions of the most expensive operations for finding ${\N_0=1024}$ eigenpairs of the 3D test cases illustrated in \fig{\ref{fig.timevsBS3D}}.
The optimal $\BS$ of 16 elements is considered for comparisons with the rIGA discretizations. 
Note that for ${\ne=32}$, a higher improvement in matrix factorization is achievable with $\BS$ of 8 elements. However, the overall eigenanalysis is better improved with $\BS$ of 16 elements.}
\label{tab:T1_3D}
\footnotesize
\begin{tabular}{@{}llllllllllllll@{}}
\toprule
\multicolumn{3}{l}{Discretization}  & \multicolumn{3}{l}{Factorization} & \multicolumn{3}{l}{\Fb~elimination} & \multicolumn{3}{l}{Mat--vec product} & \multicolumn{2}{l}{Total} \\ 
\cmidrule(r){1-3} \cmidrule(lr){4-6} \cmidrule(lr){7-9} \cmidrule(lr){10-12} \cmidrule(l){13-14} 
Mesh & $\p$ & Method & $\Nfa$ &
\begin{tabular}[c]{@{}l@{}}$\tiax{\text{fact}}$\\(sec)\end{tabular} & 
\begin{tabular}[c]{@{}l@{}}Improved\\by\end{tabular} & $\Nfb$ & \begin{tabular}[c]{@{}l@{}}$\tiax{\fb}$\\(sec)\end{tabular} & 
\begin{tabular}[c]{@{}l@{}}Improved\\by\end{tabular} & $\Nmv$ & \begin{tabular}[c]{@{}l@{}}$\tiax{\text{m--v}}$\\(sec)\end{tabular} &
\begin{tabular}[c]{@{}l@{}}Degraded\\by\end{tabular} & \begin{tabular}[c]{@{}l@{}}$\tia$\\(sec)\end{tabular} & 
\begin{tabular}[c]{@{}l@{}}Improved\\by\end{tabular} \\ \midrule
\multirow{8}{*}{$32^{\,3}$} 
& \multirow{2}{*}{2} & IGA & 23 & 0.007 & \multirow{2}{*}{1.994} & 4425 & 0.106 & \multirow{2}{*}{1.648} & 13213 & 0.072 & \multirow{2}{*}{0.949} & 0.191 & \multirow{2}{*}{1.273} \\ 
				 &  & rIGA & 22 & 0.003 &						 & 4376 & 0.064 &					     & 13113 & 0.076 &  & 0.150 &  \\ [2pt] 
& \multirow{2}{*}{3} & IGA & 23 & 0.031 & \multirow{2}{*}{4.234} & 4584 & 0.223 & \multirow{2}{*}{2.030} & 13673 & 0.186 & \multirow{2}{*}{0.898} & 0.457 & \multirow{2}{*}{1.340} \\ 
				 &  & rIGA & 24 & 0.007 &						 & 4556 & 0.110 &					     & 13603 & 0.207 &  & 0.340 &  \\ [2pt] 
& \multirow{2}{*}{4} & IGA & 25 & 0.079 & \multirow{2}{*}{3.831} & 4725 & 0.429 & \multirow{2}{*}{1.927} & 14174 & 0.402 & \multirow{2}{*}{0.876} & 0.940 & \multirow{2}{*}{1.277} \\ 
				 &  & rIGA & 23 & 0.020 &						 & 4633 & 0.222 &					     & 13898 & 0.458 &  & 0.736 &  \\ [2pt] 
& \multirow{2}{*}{5} & IGA & 24 & 0.335 & \multirow{2}{*}{6.077} & 4735 & 0.952 & \multirow{2}{*}{2.680} & 14203 & 0.861 & \multirow{2}{*}{0.933} & 2.200 & \multirow{2}{*}{1.563} \\ 
				 &  & rIGA & 24 & 0.055 &						 & 4607 & 0.355 &					     & 13819 & 0.923 &  & 1.407 &  \\ \cmidrule(){1-14}
\multirow{8}{*}{$64^{\,3}$} 
& \multirow{2}{*}{2} & IGA & 25 & 1.305 & \multirow{2}{*}{6.599} & 4595 & 4.061 & \multirow{2}{*}{2.237} & 13707 & 1.010 & \multirow{2}{*}{0.929} & 6.486 & \multirow{2}{*}{2.034} \\ 
				 &  & rIGA & 23 & 0.197 &						 & 4563 & 1.815 &					     & 13567 & 1.087 &  & 3.187 &  \\ [2pt] 
& \multirow{2}{*}{3} & IGA & 22 & 6.517 & \multirow{2}{*}{10.14} & 4653 & 9.272 & \multirow{2}{*}{3.466} & 13958 & 2.526 & \multirow{2}{*}{0.857} & 18.87 & \multirow{2}{*}{2.871} \\ 
				 &  & rIGA & 22 & 0.642 &						 & 4609 & 2.675 &					     & 13827 & 2.946 &  & 6.573 &  \\ [2pt] 
& \multirow{2}{*}{4} & IGA & 24 & 22.14 & \multirow{2}{*}{15.86} & 4619 & 16.38 & \multirow{2}{*}{4.536} & 13856 & 5.192 & \multirow{2}{*}{0.812} & 44.00 & \multirow{2}{*}{3.727} \\ 
				 &  & rIGA & 23 & 1.396 &						 & 4757 & 3.612 &					     & 14270 & 6.391 &  & 11.80 &  \\ [2pt] 
& \multirow{2}{*}{5} & IGA & 24 & 48.57 & \multirow{2}{*}{24.74} & 4733 & 27.83 & \multirow{2}{*}{5.961} & 14197 & 10.38 & \multirow{2}{*}{0.756} & 87.24 & \multirow{2}{*}{4.144} \\ 
				 &  & rIGA & 22 & 1.963 &						 & 4661 & 4.669 &					     & 13981 & 13.72 &  & 21.04 &  \\ \cmidrule(){1-14}
\multirow{4}{*}{$128^{\,3}$}
& \multirow{2}{*}{2} & IGA & 24 & 160.9 & \multirow{2}{*}{5.593} & 4781 & 91.16 & \multirow{2}{*}{2.680} & 14342 & 5.525 & \multirow{2}{*}{0.911} & 258.6 & \multirow{2}{*}{3.669} \\ 
 			     &  & rIGA & 24 & 28.76 &						 & 4607 & 34.00 &					     & 13820 & 6.061 &  & 70.48 &  \\ [2pt] 
& \multirow{2}{*}{3} & IGA & 25 & 417.1 & \multirow{2}{*}{9.254} & 4684 & 148.3 & \multirow{2}{*}{3.636} & 14051 & 20.21 & \multirow{2}{*}{0.891} & 590.1 & \multirow{2}{*}{5.275} \\ 
    		 	 &  & rIGA & 23 & 45.07 &						 & 4673 & 40.79 &					     & 14018 & 22.68 &  & 111.8 &  \\ \bottomrule \\[10pt]
\end{tabular}
\end{table}

\begin{figure}[!h]
	\begin{overpic}[height=0.31\textheight,trim={0cm 0cm 0cm 0cm} ,clip]{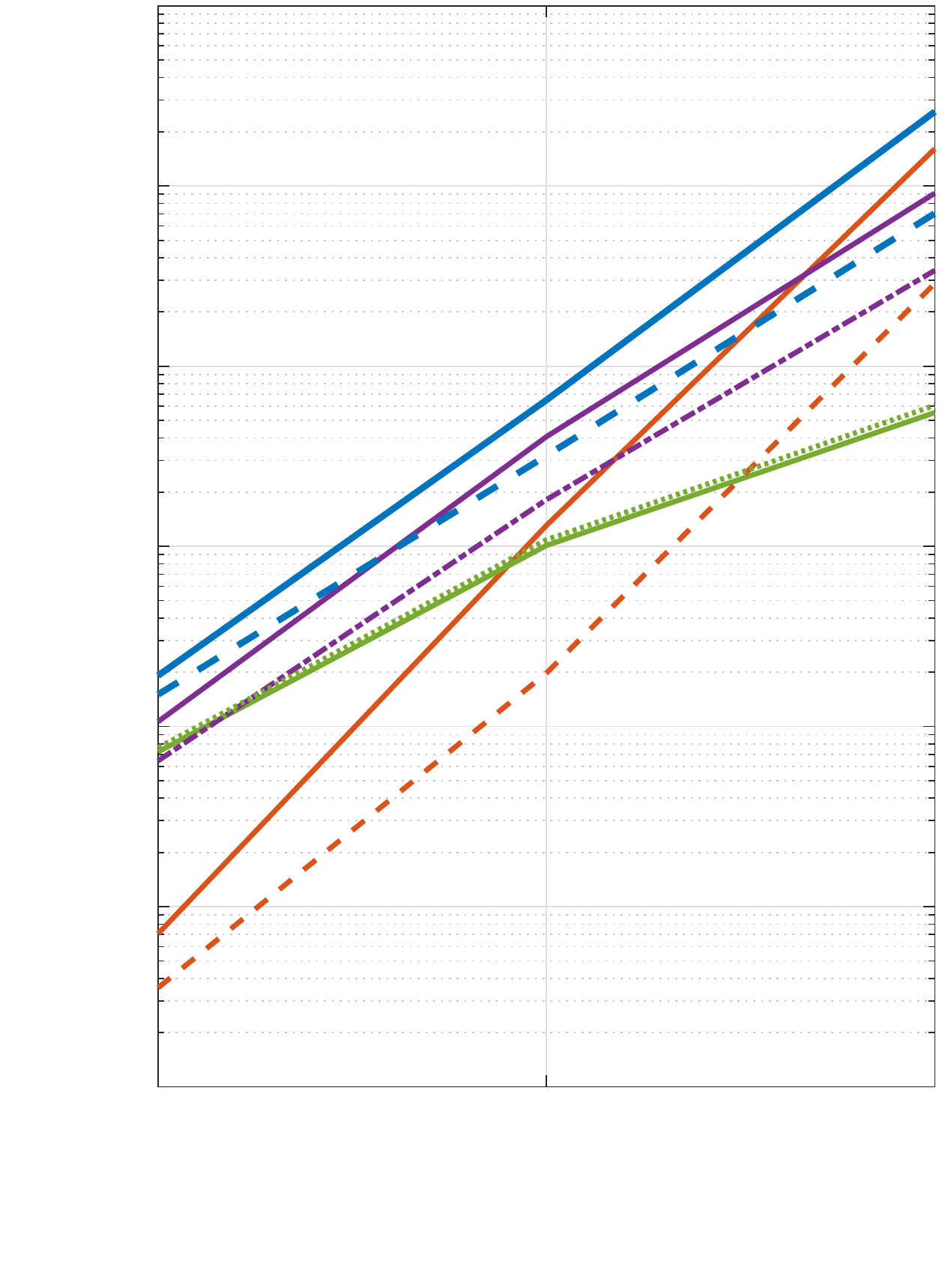}
		\put(27,6){\footnotesize No. of elements in}
		\put(28,2){\footnotesize each direction, $\ne$}
		\put(34,-6){\small (a) $\p=2$}
		\put(11.00,11.00){\footnotesize $2^{5}$}
		\put(40.50,11.00){\footnotesize $2^{6}$}
		\put(70.00,11.00){\footnotesize $2^{7}$}
		\put(4.50,15.00){\footnotesize $10^{-3}$}
		\put(4.50,27.83){\footnotesize $10^{-2}$}
		\put(4.50,41.67){\footnotesize $10^{-1}$}
		\put(5.00,55.50){\footnotesize $10^{\,0}$}
		\put(5.00,69.33){\footnotesize $10^{\,1}$}
		\put(5.00,83.17){\footnotesize $10^{\,2}$}
		\put(5.00,97.00){\footnotesize $10^{\,3}$}
		\begin{turn}{90}
		    \put(34.00,0){\footnotesize Time per eigenvalue, $\tia$ (sec)}
		\end{turn}
	\end{overpic}\,\,\,
	\begin{overpic}[height=0.31\textheight,trim={0cm 0cm 0cm 0cm} ,clip]{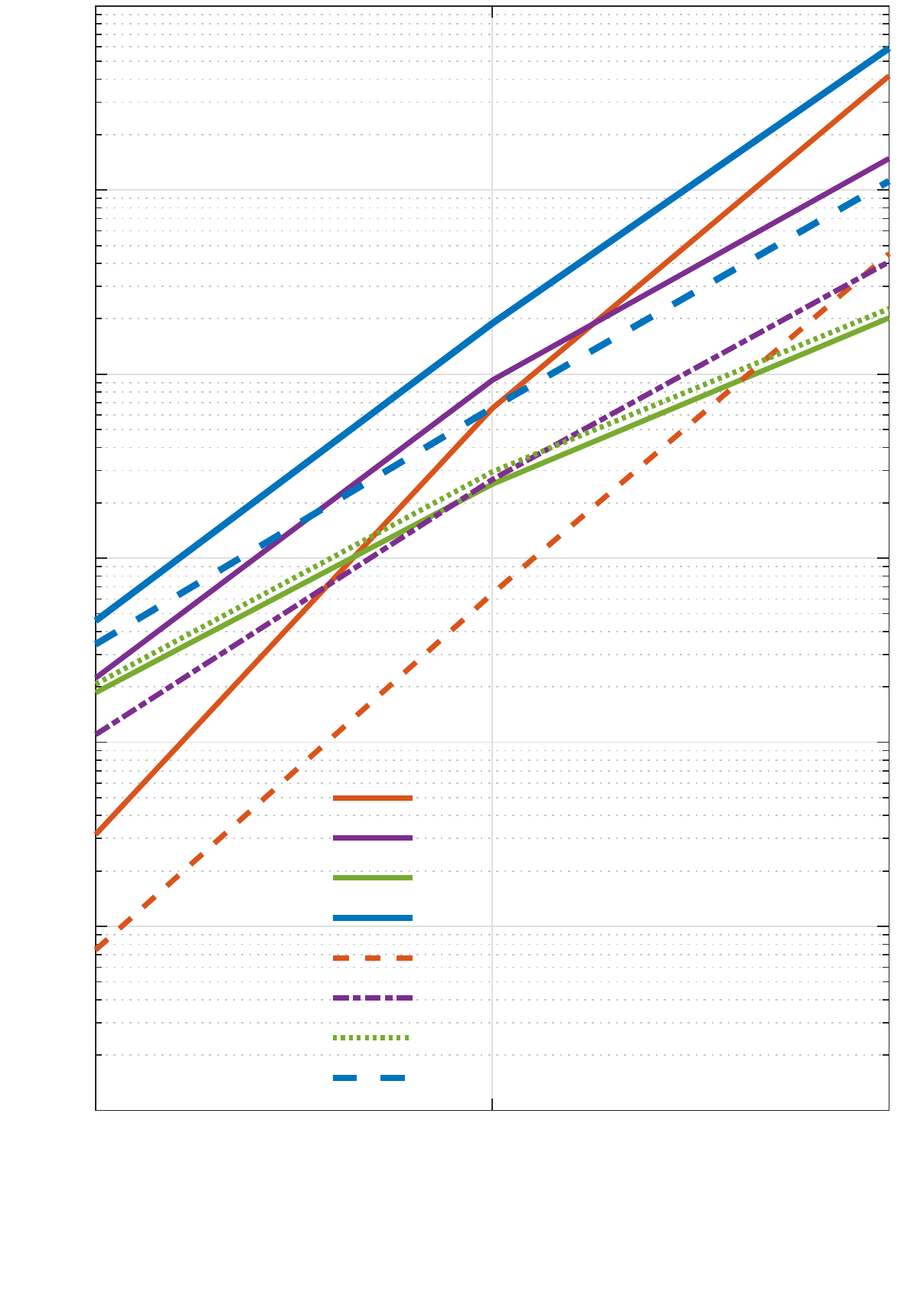}
		\put(22,6){\footnotesize No. of elements in}
		\put(23,2){\footnotesize each direction, $\ne$}
		\put(29,-6){\small (b) $\p=3$}
		\put(6.00,11.00){\footnotesize $2^{5}$}
		\put(35.50,11.00){\footnotesize $2^{6}$}
		\put(65.00,11.00){\footnotesize $2^{7}$}
		\put(-0.50,15.00){\footnotesize $10^{-3}$}
		\put(-0.50,27.83){\footnotesize $10^{-2}$}
		\put(-0.50,41.67){\footnotesize $10^{-1}$}
		\put(0.75,55.50){\footnotesize $10^{\,0}$}
		\put(0.75,69.33){\footnotesize $10^{\,1}$}
		\put(0.75,83.17){\footnotesize $10^{\,2}$}
		\put(0.75,97.00){\footnotesize $10^{\,3}$}
		\put(32.00,38.50){\scriptsize  IGA, factorization}
		\put(32.00,35.43){\scriptsize  IGA, \fb~elimination}
		\put(32.00,32.36){\scriptsize  IGA, mat--vec product}
		\put(32.00,29.29){\scriptsize  IGA, total time}
		\put(32.00,26.21){\scriptsize rIGA, factorization}
		\put(32.00,23.14){\scriptsize rIGA, \fb~elimination}
		\put(32.00,20.07){\scriptsize rIGA, mat--vec product}
		\put(32.00,17.00){\scriptsize rIGA, total time}
	\end{overpic}\,\,\,
	\begin{overpic}[height=0.31\textheight,trim={0cm 0cm 0cm 0cm} ,clip]{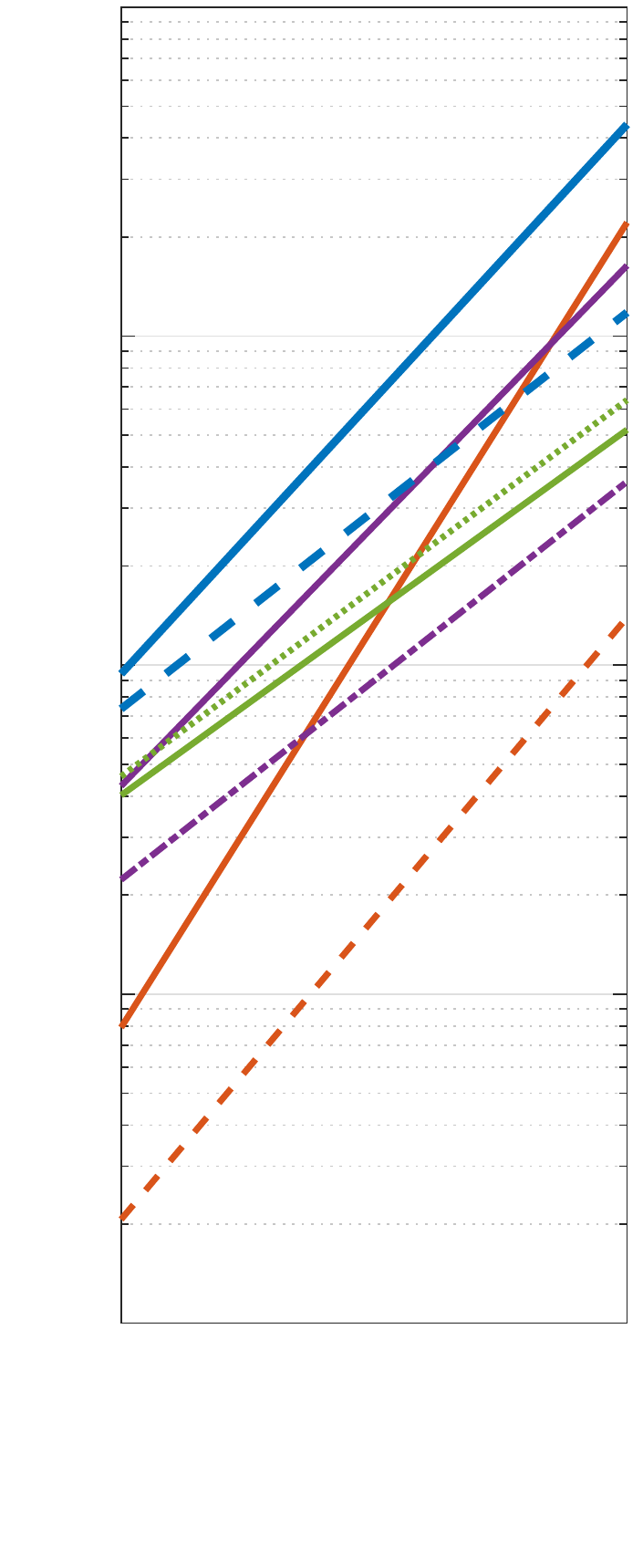}
		\put(9,6){\footnotesize No. of elements in}
		\put(10,2){\footnotesize each direction, $\ne$}
		\put(16,-6){\small (c) $\p=4$}
		\put(6.00,11.00){\footnotesize $2^{5}$}
		\put(37.00,11.00){\footnotesize $2^{6}$}
		\put(-0.50,15.00){\footnotesize $10^{-2}$}
		\put(-0.50,34.75){\footnotesize $10^{-1}$}
		\put(0.50,55.50){\footnotesize $10^{\,0}$}
		\put(0.50,76.25){\footnotesize $10^{\,1}$}
		\put(0.50,97.00){\footnotesize $10^{\,2}$}
	\end{overpic}\,\,\,
	\begin{overpic}[height=0.31\textheight,trim={0cm 0cm 0cm 0cm} ,clip]{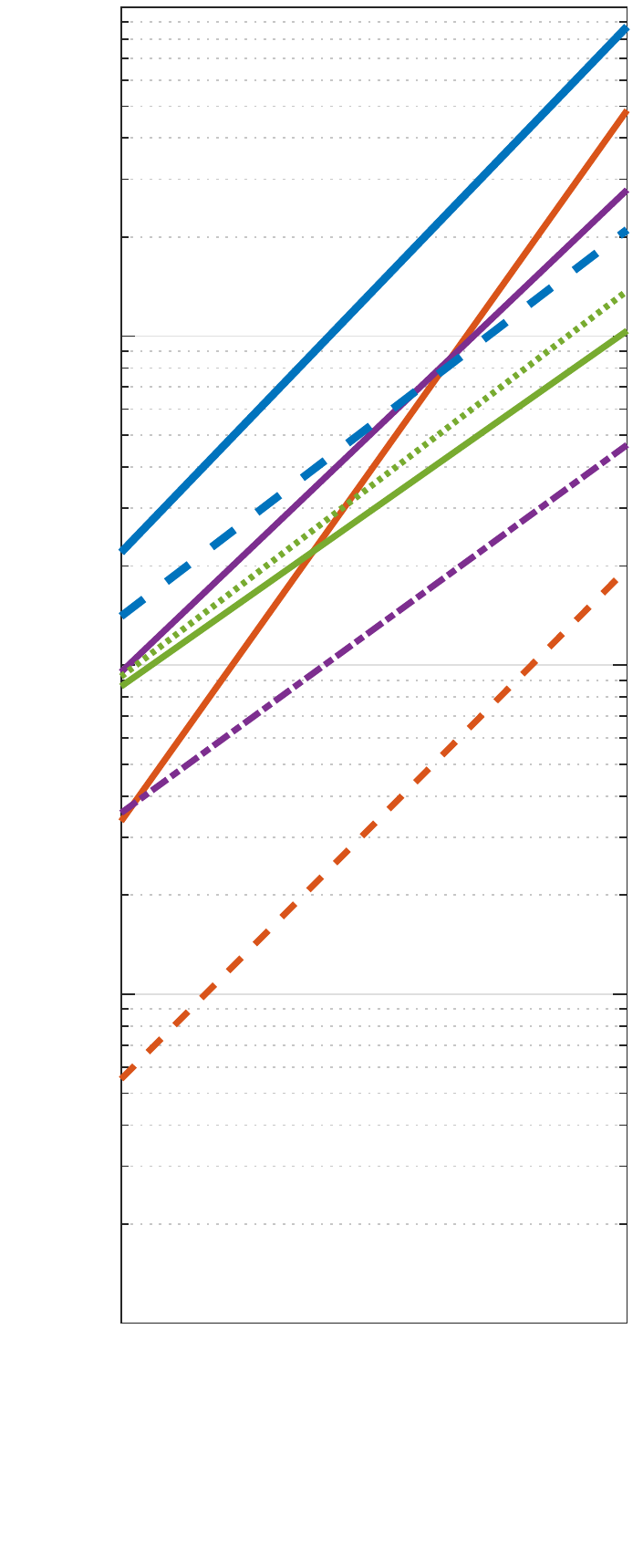}
		\put(9,6){\footnotesize No. of elements in}
		\put(10,2){\footnotesize each direction, $\ne$}
		\put(16,-6){\small (d) $\p=5$}
		\put(6.00,11.00){\footnotesize $2^{5}$}
		\put(37.00,11.00){\footnotesize $2^{6}$}
		\put(-0.50,15.00){\footnotesize $10^{-2}$}
		\put(-0.50,34.75){\footnotesize $10^{-1}$}
		\put(0.50,55.50){\footnotesize $10^{\,0}$}
		\put(0.50,76.25){\footnotesize $10^{\,1}$}
		\put(0.50,97.00){\footnotesize $10^{\,2}$}
	\end{overpic} \\
	\caption{Average computational time of each numerical operation per eigenvalue, ${\tia=\tiN/\N_0}$, versus the number of elements in each direction, $\ne$. 
	3D eigenproblems under maximum-continuity IGA and optimal rIGA discretizations with $\BS$ of 16 elements in each direction
	(for degrees ${\p=4~\rm{and}~5}$, we only test up to ${\ne=64}$).}
	\label{fig.T13D}
\end{figure}

\begin{figure}[!h]
	\begin{overpic}[height=0.31\textheight,trim={0cm 0cm 0cm 0cm} ,clip]{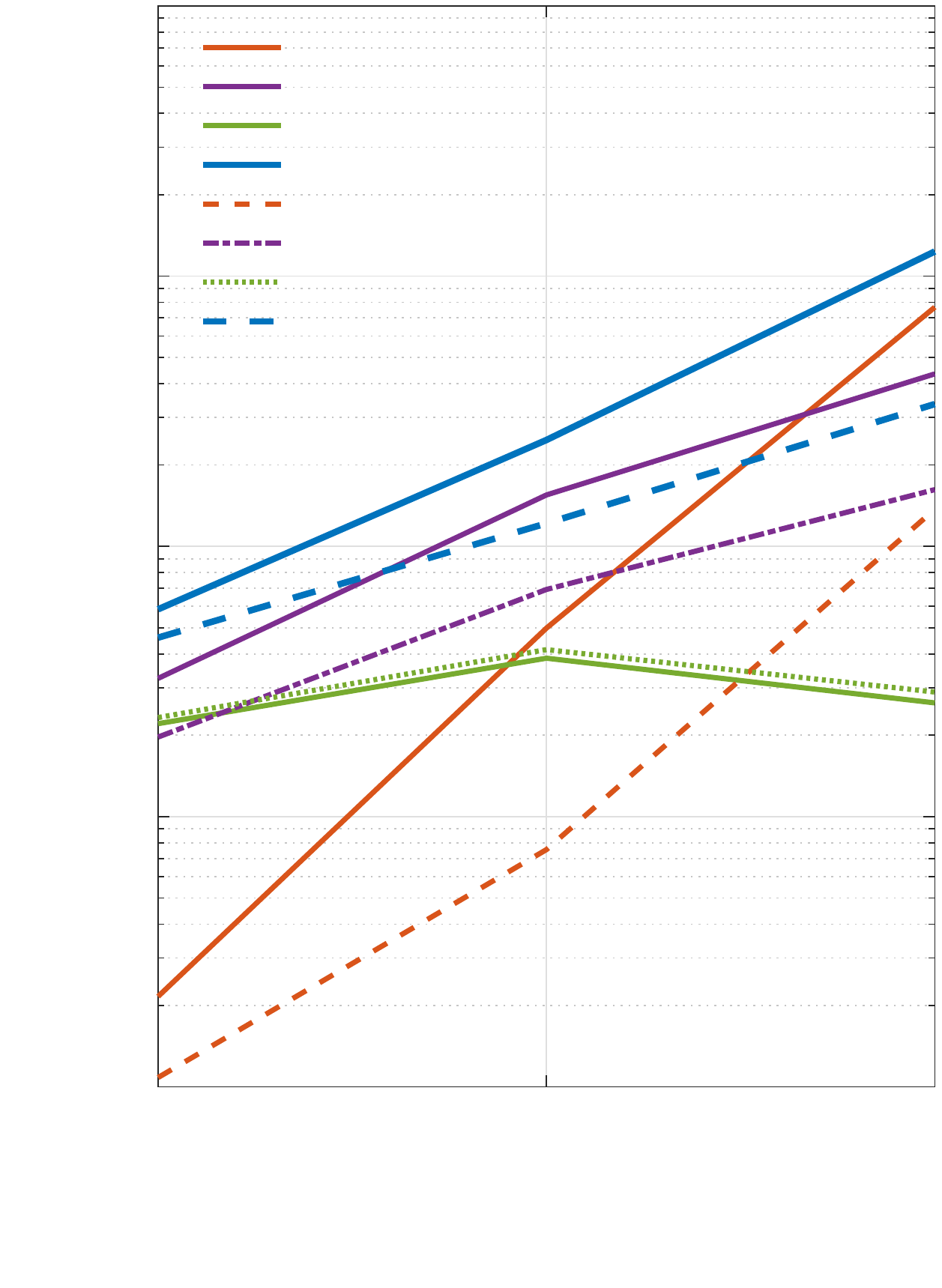}
		\put(27,6){\footnotesize No. of elements in}
		\put(28,2){\footnotesize each direction, $\ne$}
		\put(34,-6){\small (a) $\p=2$}
		\put(11.00,11.00){\footnotesize $2^{5}$}
		\put(40.50,11.00){\footnotesize $2^{6}$}
		\put(70.00,11.00){\footnotesize $2^{7}$}
		\put(4.5,15.00){\footnotesize $10^{-7}$}
		\put(4.5,34.75){\footnotesize $10^{-6}$}
		\put(4.5,55.50){\footnotesize $10^{-5}$}
		\put(4.5,76.25){\footnotesize $10^{-4}$}
		\put(4.5,97.00){\footnotesize $10^{-3}$}
		\put(23.00,96.00){\scriptsize  IGA, factorization}
		\put(23.00,92.93){\scriptsize  IGA, \fb~elimination}
		\put(23.00,89.86){\scriptsize  IGA, mat--vec product}
		\put(23.00,86.79){\scriptsize  IGA, total time}
		\put(23.00,83.71){\scriptsize rIGA, factorization}
		\put(23.00,80.64){\scriptsize rIGA, \fb~elimination}
		\put(23.00,77.57){\scriptsize rIGA, mat--vec product}
		\put(23.00,74.50){\scriptsize rIGA, total time}
		\begin{turn}{90}
		    \put(37.00,0){\footnotesize Normalized time, $\tih$ (sec)}
		\end{turn} 
	\end{overpic}\,\,\,
	\begin{overpic}[height=0.31\textheight,trim={0cm 0cm 0cm 0cm} ,clip]{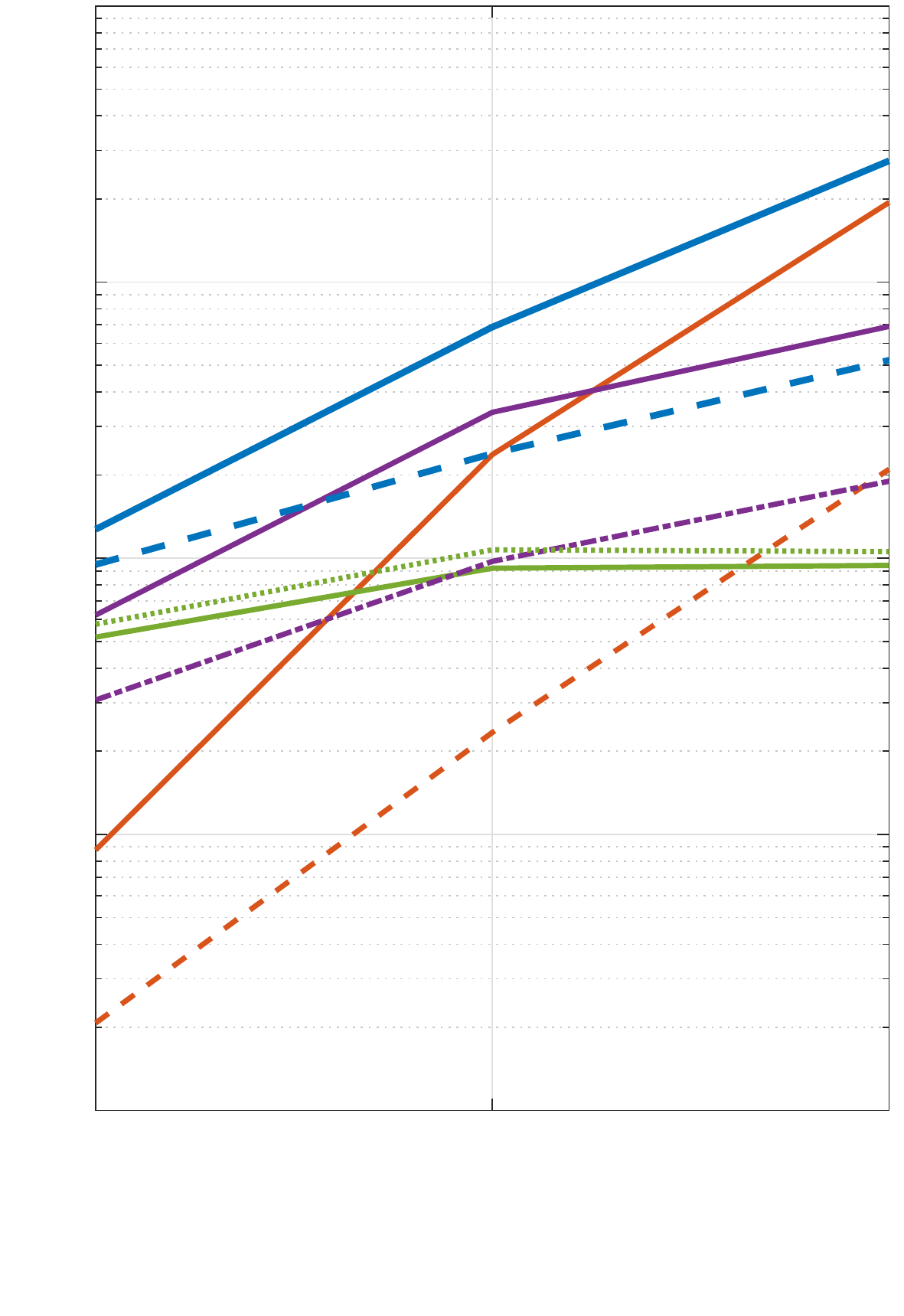}
		\put(22,6){\footnotesize No. of elements in}
		\put(23,2){\footnotesize each direction, $\ne$}
		\put(29,-6){\small (b) $\p=3$}
		\put(6.00,11.00){\footnotesize $2^{5}$}
		\put(35.50,11.00){\footnotesize $2^{6}$}
		\put(65.00,11.00){\footnotesize $2^{7}$}
		\put(-0.50,15.00){\footnotesize $10^{-7}$}
		\put(-0.50,34.75){\footnotesize $10^{-6}$}
		\put(-0.50,55.50){\footnotesize $10^{-5}$}
		\put(-0.50,76.25){\footnotesize $10^{-4}$}
		\put(-0.50,97.00){\footnotesize $10^{-3}$}
	\end{overpic}\,\,\,
	\begin{overpic}[height=0.31\textheight,trim={0cm 0cm 0cm 0cm} ,clip]{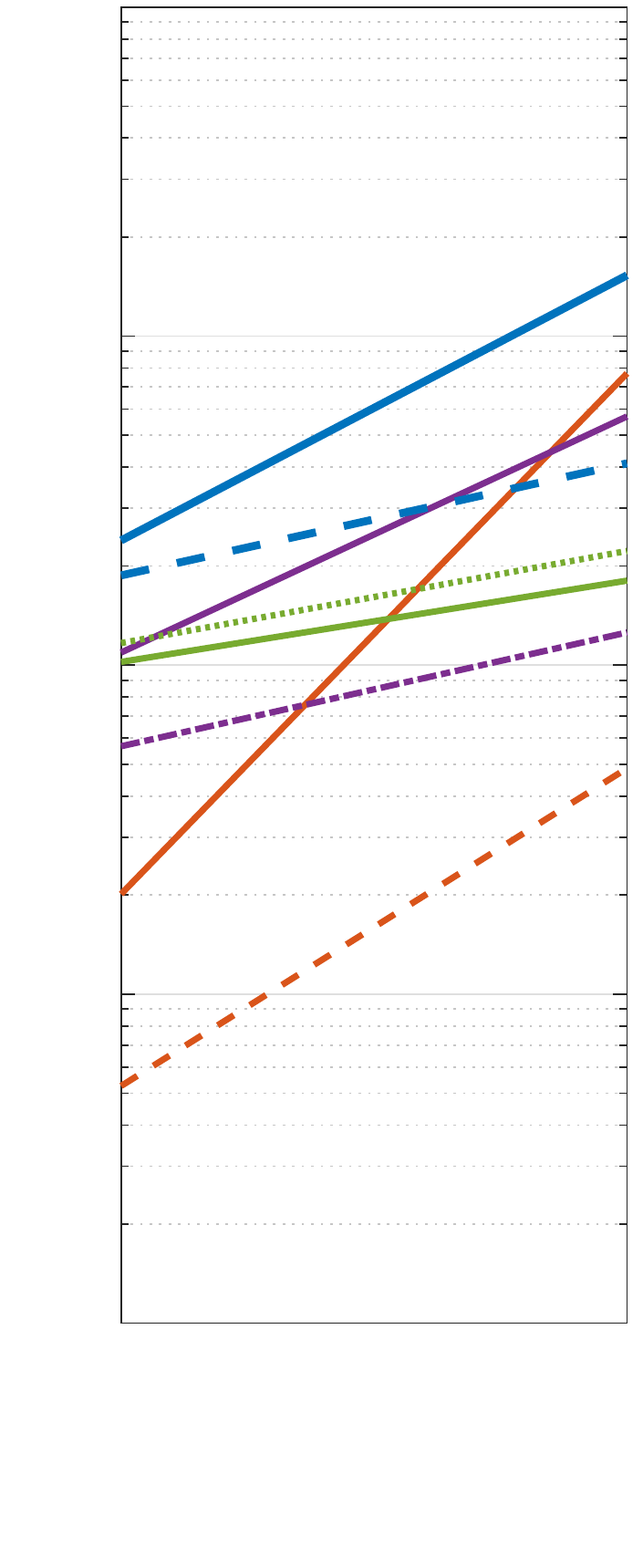}
		\put(9,6){\footnotesize No. of elements in}
		\put(10,2){\footnotesize each direction, $\ne$}
		\put(16,-6){\small (c) $\p=4$}
		\put(6.00,11.00){\footnotesize $2^{5}$}
		\put(37.00,11.00){\footnotesize $2^{6}$}
		\put(-0.50,15.00){\footnotesize $10^{-7}$}
		\put(-0.50,34.75){\footnotesize $10^{-6}$}
		\put(-0.50,55.50){\footnotesize $10^{-5}$}
		\put(-0.50,76.25){\footnotesize $10^{-4}$}
		\put(-0.50,97.00){\footnotesize $10^{-3}$}
	\end{overpic}\,\,\,
	\begin{overpic}[height=0.31\textheight,trim={0cm 0cm 0cm 0cm} ,clip]{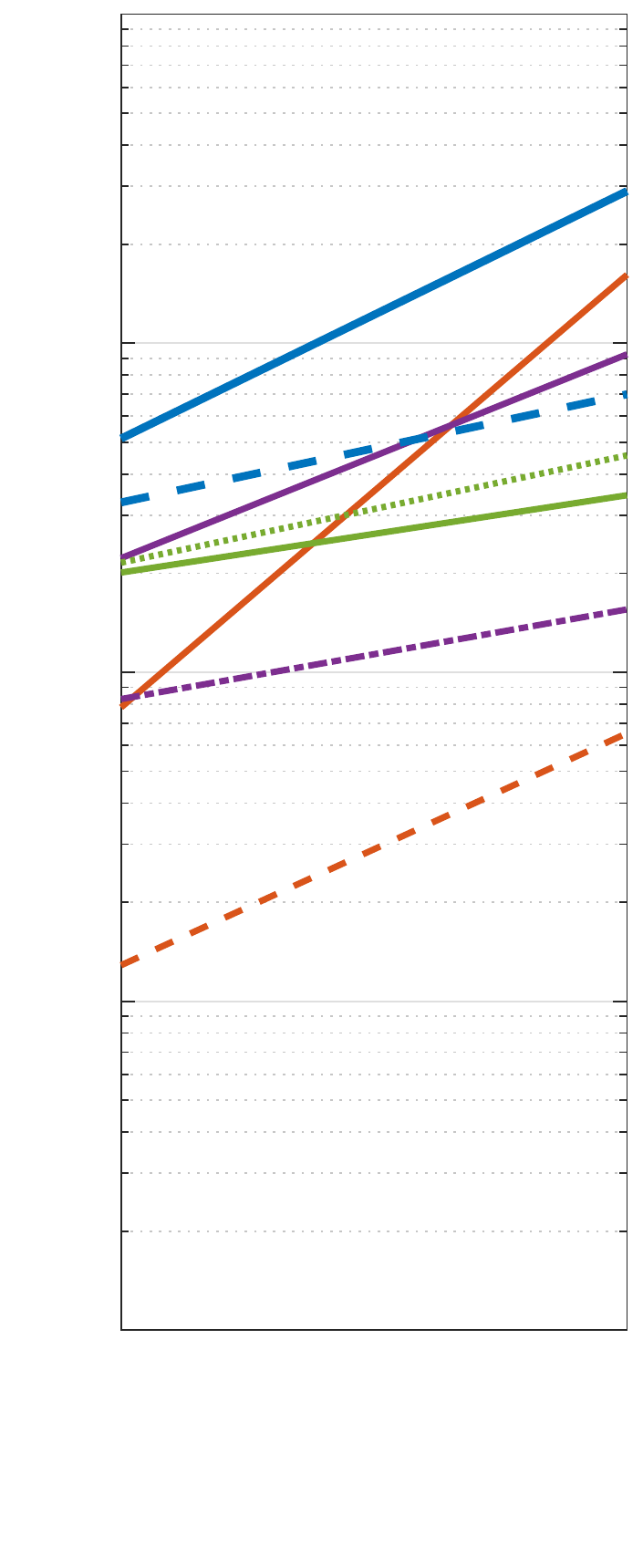}
		\put(9,6){\footnotesize No. of elements in}
		\put(10,2){\footnotesize each direction, $\ne$}
		\put(16,-6){\small (d) $\p=5$}
		\put(6.00,11.00){\footnotesize $2^{5}$}
		\put(37.00,11.00){\footnotesize $2^{6}$}
		\put(-0.50,15.00){\footnotesize $10^{-7}$}
		\put(-0.50,34.75){\footnotesize $10^{-6}$}
		\put(-0.50,55.50){\footnotesize $10^{-5}$}
		\put(-0.50,76.25){\footnotesize $10^{-4}$}
		\put(-0.50,97.00){\footnotesize $10^{-3}$}
	\end{overpic} \\
	\caption{Normalized average time per eigenvalue, ${\tih=\tiN/\N_0\.\N}$, versus the number of elements in each direction, $\ne$.
	3D eigenproblems under maximum-continuity IGA and optimal rIGA discretizations with $\BS$ of 16 elements in each direction
	(for degrees ${\p=4~\rm{and}~5}$, we only test up to ${\ne=64}$).}
	\label{fig.That3D}
\end{figure}

\figs{\ref{fig.T13D}}{\ref{fig.That3D}} show the average time per computed eigenvalue, $\tia$, and the normalized time, $\tih$.
As in the 2D case, 
the f/b elimination dominates the total cost in small problems.
Whereas for large problems, the matrix factorization is the most expensive procedure.
However,
to find a fixed number of eigenpairs, 3D problems employ more iterations than 2D ones (see \fig{\ref{fig.mconverge.b}}),
resulting in more \fb~eliminations and mat--vecs for each spectral transform (i.e., shift).
In our case,
the number of \fb~eliminations and mat--vec multiplications is close to ${\Nfb\approx 200~\Nfa}$ and ${\Nmv\approx 600~\Nfa}$ in 3D eigenproblems (see \tab{\ref{tab:T1_3D}}).


\section{Accuracy assessment of eigensolution using rIGA} 
\label{sec.Accuracy}

\begin{figure}[!h]
	\centering
	\begin{subfigure}{0.49\textwidth}\centering
		\includegraphics{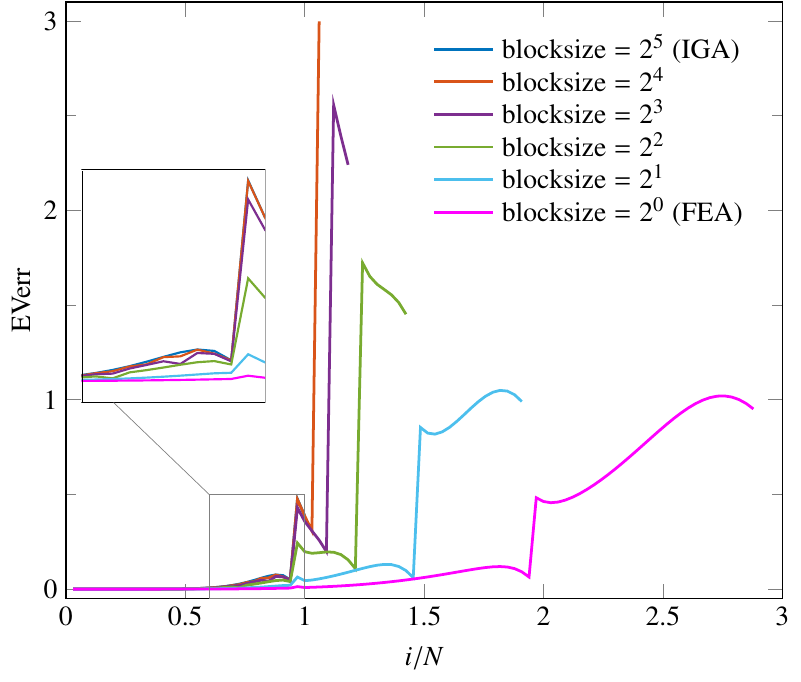}
		\caption{Eigenvalue error}
		\label{fig.EVerr1D.a}
	\end{subfigure}
	\begin{subfigure}{0.49\textwidth}\centering
		\includegraphics{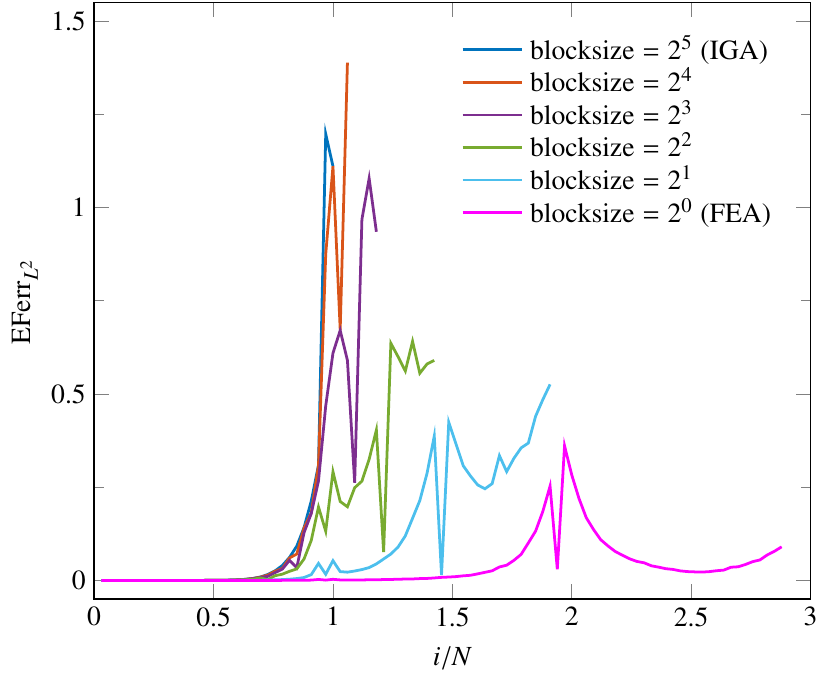}
		\caption{Eigenfunction $L^2$ error}
		\label{fig.EVerr1D.b}
	\end{subfigure}
	\caption{Eigenvalue and $L^2$ eigenfunction errors of the 1D eigenproblem discretized by ${\ne=32}$ with cubic bases: maximum-continuity IGA versus rIGA of different blocksizes.}
	\label{fig.EVerr1D}
\end{figure}

This section investigates the effect of employing rIGA discretizations on the accuracy of eigenanalysis
we utilize the eigenvalue error, $\EV$, and eigenfunction $L^2$ and energy norm errors, $\EFL$ and $\EFE$, respectively, as expressed by \eq{\eqref{eq.eigerrors}}.
The knot insertion steps of the rIGA approach add new control variables and, therefore, enriches the Galerkin space, modifying the spectral approximation properties of the IGA approach. 
In order to investigate how rIGA affects the accuracy of eigenpairs throughout the entire spectrum, 
we introduce a 1D eigenproblem with ${\ne=32}$ and ${\p=3}$.
\fig{\ref{fig.EVerr1D}} depicts the eigenvalue and $L^2$ eigenfunction errors of maximum-continuity IGA versus those obtained by different partitioning levels of the rIGA approach.
The abscissa of this figure shows the eigenmode numbers $i$ normalized with respect to $\N$, the total number of eigenmodes of the IGA discretization.
Since the rIGA-discretized system has more discrete eigenmodes, the spectra plots extend to ${i/\N > 1}$.
The eigenvalue errors are computed by comparing the approximated values, $\lmh$, with exact ones, ${\lme=i^{\.2}\pi^{\.2}}$, 
while the approximated eigenfunctions, $\uh$, are compared with ${\ue=\sqrt{2}\sin\pra{i\.\pi\x}}$.
In terms of eigenvalues (see \fig{\ref{fig.EVerr1D.a}}),
rIGA discretizations of lower $\BS$ reach a lower error at the same mode number
in the range of ${i/\N \leq 1}$,
improving the outliers behavior of the original IGA-discretized system.
It is justified in such a way that by decreasing the \BS, the eigenvalue errors converge to the
acoustic branch of the FEA spectrum and we achieved a better approximation (see, e.g.,~\cite{ Puzyrev2017}). 
However, for ${i/\N > 1}$ 
we observe larger errors for the outliers.
We obtain similar conclusions in terms of eigenvector errors (see \fig{\ref{fig.EVerr1D.b}}).

\begin{figure}[!h]
	\centering
	\begin{subfigure}{0.995\textwidth}\centering
		\includegraphics[width=0.325\textwidth]{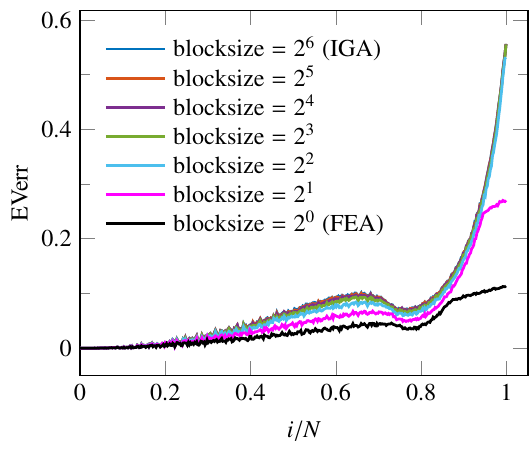}
		\includegraphics[width=0.325\textwidth]{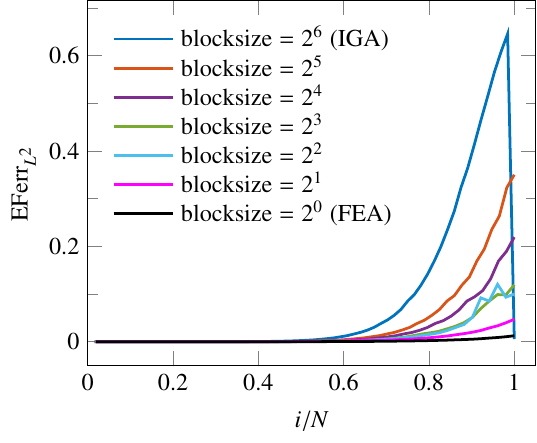}
		\includegraphics[width=0.325\textwidth]{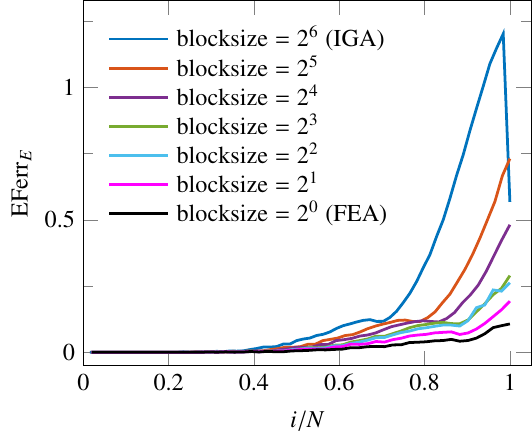}
		\caption{$\p=2$}
	\end{subfigure} \\[2pt]
	\begin{subfigure}{0.995\textwidth}\centering
		\includegraphics[width=0.325\textwidth]{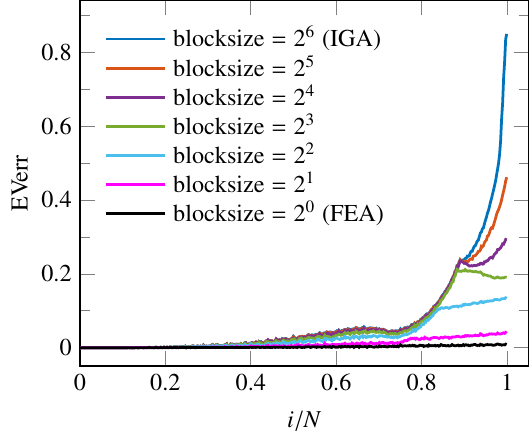}
		\includegraphics[width=0.325\textwidth]{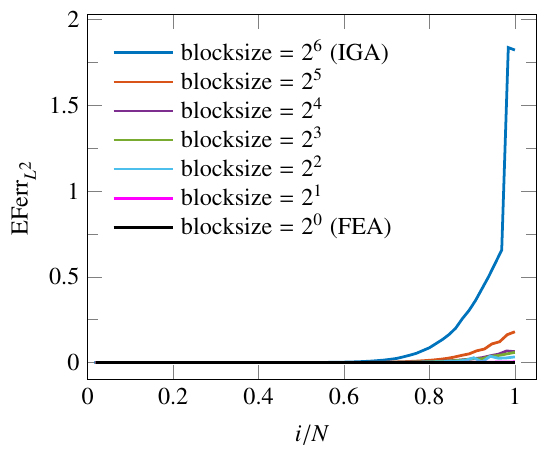}
		\includegraphics[width=0.325\textwidth]{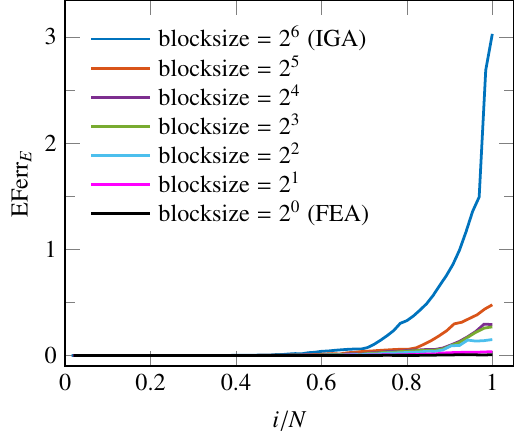}
		\caption{$\p=3$}
	\end{subfigure}	\\[2pt]
	\begin{subfigure}{0.995\textwidth}\centering
		\includegraphics[width=0.325\textwidth]{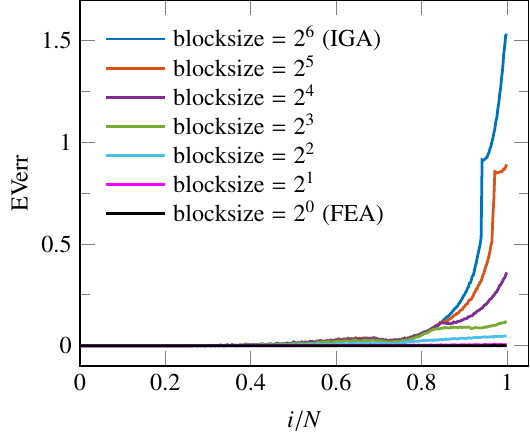}
		\includegraphics[width=0.325\textwidth]{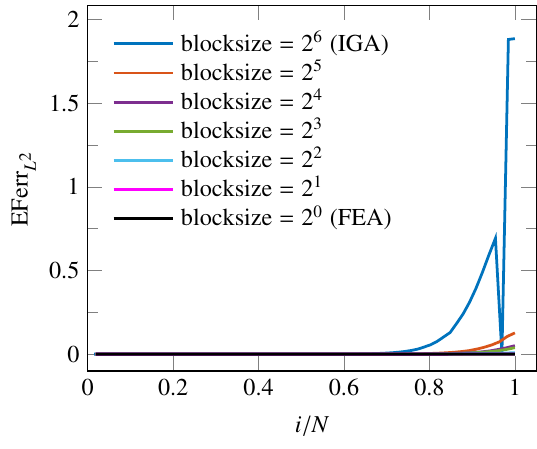}
		\includegraphics[width=0.325\textwidth]{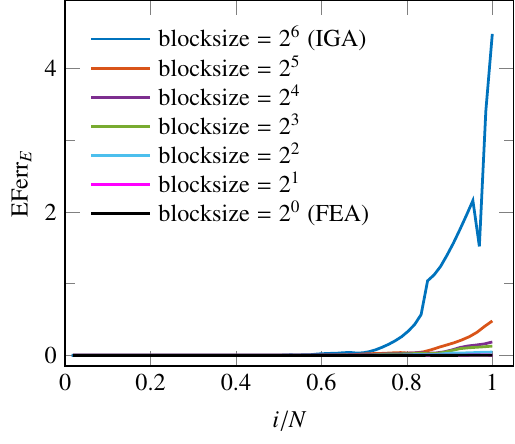}
		\caption{$\p=4$}
	\end{subfigure}	\\[2pt]
	\begin{subfigure}{0.995\textwidth}\centering
		\includegraphics[width=0.325\textwidth]{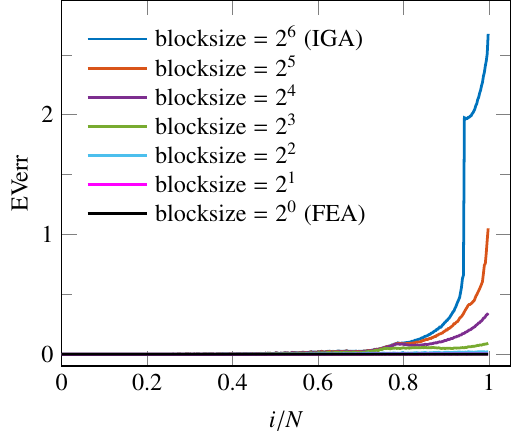}
		\includegraphics[width=0.325\textwidth]{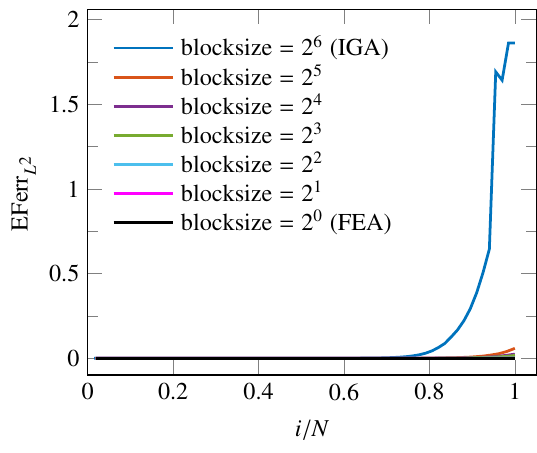}
		\includegraphics[width=0.325\textwidth]{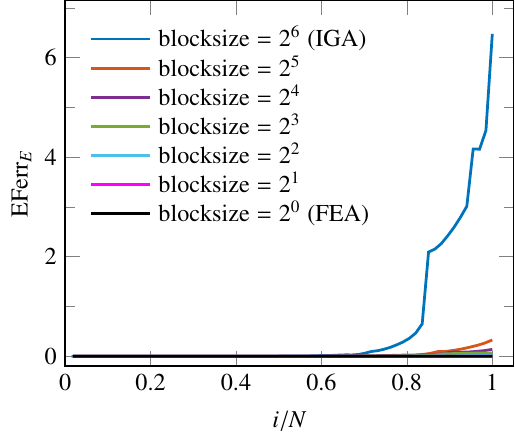}
		\caption{$\p=5$}
	\end{subfigure}	
	\caption{Eigenvalue error, $\EV$, and eigenfunction $L^2$ and energy norm errors, $\EFL$ and $\EFE$, of the 2D eigenproblem discretized by ${\ne=64}$ elements in each direction and polynomial degrees of ${\p=2,3,4,5}$, obtained by maximum-continuity IGA and rIGA of different \BS s.
	}
	\label{fig.EVerr64}
\end{figure}
\begin{figure}[!h]
	\centering
	\begin{subfigure}{0.99\textwidth}\centering
		\includegraphics[width=0.325\textwidth]{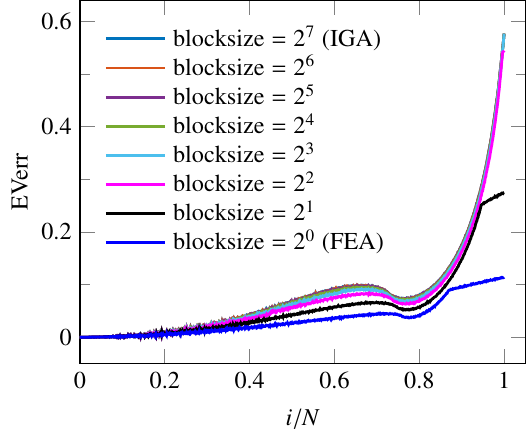}
		\includegraphics[width=0.325\textwidth]{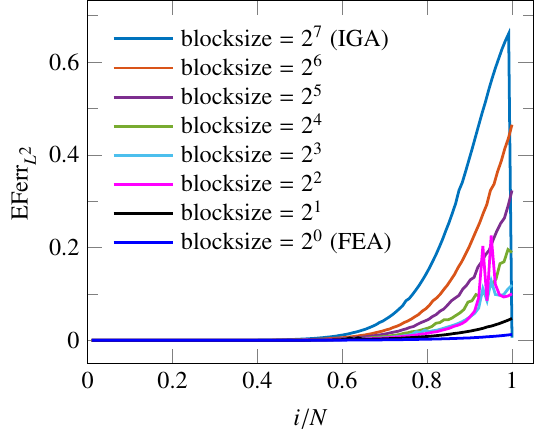}
		\includegraphics[width=0.325\textwidth]{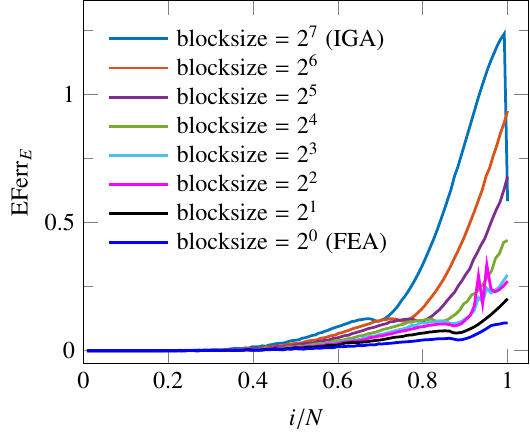}
		\caption{$\p=2$}
	\end{subfigure} \\[3pt]
	\begin{subfigure}{0.99\textwidth}\centering
		\includegraphics[width=0.325\textwidth]{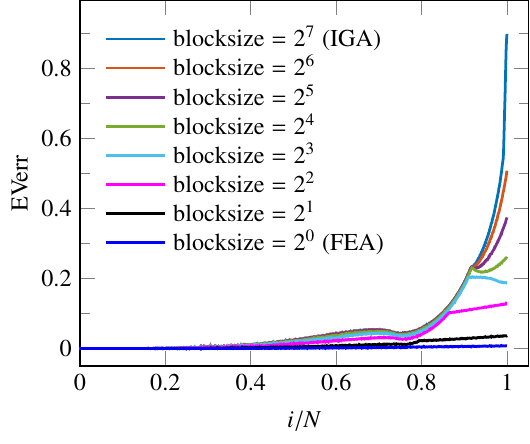}
		\includegraphics[width=0.325\textwidth]{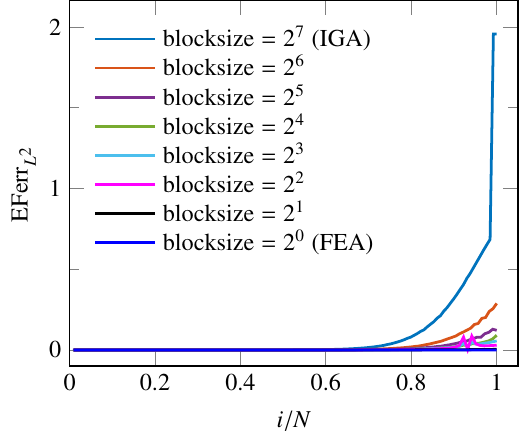}
		\includegraphics[width=0.325\textwidth]{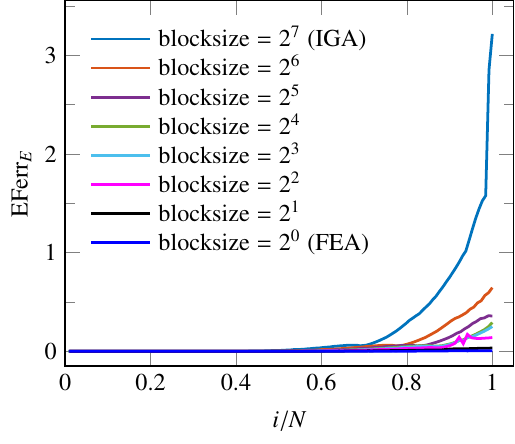}
		\caption{$\p=3$}
	\end{subfigure} \\[4pt]	
	\begin{subfigure}{0.99\textwidth}\centering
		\includegraphics[width=0.325\textwidth]{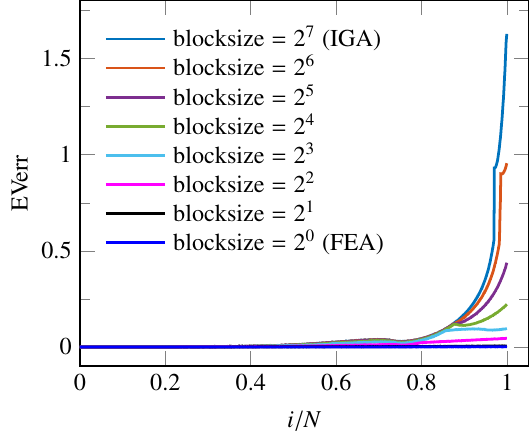}
		\includegraphics[width=0.325\textwidth]{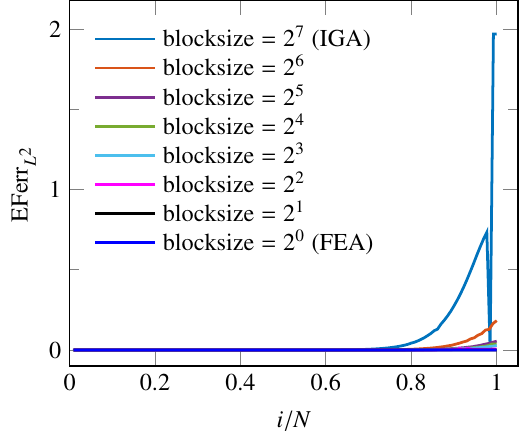}
		\includegraphics[width=0.325\textwidth]{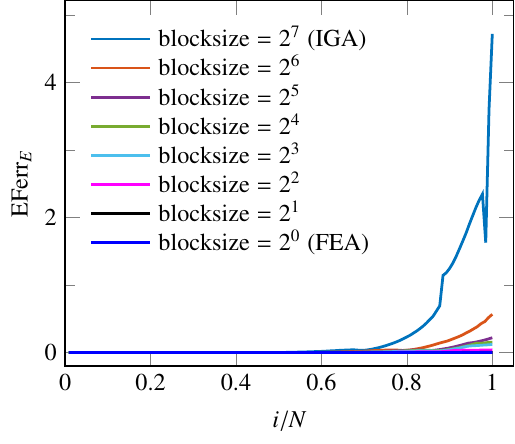}
		\caption{$\p=4$}
	\end{subfigure} \\[4pt]	
	\begin{subfigure}{0.99\textwidth}\centering
		\includegraphics[width=0.325\textwidth]{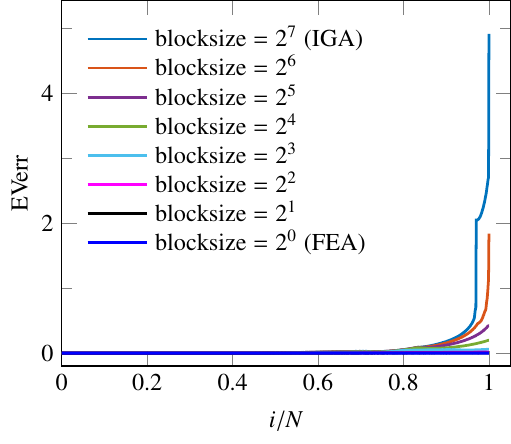}
		\includegraphics[width=0.325\textwidth]{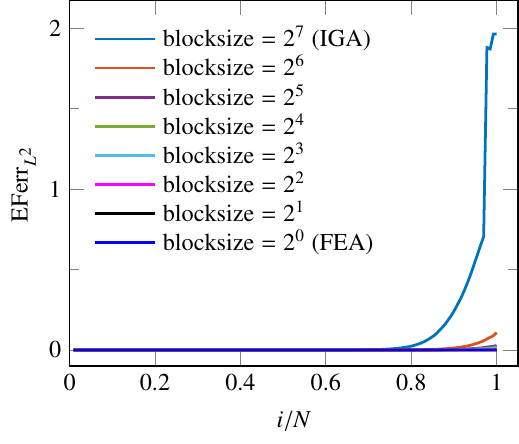}
		\includegraphics[width=0.325\textwidth]{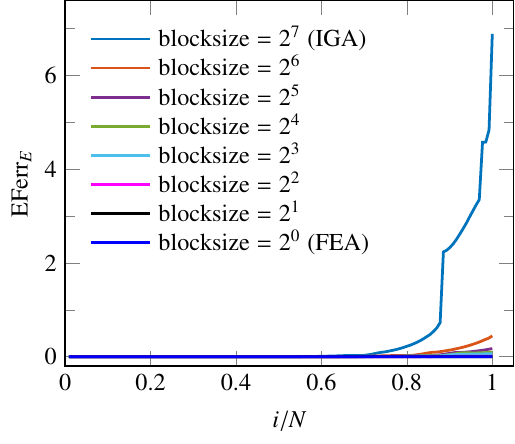}
		\caption{$\p=5$}
	\end{subfigure}	
	\caption{Eigenvalue error, $\EV$, and eigenfunction $L^2$ and energy norm errors, $\EFL$ and $\EFE$, of the 2D eigenproblem discretized by ${\ne=128}$ elements in each direction and polynomial degrees of ${\p=2,3,4,5}$, obtained by maximum-continuity IGA and rIGA of different \BS s.
	}
	\label{fig.EVerr128}
\end{figure}

We now consider the accuracy of the eigenanalysis for 2D systems discretized with ${\ne=64}$ and 128 elements in each direction and polynomial degrees of ${\p=2,3,4,5}$.
The exact eigenvalue and eigenfunctions are expressed by \eq{\eqref{eq.exact}}.
We use \Rem{\ref{rem.matvec}} to compare the approximate eigenpairs $\lmh_i$ and $\uh_i$ with the tensor-represented exact ones $\lme_{ij}$ and $\ue_{ij}$\..
Note that
for any off-diagonal combination of $i$ and $j$ (i.e., ${i\neq j}$),
\mbox{$\ue_{ij}$ and $\ue_{ji}$} are generally referred to as
\textit{degenerate} eigenfunctions (see, e.g.,~\cite{ Korsch1983}).
Since degenerate eigenfunctions are not unique,
we restrict our eigenfunction accuracy assessments only to diagonal eigenfunctions.
\figs{\ref{fig.EVerr64}}{\ref{fig.EVerr128}} describe the eigenvalue error as well as eigenfunction $L^2$ and energy norm errors of the above-mentioned systems.
In both cases, we find the first ${\N_0=\N_{\rm{IGA}}}$ eigenpairs of the problem.
As for the 1D systems, rIGA discretizations improve the accuracy of eigenpairs approximation compared to the maximum-continuity IGA when ${i/\N\leq 1}$. 
(although we may see higher errors in the outliers for ${i/\N_0>1}$). 
These error differences are more noticeable for higher $\p$.


\section{Conclusions} 
\label{sec.Conclusions}

This paper proposes the use of refined isogeometric analysis (rIGA) discretizations to solve generalized Hermitian eigenproblems (GHEP). 
We compare the computational time of rIGA versus that of maximum-continuity IGA when computing a fixed number of eigenpairs using a Lanczos eigensolver with a shift-and-invert spectral transform approach.

We consider two cases attending to the problem size. 
For large problems, namely, ${\N\geq256^2}$ in 2D and ${\N\geq64^3}$ in 3D, 
the most expensive operation is the matrix factorization. 
This is followed by matrix--vector operations. 
We compute the Cholesky factorization ${\O(\p^{\.2})}$ times faster with an rIGA discretization than with an IGA one.
As a result, in the asymptotic regime, we theoretically reach an improvement of ${\O(\p^{\.2})}$ in the total computational time of the eigenanalysis.
In our computations with $\N$ up to ${2048^2}$ in 2D, and ${128^3}$ in 3D, savings associated with the implementation of rIGA limit to ${\O(\p)}$. 
This occurs because we need to consider larger matrices to arrive at the asymptotic limit. 
In a 2D mesh with ${2048^2}$ elements and ${\p=5}$, we need an average of 141 seconds to compute each eigenpair when using IGA, and only 26 seconds for rIGA.
In a 3D mesh with ${128^3}$ elements and ${\p=3}$, rIGA reduces the average required time per eigenmode from 590 seconds to 112 seconds.

For smaller problems (and with lower polynomial degrees), the forward/backward elimination is the most expensive numerical operation.
This operation,
which is called almost 150 and 200 times per each matrix factorization in our 2D and 3D cases,
is theoretically up to ${\O(\p)}$ faster with rIGA than with IGA for sufficiently large problems.
For small problems, our improvement hardly reaches the expected rates.
As a result, 
we suggest to use the maximum-continuity IGA discretization for small problems.

Finally, 
we obtain a better accuracy in spectral approximation using rIGA discretizations
when computing the first $\N_0$ eigenpairs,
being $\N_0$ the size of the IGA-discretized system.
This occurs because the continuity reduction of basis functions enriches the Galerkin space and modifies the approximation properties of the IGA approach.
However, rIGA eigenvalues beyond $\N_0$ show important inaccuracies. This may be detrimental for solving nonlinear problems.


\section*{Acknowledgment}

This work has received funding from 
the European Union's Horizon 2020 research and innovation program under the Marie Sklodowska-Curie grant agreement No 777778 (MATHROCKS), 
the European POCTEFA 2014-2020 Project PIXIL (EFA362/19) by the European Regional Development Fund (ERDF) through the Interreg V-A Spain-France-Andorra program, 
the Project of the Spanish Ministry of Science and Innovation with reference PID2019-108111RB-I00 (FEDER/AEI),
the BCAM ``Severo Ochoa" accreditation of excellence (SEV-2017-0718) and the Basque Government through the BERC 2018-2021 program, 
the two Elkartek projects ArgIA (KK-2019-00068) and MATHEO (KK-2019-00085),
the grant ``Artificial Intelligence in BCAM number EXP. 2019/00432", 
and the Consolidated Research Group MATHMODE (IT1294-19) given by the Department of Education.
The authors also acknowledge the Group of Applied Mathematical Modeling, Statistics, and Optimization \mbox{(MATHMODE)} at the 
University of the Basque Country \mbox{(UPV/EHU)}
for providing HPC resources that have contributed to the research results reported in the paper.


\section*{References}
\small
\bibliography{_eig_refs}

\end{document}